\documentclass[a4paper,11pt]{amsart}
\usepackage[LGR,X2,T2C,T1]{fontenc}
\usepackage{graphicx}
\usepackage{subcaption}
\usepackage{float}
\usepackage{amsmath,amsthm,amssymb}
\usepackage{comment}
\usepackage{url}
\newcommand{\mychoice}[3]{#1
}
\mychoice{
\newcommand{\plabel}[1]{ \label{#1}}
\newcommand{\gbibitem}[1]{ \bibitem{#1}}
\newcommand{\snewpage}{}

\specialcomment{commentx}{}{}\excludecomment{commentx}
\specialcomment{commenty}{}{}\excludecomment{commenty}
\specialcomment{commentz}{}{}

}{
\usepackage{xcolor}
\newcommand{\plabel}[1]{ \label{#1}\rlap{\smash{${}^{^{[#1]}}$}}}
\newcommand{\gbibitem}[1]{ \bibitem{#1}\rlap{\smash{${}^{^{[#1]}}$}}}
\newcommand{\snewpage}{\newpage}

\newenvironment{commentx}{\color{magenta} }{\color{black} }
\newenvironment{commenty}{\color{blue} }{\color{black} }

}{
\usepackage{xcolor}
\newcommand{\plabel}[1]{ \label{#1}}
\newcommand{\gbibitem}[1]{ \bibitem{#1}}
\newcommand{\snewpage}{}

\specialcomment{commentx}{}{}\excludecomment{commentx}

\specialcomment{commentz}{}{}\excludecomment{commentz}

}
\DeclareMathOperator{\sgn}{sgn}
\DeclareMathOperator{\dist}{dist}
\DeclareMathOperator{\artanh}{artanh}
\DeclareMathOperator{\arcosh}{arcosh}
\DeclareMathOperator{\arsinh}{arsinh}
\DeclareMathOperator{\arcoth}{arcoth}
\DeclareMathOperator{\Len}{\mathbf L}
\DeclareMathOperator{\Area}{\mathbf A}
\newcommand{\ppp}{{}^{+}\!}
\newcommand{\nnn}{{}^{-}\!}
\newcommand{\hyp}{\mathrm{hyp}}
\newcommand{\SdS}{\mathrm{SdS}}
 \newcommand{\leaveout}[1]{}
\newcommand{\marginextend}[1]{ \addtolength{\oddsidemargin}{-#1}  \addtolength{\evensidemargin}{-#1}
  \addtolength{\textwidth}{#1}\addtolength{\textwidth}{#1}}
\newcommand{\updownextend}[1]{ \addtolength{\topmargin}{-#1}  \addtolength{\textheight}{#1}
\addtolength{\textheight}{#1}}
\marginextend{1cm}
\updownextend{0cm}
\allowdisplaybreaks[4]
\title{Hyperbolic elliptic parabolic disks approximated by half distance bands}
\author{Gyula Lakos}
\email{gyula.lakos@uni-miskolc.hu}
\address{Institute of Mathematics, Department of Analysis, University of Miskolc, H-3515 Miskolc-Egyetemváros, Hungary}
\begin{document}
\begin{abstract}
Hyperbolic elliptic parabolic disks can be described by the inequality $\frac{x^2}{C^2}+2y^2-2y\leq0$  ($0<C<1$) in the unit disk
 based Beltrami--Cayley--Klein model of the hyperbolic geometry, up to hyperbolic congruences.
The hyperbolic elliptic parabolic disks considered above are sort of close  to their supporting half distance bands given
 by the inequalities $\frac{x^2}{C^2}+ y^2-1\leq0$ and $y\geq0$.
Here we consider what `close' might mean, and we look for even more precise approximations, in terms of area and circumference.
\end{abstract}
\maketitle
\snewpage
\section*{Introduction}
It is said that, in geometry, everything important comes in three versions: elliptic, parabolic, and hyperbolic.
It is therefore particularly pleasant that some objects mentioned in the title of this paper possess  all these three adjectives.
Speaking more seriously, however, in hyperbolic geometry, elliptic parabolas are a kind of  hyperbolic  quadrics,
 and whose convex closures are the hyperbolic elliptic parabolic disks in question.

\begin{figure}[htp]
   \begin{subfigure}[b]{2.5in}
    \includegraphics[width=2.5in]{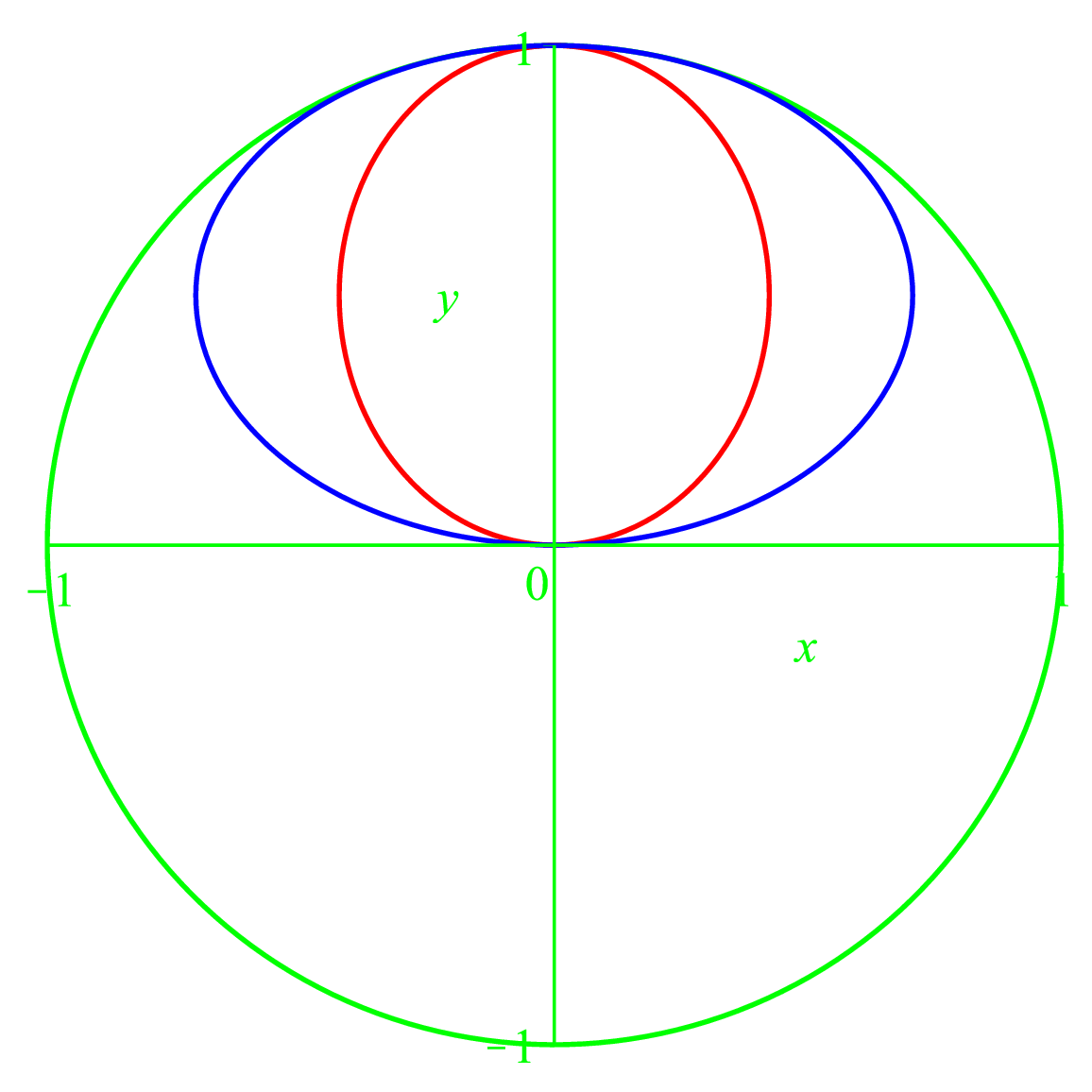}
    \caption*{\ref{fig:figHEP01}(a)  Supporting horocycle [blue].}
  \end{subfigure}
  \begin{subfigure}[b]{2.5in}
    \includegraphics[width=2.5in]{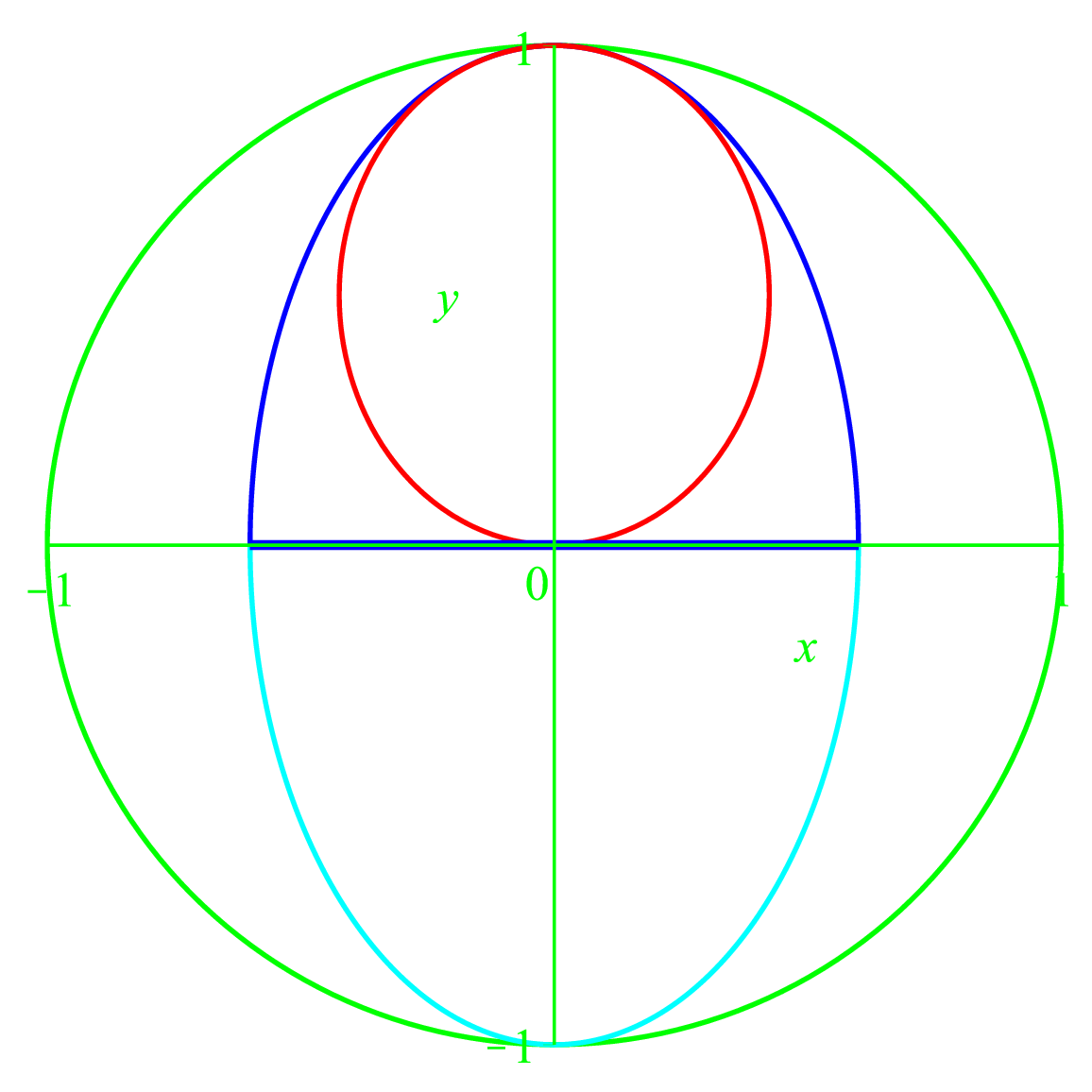}
    \caption*{\ref{fig:figHEP01}(b)  Supporting half distance band [blue].}
  \end{subfigure}
  \\

   \begin{subfigure}[b]{2.5in}
    \includegraphics[width=2.5in]{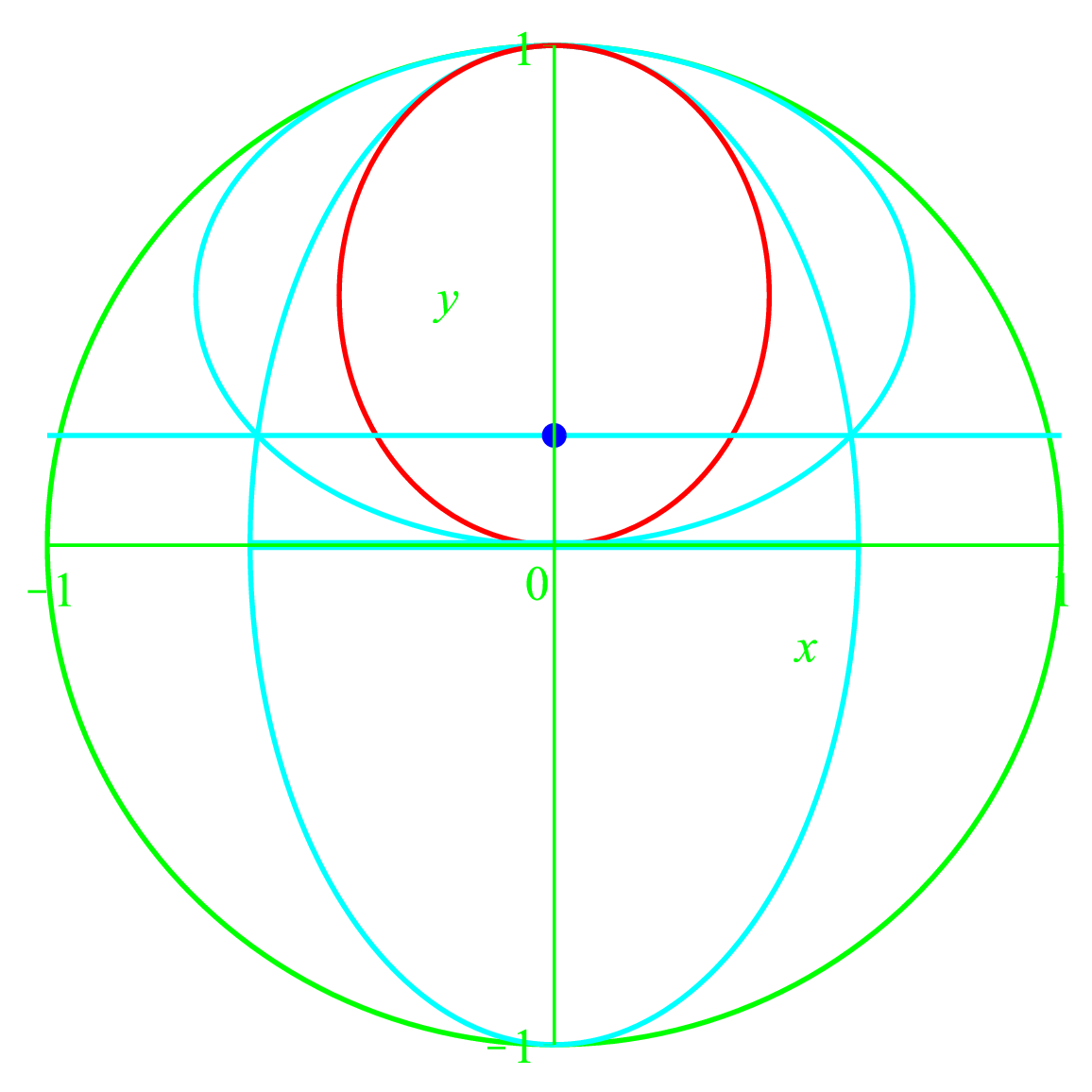}
    \caption*{ \ref{fig:figHEP01}(c)  The focus I [blue].}
  \end{subfigure}
  \begin{subfigure}[b]{2.5in}
    \includegraphics[width=2.5in]{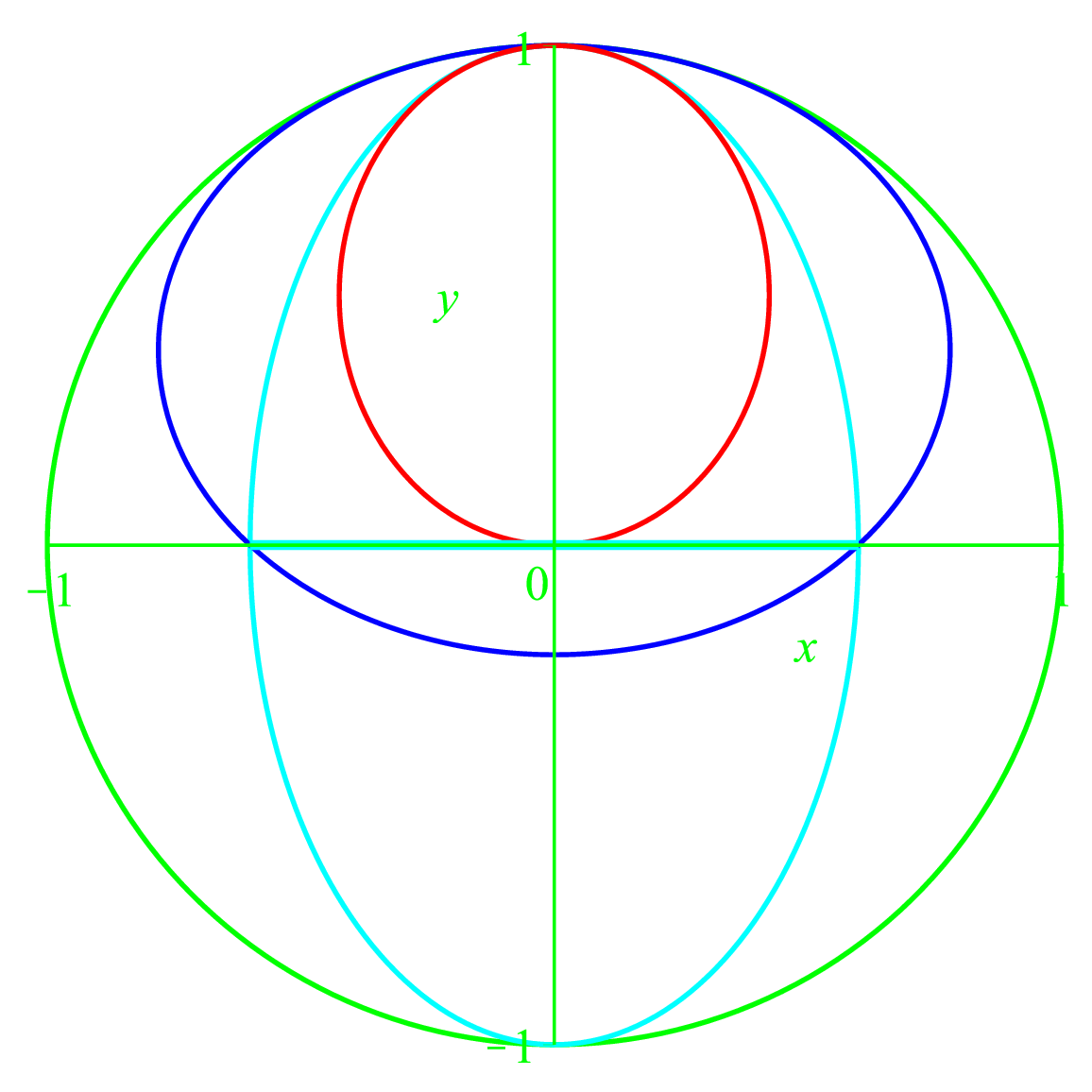}
    \caption*{\ref{fig:figHEP01}(d)  Directrix horocycle  [blue].}
  \end{subfigure}
   \\

   \begin{subfigure}[b]{2.5in}
    \includegraphics[width=2.5in]{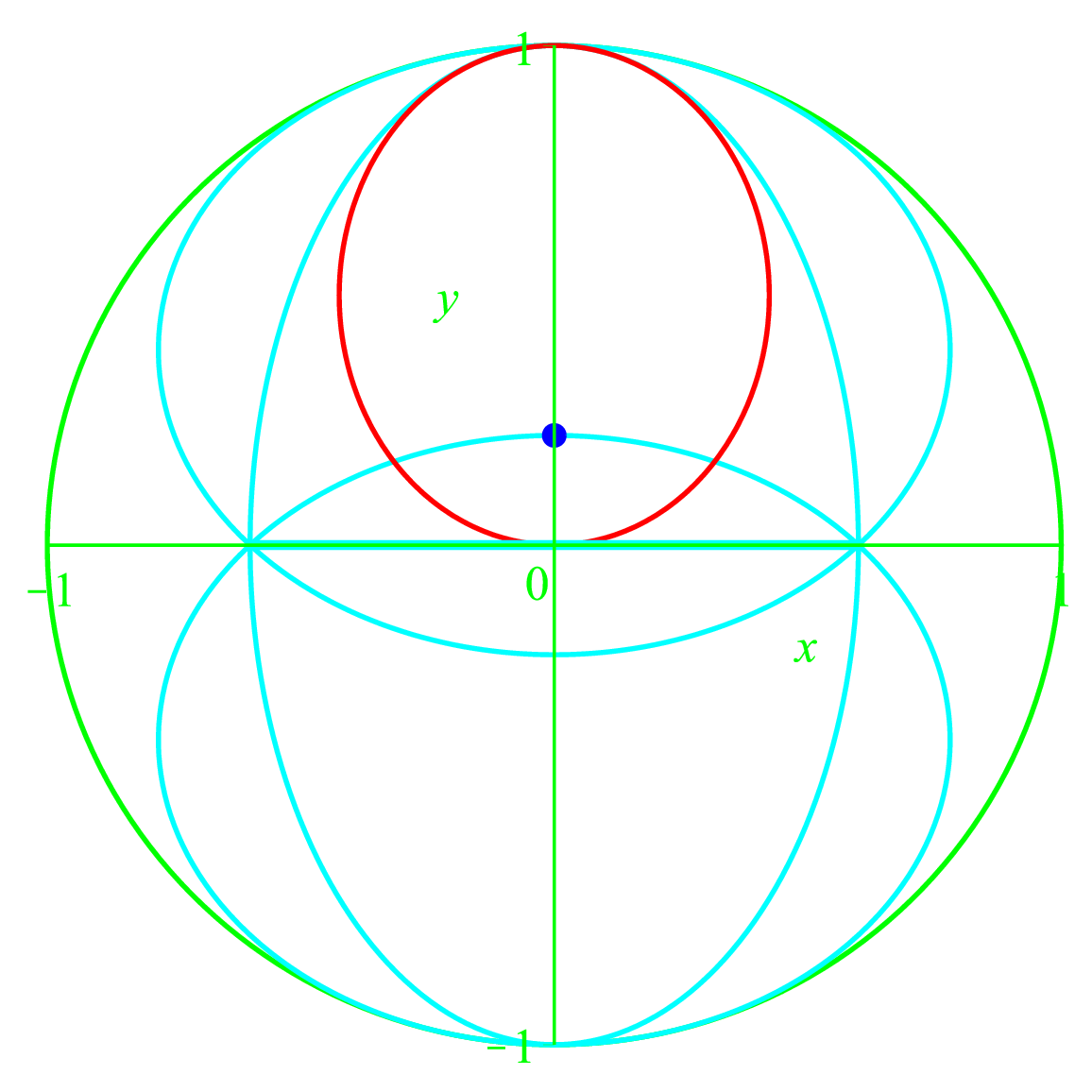}
    \caption*{ \ref{fig:figHEP01}(e)  The focus II   [blue].}
  \end{subfigure}
  \begin{subfigure}[b]{2.5in}
    \includegraphics[width=2.5in]{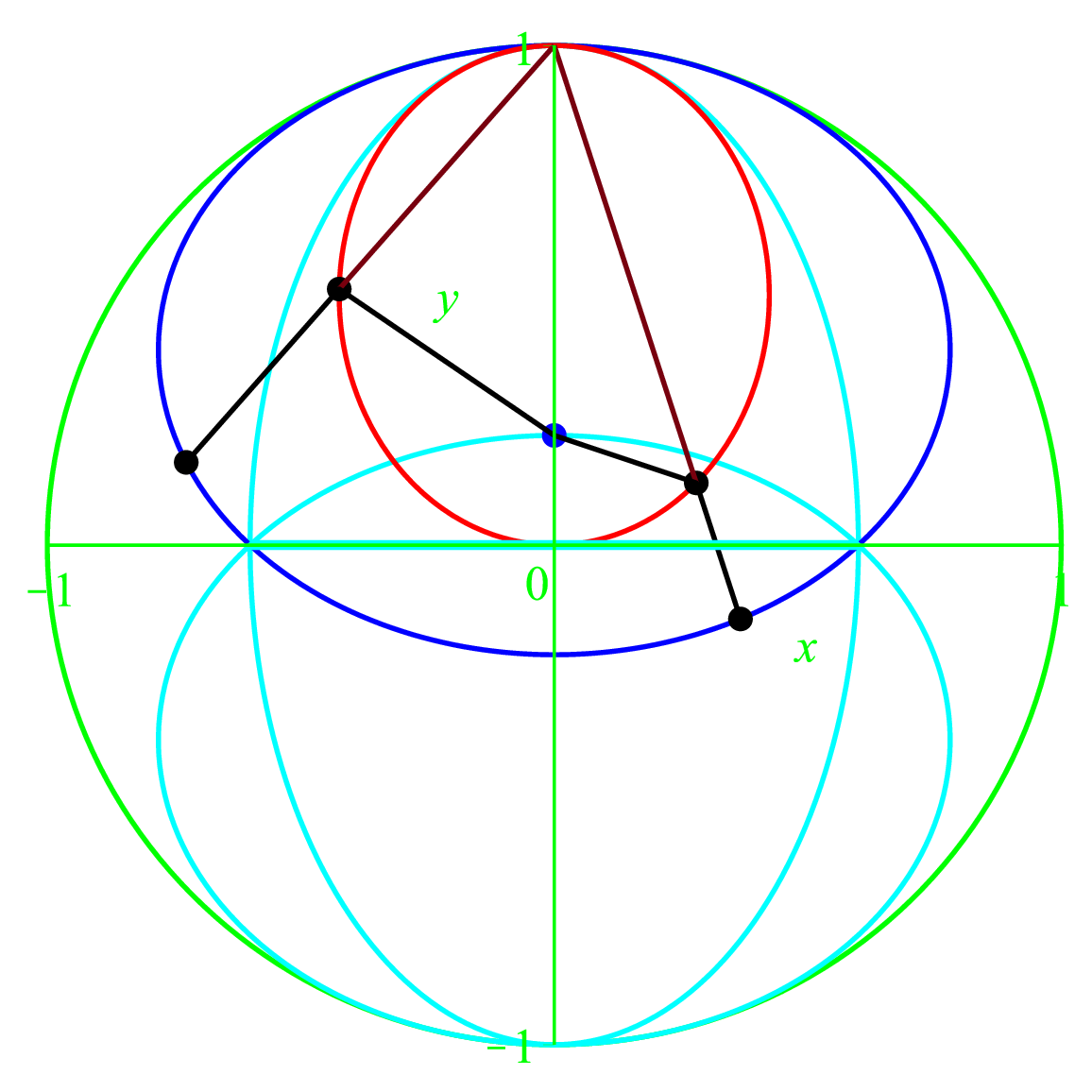}
    \caption*{\ref{fig:figHEP01}(f)  The focal distance property [blacks].}
  \end{subfigure}
   \caption*{Fig.~\ref{fig:figHEP01} Geometric elements related the $h$-elliptic parabola [red].}
\phantomcaption
\plabel{fig:figHEP01}
\end{figure}

Hyperbolic elliptic parabolic disks can be described by the inequality
\begin{equation}
\frac{x^2}{C^2}+2y^2-2y\leq0
\plabel{eq:amat1}
\end{equation}
($0<C<1$) in the unit disk
 based Beltrami--Cayley--Klein model of the hyperbolic geometry, up to hyperbolic congruences.
(In this canonical presentation, it looks like the supporting horodisk squeezed by ratio $C$, see Figure \ref{fig:figHEP01}(a).)
The above mentioned hyperbolic elliptic parabolic disks are quite close, for example,  to their supporting half distance bands given
 by the inequalities
 \begin{equation}
 \frac{x^2}{C^2}+ y^2-1\leq0\qquad\text{ and }\qquad y\geq0.
 \plabel{eq:amat2}
\end{equation}
(In this canonical presentation the supporting distance band looks like the model disk squezed by ratio $C$, and
 the supporting half-plane is produced by the tangent at the vertex of $h$-ellipic parabola, see Figure \ref{fig:figHEP01}(b).)
This phenomenon is very hard to miss if one considers hyperbolic elliptic parabolas even just very briefly.
(Especially if it is paired with a synthetic viewpoint.)
Both the hyperbolic elliptic parabolic disks and their supporting half distance bands are convex
 (the Beltrami--Cayley--Klein model is faithful in this respect), and the former ones are contained in the latter ones
 (and their boundaries are seen in $360^\circ$ from the interior points of the former ones), thus the ``approximation'' must be
 relatively crude both in terms of area and circumference:
hyperbolic elliptic parabolic disks should have a ``lesser amount'' (even if infinity) of those than the supporting half distance bands,
 but comparable (even if infinity).
In this paper we are going to look after this situation and make these statements more precise; just for fun.

Now, the focus point of the above considered hyperbolic elliptic parabolas is given by the point $\left(0,\frac{C^2}{2-C^2}\right)$.
(See Figure \ref{fig:figHEP01}(c).)
Therefore an educated guess would be to approximate the hyperbolic elliptic parabolic disks by the half distance bands given
 by the inequalities
\[\frac{x^2}{C^2}+ y^2-1\leq0\qquad\text{ and }\qquad y\geq\frac{C^2}{2-C^2}.\]
This does not work out particularly well either, in general; but in the process we can, after all, find replacement objects
 for the hyperbolic elliptic parabolic disks in area or circumference.

This paper is a kind of propaganda for hyperbolic conics.
It also conveys a picture about (elementary) hyperbolic geometry.
However, it is not an introduction to hyperbolic geometry.
It is written for one who has already developed an affinity toward hyperbolic geometry,
 or for a student who has taken a course in hyperbolic geometry, and wants to check whether
 his or her knowledge is comprehensive enough to deal with a rather minor problem in hyperbolic geometry.

\snewpage
\section{Background --- hyperbolic geometry}\plabel{sec:Bg}
\textbf{Hyperbolic geometry.}
After many previous thoughts, most notably, by Euclid, Saccheri, Lambert, Schweikart, Taurinus, Gauss,
 see Stäckel, Engel \cite{SE}, hyperbolic geometry
 enters to the scene by the independent work of Lobachevsky \cite{Lob0} (1829) and Bolyai \cite{Bolyai0} (1831), see Bonola \cite{Bon}.
The original editions are not easy to access (in particular, for being in Russian and in Latin, respectively),
 but there are  commented translations which can be read very profitably.
Regarding Lobachevsky's  original text \textit{Elements} (or \textit{Principles)} \cite{Lob0},
 one can refer to the German translation by Engel \cite{Lob1}.
Moreover, Lobachevsky has also published follow-up works (of various depth), some of which have appeared in French or German
 relatively soon, see \cite{Lob01} (1837), \cite{Lob02} (1840) (cf.~\cite{Lob03}).
For us, the relatively recent English  edition of Lobachevsky's \textit{Pangeometry}
 (Russian version: 1855, French version: 1856) by Papadopoulos \cite{Lob2} can be the most relevant.
This is originally Lobachevsky's last and probably most polished work on hyperbolic geometry, yet it is quite close to his original one.
I find the well-commented edition by Papadopoulos  extremely useful; especially because it also offers a
 good guide to the existing literature on classical hyperbolic geometry, including other editions of the basic works.
%Indeed, for those references, the reader directed to \cite{Lob2}.
(There are some occasional typos due to polyglotte environment, though.)
Regarding Bolyai's original text \textit{Appendix} \cite{Bolyai0}, one has the nice English translation by Halsted \cite{Bolyai1}, but
 the well-commented edition by Kárteszi \cite{Bolyai4} is very instructive.
(In general, for the sake accessibility, we will make references in English primarily; if it is not possible, then in German or French;
 then in Italian and Russian; and, at last, in any other languages, language indicated.
We do so even if in the future language barriers will predictably be eased further by technology.)
The works of Lobachevsky and Bolyai are of synthetic, but only of semi-axiomatic nature.
This is, of course, not surprising as no precise axiomatics have been developed at that time yet;
 the original (1899) edition of Hilbert's \textit{Grundlagen} \cite{Hil1} being published only much later.
However, both Lobachevsky and Bolyai  were displeased by the lack a proof of consistency for hyperbolic
 geometry relative to Euclidean geometry.

The demonstration of consistency was achieved by  models of the hyperbolic geometry.
(Here we note that regarding the genesis of the very first models, Stillwell \cite{Sti} provides not only good
 translations but insightful comments.)
The disc based Beltrami--Cayley--Klein model was exhibited by Beltrami \cite{Bel1} (1868).
Of course, in its viewpoint, Beltrami's work tilts toward differential geometry.
(His famous pseudosphere also appears there.)
Klein \cite{Kle1} (1871) (who is aware of Beltrami \cite{Bel1}, \cite{Bel2})
 reworks the disc based Beltrami--Cayley--Klein model to purely projective terms.
His work was inspired by Cayley \cite{Cay}
 (hence the name appears; in those times the unity of the geometries with constant curvature was a major point).
Klein's work strongly tilts toward projective geometry, he also considers more general base quadrics for the projective model.
In any case, both Beltrami and Klein have the distance formula \eqref{eq:kdist}.
In a continuation, Beltrami \cite{Bel2} (1868) proceeds to study the conformal ball and half-space models,
 earlier and in higher dimensions than Poincaré \cite{Poi2}.
In a continuation, Klein \cite{Kle2} (1873) elaborates further on certain points related to projective geometry.
The hyperboloid model appears in Killing \cite{Kil1} (1879) tangentially, then, in greater detail, in Killing \cite{Kil2} (1880)
 (cf.~also Killing \cite{Kil3}).
Killing attributes this general way of coordinatization to Weierstrass from 1872 (which seems certainly be so in the case of positive curvature).
This is somewhat before Poincaré \cite{Poi1} and Cox \cite{Cox}.
Although these works are also still well before the reaxiomatization by Hilbert, they leave little doubt about
 the relative consistency of hyperbolic geometry.

\snewpage
It must be said, however, that Lobachevsky and Bolyai already achieve a sort of
 computational universality in hyperbolic geometry.
In particular, they use a kind of ``rectangular coordinate system'' with coordinates $(p,q)\in\mathbb R^2$
 such that $q$ is the oriented distance of the point from a fixed distance-gauged line, and $p$ is
 the gauge value of the perpendicular projection of the point to the fixed line.
There they find that the infinitesimal line element in that coordinate system is
$(\mathrm ds)^2=(\cosh\frac qk)^2(\mathrm dp)^2+(\mathrm dq)^2$.
Cf.~\cite{Lob1}, \S21.~(34) (note: $\cosh\frac qk=\frac1{\sin\Pi(q)}$) and \cite{Bolyai4} \S32.~(I--II).
(Here $k>0$ is the scaling constant of the metric relative to a distingushed one, see more about that later later.)
In particular, in modern viewpoint, they construct the hyperbolic plane as a Riemannian space.
Except, of course, at that time Riemannian geometry was not yet developed.
Riemann delivered his famous lecture, itself partly inspired by hyperbolic geometry, only in 1854, appearing
 in print only in 1868  \cite{Rie1}, cf.~\cite{Rie2} (extensively commented in Spivak \cite{Spi}, Jost \cite{Rie3}).
(It is notable that how much Riemann's work, ultimately, inspired Beltrami,
 see  \cite {Lob2} and also Alcalá Vicente \cite{Alc}.)
Another matter is, however, that the ``rectangular coordinate system'' is not particularly convenient to study the congruences of the
 hyperbolic plane.

Comprehensive introductions to hyperbolic geometry are Coolidge \cite{Coo} (1909),
 Som\-mer\-ville \cite{Som} (1914), Liebmann \cite{Lib1} (1905), \cite{Lib2} (1912), \cite{Lib3} (1923), and the posthumous book of Klein \cite{Kle3} (1928) based on  \cite{Kle31}, \cite{Kle32} (1893) but later reworked so that even Hilbert \cite{Hil2} (1901) is referenced.
These books may  look like   dated but this is not quite so:
In general, hyperbolic geometry in itself was in the relative forefront of mathematical attention
 only before World War I; then attention turned toward more general geometries,
 not in the least due to the work of people like Riemann, Beltrami, Klein, Killing, Engel, Friedrich Schur, Hilbert, Poincaré;
 who were at least inspired by hyperbolic geometry.
On the other hand,  hyperbolic geometry as a standard tool or setting have became ubiquitous in several branches of mathematics.
Thus, in terms of modern mathematics (synthetic) hyperbolic geometry came to the same position as
 (synthetic) Euclidean geometry; and, in fact, even more so.
Therefore, more recent introductions to hyperbolic geometry have students in mind, telling a good story
 (sometimes quite passionately, cf.~G.~Horváth \cite{GH1}), or they are shorter model based expositions
 which are often subordinated to the needs of particular mathematical areas.
(Although also to tell a somewhat different story:
After Word War II, the work of Bolyai was judged in Hungary to be of sufficiently national and international nature
 to take pride in and be promoted academically.
This led to modern critical translations \cite{Bolyai2}, \cite{Bolyai3} of Bolyai;
 to the edition of Lobachevsky \cite{Lob05};
 to the publication of the translations Kerékjártó \cite{Ker1}, \cite{Ker2}
 about reworking   Hilbert \cite{Hil1} and F. Schur \cite{Sch} respectively;
 and to the works Rédei \cite{Red} and Szász \cite{Sza}, written roughly in the spirit of Klein
 and Liebmann respectively.
In fact, hyperbolic geometry became a feature of high school math teacher's training,
  cf.~Hajós, Strohmajer \cite{HS}, actually  much following Bolyai's elegant geometric arguments, or Szenthe, Juhász \cite{SzJ}.)

Nevertheless, I agree that modern economical expositions to
 hyperbolic geometry (like Berger \cite{Ber}) are sufficient to contain the description of a relevant couple of models
 with the descriptions of the canonical transition maps between them.
Then the exposition can be amended with more details according to needs, most typically,
 with the discussion of the congruence groups, distance or angle formulas, the description of cycles,
 trigonometry, the relationship between area and angular defects,
 infinitesimal arc length and area elements (which I would call standard topics).
On the other hand, experience tells that practice with hyperbolic geometry sooner or later also develops
 a kind of synthetic viewpoint.

\snewpage
In any case, for many people, hyperbolic geometry signifies a cultural experience.
And indeed, sometimes people just get amused with various features of hyperbolic geometry, without any particular consequence.
In terms of recreational mathematics, however, hyperbolic geometry is not as easygoing as Euclidean geometry:
 elementary hyperbolic geometry is understood to involve more analysis (and computations in general) than elementary Euclidean geometry.
(And computations can indeed be quite effective, but possibly sparing us from  geometric insight).
Combined with the fact that practitioners of hyperbolic geometry are more likely to be judged
 by mathematicians' standards, there is some unease in (finding places for) publishing curiosities in hyperbolic geometry.
Even so, in the recent decades, hyperbolic probably got more popular.

\textbf{The hyperbolic plane and the Beltrami--Cayley--Klein model.}
Hyperbolic plane geometry can be introduced through the unit disk based Beltrami--Cayley--Klein model in an extremely
 straighforward manner:
Let the $h$-points be the points of the open unit disk in the plane; let the
 $h$-lines be the nonempty intersections ordinary lines and the open unit disk
 (in other terminology: the nonempty traces of ordinary lines on the open unit disk).
In what follows, `$h$-' means  `in the viewpoint of hyperbolic geometry'.
Then this incidence geometry is the unit disk based Beltrami--Cayley--Klein model;
 and, more generally, the geometry of the hyperbolic plane is any incidence geometry
 abstractly isomorphic to this incidence model.
The automorphisms of the incidence model are the collineations (or $h$-collineations for the sake of clarity).
More generally, the geometry of the hyperbolic plane is everything which can naturally be built up from its
 incidence structure, that is up to $h$-collineations.
(Curiously, hyperbolic geometry is a pure collineation geometry.
This is in contrast to the Euclidean geometry, whose collineation geometry is the affine geometry,
 and to the spherical geometry, whose collineation geometry is the doubly covered projective geometry.)

Natural elements of the build-up are as follows:
Although we cannot naturally reconstruct the ambient Euclidean space of the unit disk from the incidence geometry,
 we can reconstruct the ambient projective plane (i.~e.~the projective closure of the Euclidean plane).
One argument for that is that using the theorem of Desargues, we can test the projective concurrency of $h$-lines
 in the $h$-plane, therefore traces of the pencils can be characterized, leading to the points of the ambient projective
 plane abstractly; etc.
In particular, $h$-collineations will extend to some collineations of the ambient projective plane abstractly.
Due to the fundamental theorem of projective geometry we can overview the ($h$-)collineations easily.
In the view of BCK modell, we find that the $h$-collineations are exactly the traces of those collineations
 which take the unit circle into itself.
(In terms of the abstract reconstruction of the projective plane the unit circle corresponds to the so-called ``absolute'',
 the quadric of the ``asymptotic points'').
Thus the group of $h$-collineations is naturally isomorphic to $\mathrm{PO(2,1)}\simeq \mathrm{O}^{\uparrow}(2,1)$.
This overview of the $h$-collineations is sufficient to allow us the use of projective tools in hyperbolic geometry.
In particular, using the projective structure, we can naturally topologize the $h$-plane.
(Although topologization is very straightforward even directly from the incidence structure, as open half-planes are very easy to define.)
Using the projective structure, and more precisely (the quadric of) the absolute,
 we can define natural metrics on the $h$-plane depending on an arbitrarily chosen scaling parameter $k>0$,
 such that in BCK modell one has
\begin{equation}
\mathrm{dist}^{[k]}_{\mathrm{BCK}}\left((x_1,y_1),(x_2,y_2)\right)
=k\arcosh\frac{1-x_1x_2-y_1y_2}{\sqrt{1-(x_1)^2-(y_1)^2}\sqrt{1-(x_2)^2-(y_2)^2}}.
\plabel{eq:kdist}
\end{equation}
In fact, these metrics can be characterized as those $h$-collineation invariant metrics, where
 the distance of the $h$-midpoint ($h$-symmetry center) of two points from any of the original points is
 half of the distance of the two original points from each other.
The choice $k$ leads to inequivalent metrics.
In a natural but complicated way, $\mathrm{dist}^{[k]}_{\mathrm{BCK}}$ can be characterized further as the metric which leads to
 constant Gauss curvature $-\frac1{k^2}$; and a in much simpler but less natural way, $\mathrm{dist}^{[k]}_{\mathrm{BCK}}$
 can be characterized further as the metric
 where, in an asymptotic triangle, the distance of the center (the intersection of the symmetry lines) from the sides is $k\artanh\frac12$.
However, $k$ is a purely a scaling constant, and while it is sometimes retained in analogy to
 spherical geometry, the choice $k=1$ is most natural; it can be said to signify \textit{the} natural hyperbolic metric.
In what follows we use $k\equiv1$.
(Note, however, that in older texts the idea of ``imaginary geometry'' is taken seriously.
In those texts $k_{\mathrm{old}}=\mathrm ik$ is used, and for example, $\cos\frac{d}{k_{\mathrm{old}}}$ is used instead of $\cosh\frac dk$.
For example, Coolidge \cite{Coo} (1909) is quite dismissive toward using hyperbolic cosine,
 and even Klein \cite{Kle3} (1928) uses it only half-heartedly; however Liebmann \cite{Lib1} (1905) and
 Sommerville \cite{Som} (1914) are already according to modern usage.
Nevertheless those older works on ``imaginary geometry'' are strangely modern and of mathematician-like in a certain sense.)
In the BCK model, the $h$-segments connecting $h$-points (either inferred from the topology of $h$-lines or $h$-metrically)
 correspond to ordinary Euclidean segments connecting the respective points of the unit disk, etc.
Ultimately, notions which can be defined naturally in the unit disk based BCK model yield a model of hyperbolic plane geometry
 once the axiomatics are clarified, and apparently even before that; and the higher dimensional generalizations are also straightforward.

Other features associated to the BCK model will be recalled later in text as needed.
However, the reader will, in general, be expected to be familiar with the BCK model and be able to do computations there.
Anyway, the BCK model offers a relatively simple way to present the hyperbolic geometry,
 and also to carry out at least certain computations  effectively.

\snewpage
\textbf{Hyperbolic conics.}
The systematic study of hyperbolic conics was apparently started by Story \cite{Sto} (1882).
The analytic approach is augmented by more synthetic interpretations by Killing \cite{Kil3} (1885).
Then hyperbolic conics are discussed by
 D'Ovidio \cite{DO1}, \cite{DO2}, \cite{DO3} (1891), Barbarin \cite{Bar} (1901), Liebmann \cite{Lib0} (1902),
 Massau \cite{Mas} (1905), Liebmann \cite{Lib1}, Coolidge \cite{Coo}.
A very detailed study of the conics is conducted by V\"or\"os \cite{V1} / \cite{V2} (1909/10).
(For the sake of comparison, a parallel record for quadratic surfaces is provided by
 with respect to Barbarin \cite {Bar}, Coolidge \cite{CooP} (1903), Bromwich \cite{Bro} (1905), Coolidge \cite{Coo},
 and V\"or\"os \cite{V3} (1912).)
Thus, it can be said that by 1910 the hyperbolic conics were rather well-explored;
 although not all of the sources cited above can be considered equally accessible.
Ultimately, the point is that $h$-conics are simplest to be defined as conics of
 the projective models (most simply: of the Beltrami--Cayley--Klein model).
Then they enjoy similar synthetic presentations (essentially: focal descriptions) and properties as in the Euclidean case.
Sommerville \cite{Som} and Klein \cite{Kle3} also discuss hyperbolic conics; but,
 in the lack of particular applications, the general interest in them declines.
(Compare Liebmann \cite{Lib1}, \cite{Lib2}, \cite{Lib3}.)
Among later works some notable ones are as follows:
Arithmetic aspects of the classification of the hyperbolic conics are discussed further by Fladt \cite{F1}, \cite{F2}, \cite{F3}.
Hyperbolic conics are readdressed using various approaches by
 Epshtein \cite{Eps}, Pevzner \cite{Pev0}, \cite{Pev1}, \cite{Pev2},  Moln\'ar \cite{Mol1}, \cite{Mol2}, \cite{Mol}.
(Also see the literature cited therein.)
Recently,  hyperbolic conics gained some moderate publicity, cf.~G.~Horváth \cite{GH2},
 Weiss \cite{W}, Izmestiev \cite{Izm}, Csima, Szirmai \cite{CsSz}, Božić Dragun, Koncul \cite{BDK}, and the literature cited therein.
Actually, for those who are not specialists of hyperbolic geometry, the
 exposition by Izmestiev \cite{Izm} may be a good choice for a start.
However, beyond cycles, in this paper we will deal only with one kind of hyperbolic conic: the $h$-elliptic parabola, which we explain.
But, familiarity with cycles  will be required (as part of the standard topics).
\snewpage

\section{Hyperbolic elliptic parabolas, an introduction}\plabel{sec:hep}
The name and analytic description of $h$-elliptic parabolas is due to Story \cite{Sto}.
Considering the classification of hyperbolic conics the name is quite apt.
(One can browse Fladt \cite{F1} for the classification but there Figures 11, 12, 13, 14
 should be relabeled to Figures 12, 14, 11, 13 respectively; or, alternatively, the collection of the
 actual pictures in Figures 11--14 should be rotated clockwise.
Or, one can also see Figures 24--26 in Izmestiev \cite{Izm}.)
Next, we will dicuss the focus (or foci) of $h$-elliptic parabolas.

\textbf{Analytic approach to the foci.}
Let us consider the canonical equation
\begin{equation}
\frac{x^2}{C^2}+2y^2-2y\stackrel{\leftrightarrow}=0.
\plabel{eq:nhep}
\end{equation}
(Here $\stackrel{\leftrightarrow}=$ means that we do take equality really seriously but only formally and up to nonzero scalars.)
If we pass to homogeneous coordinates $[x_1,x_2,x_3]$ (with $x=\frac{x_1}{x_3}$ and $y=\frac{x_2}{x_3}$), then
\eqref{eq:nhep} transcribes to
\begin{equation}
\begin{bmatrix} x_1\\x_2\\x_3\end{bmatrix}^\top
\begin{bmatrix} \frac1{C^2}&0&0\\0&2&-1\\0&-1&0\end{bmatrix}\begin{bmatrix} x_1\\x_2\\x_3\end{bmatrix}\stackrel{\leftrightarrow}=0.
\plabel{eq:nhepm}
\end{equation}
The dual conic, in terms of dual homogeneous coordinates $[\xi_1,\xi_2,\xi_3]'$, is
\begin{equation}
\begin{bmatrix} \xi_1\\\xi_2\\\xi_3\end{bmatrix}^\top
\begin{bmatrix} \frac1{C^2}&0&0\\0&2&-1\\0&-1&0\end{bmatrix}^{-1}
\begin{bmatrix} \xi_1\\\xi_2\\\xi_3\end{bmatrix}
\equiv
\begin{bmatrix} \xi_1\\\xi_2\\\xi_3\end{bmatrix}^\top
\begin{bmatrix} C^2&0&0\\0&0&-1\\0&-1&-2\end{bmatrix}
\begin{bmatrix} \xi_1\\\xi_2\\\xi_3\end{bmatrix}
\stackrel{\leftrightarrow}=0.
\plabel{eq:dnhepm}
\end{equation}
That is
\begin{equation}
C^2(\xi_1)^2-2\xi_2\xi_3-2(\xi_3)^2\stackrel{\leftrightarrow}=0.
\plabel{eq:dnhep}
\end{equation}
Similarly, the absolute has the normal formal equation
\begin{equation}
x^2+y^2-1\stackrel{\leftrightarrow}=0,
\plabel{eq:na}
\end{equation}
 with dual conic
\begin{equation}
(\xi_1)^2+(\xi_2)^2-(\xi_3)^2\stackrel{\leftrightarrow}=0.
\plabel{eq:dna}
\end{equation}
Then a simple calculation yields that the pencil generated by \eqref{eq:dnhep} and \eqref{eq:dna} contains only two singular conics.
These are
\begin{equation}
C^2\cdot\left((\xi_1)^2+(\xi_2)^2-(\xi_3)^2\right)-\left(C^2(\xi_1)^2-2\xi_2\xi_3-2(\xi_3)^2\right) \stackrel{\leftrightarrow}=0.
\plabel{eq:pam1}
\end{equation}
and
\begin{equation}
1\cdot\left((\xi_1)^2+(\xi_2)^2-(\xi_3)^2\right)-
\left(C^2(\xi_1)^2-2\xi_2\xi_3-2(\xi_3)^2\right) \stackrel{\leftrightarrow}=0.
\plabel{eq:pam2}
\end{equation}
The first one is a pair of lines in the dual projective plane that is a pair of points in the ambient projective plane,
 while the second one is only a pointellipse (a pair of imaginary lines) in the dual projective plane.

In fact,  \eqref{eq:pam1} simplifies to
\[(\xi_2+\xi_3)(C^2\xi_2+(2-C^2)\xi_3)\stackrel{\leftrightarrow}=0.\]
Here the factor $\xi_2+\xi_3\stackrel{\leftrightarrow}=0$ corresponds to the point $[0,1,1] =(0,1)$,
  and the factor $C^2\xi_2+(2-C^2)\xi_3\stackrel{\leftrightarrow}=0$ corresponds to the point
  $[0,C^2,2-C^2] =\left(0,\frac{C^2}{2-C^2}\right)$.
For this reason we can call $(0,1)$ and $\left(0,\frac{C^2}{2-C^2}\right)$ the foci of our hyperbolic conic,
 but only the second one is a proper hyperbolic point, while the first one can be described in simpler terms
 as the only asymptotic point of our conic.
(In contrast,   \eqref{eq:pam2} simplifies to
\[(1-C^2)(\xi_1)^2+(\xi_2+\xi_3)^2\stackrel{\leftrightarrow}=0,\]
 from which we cannot extract two projective points in the previous manner.)

It may be strange at first sight, but this hocus pocus can convince people with experience in analytic projective geometry
 about the designation of the focus (or foci) of our $h$-elliptic parabola.

[There is an analogy to the Euclidean case:
Let us consider, say, the parabola
\[x^2-2py\stackrel{\leftrightarrow}=0,\]
 where $p>0$.
The dual conic is
\begin{equation}
(\xi_1)^2-\frac2p\,\xi_2\xi_3\stackrel{\leftrightarrow}=0.
\plabel{eq:ppam1}
\end{equation}
The dual Euclidean absolute is given by
\begin{equation}
(\xi_1)^2+(\xi_2)^2\stackrel{\leftrightarrow}=0.
\plabel{eq:ppam2}
\end{equation}
(This object exists only in the dual picture.
It encodes the Euclidean metric up to a nonzero scalar multiple.
Do not worry, there is no natural distance scale in true Euclidean geometry, anyway.)
There are two singular conics in the pencil generated by \eqref{eq:ppam1} and \eqref{eq:ppam2}:
The first one is
\begin{equation}
1\cdot\left((\xi_1)^2+(\xi_2)^2\right)-\left((\xi_1)^2-\frac2p\,\xi_2\xi_3\right)\stackrel{\leftrightarrow}=0,\plabel{eq:ppam3}
\end{equation}
 and the second one is \eqref{eq:ppam2} itself.
Of these, the first one is a pair of lines on the dual projective space:
\eqref{eq:ppam3} simplifies as
\[\xi_2\left(\xi_2 +\frac2p\,\xi_3\right)\stackrel{\leftrightarrow}=0.\]
The factor $\xi_2\stackrel{\leftrightarrow}=0$ corresponds to the point $[0,1,0]$, which is the ideal point of the parabola;
 the factor $\xi_2 +\frac2p\,\xi_3\stackrel{\leftrightarrow}=0$  corresponds to the point $[0,1,\frac 2p]=\left(0,\frac p2\right)$.
These are the two foci of the parabola, but only the second one is a proper Euclidean point.]

This closes the analytic geometrical motivation for the focus of our $h$-elliptic parabolas.
 (The general situation for $h$-conics is more complicated, cf.~Story \cite{Sto}, Izmestiev \cite{Izm}).

\snewpage
\textbf{Synthetic approach to the foci.}
A synthetic interpretation for the $h$-elliptical parabola, and in particular, for the focus was given by Killing \cite{Kil3},
 but in detail by Liebmann \cite{Lib0}.
According to the synthetic interpretation, the (ordinary) points our canonical $h$-elliptic parabola
 are the $h$-points which are of equal distance from the focus point $F^C=\left(0,\frac{C^2}{2-C^2}\right)$ and
 from the horocycle of equation
\begin{equation}
%C^2\cdot(x^2+y^2-1)+(1-C^2)\cdot(x^2-2y^2-2y)\equiv
(x^2+y^2-1)+(1-C^2)(y-1)^2\equiv
 x^2 + (y-1)((2-C^2)y+C^2) =0.
\plabel{eq:dhor}
\end{equation}
(Cf.~Figure \ref{fig:figHEP01}.(f).)
This interpretation is not at all surprising as any regularized distance from an asymptic point (and, in particular, from the
 asymptotical focus) is the oriented distance from a horocycle with the given asymptotic point.

We will demonstrate this characterization:
In general the distance of a point and a horocycle is always measured through an axis
 of the horocycle, i.~e.~a line through the asymptotic point of the horocycle.
As the horocycle \eqref{eq:dhor} has asymptotic point $(0,1)$ it is reasonable
 to parametrize our $h$-elliptic parabola and the horocycle according to the intersections
 with the lines $l_t\colon x+t(y-1)=0$ where $t\in\mathbb R$.
In this parametrization the corresponding point of the
 $h$-elliptic parabola is
\[P^C(t)=\left(\frac{2tC^2}{2C^2+t^2},\frac{t^2}{2C^2+t^2} \right),\]
 and the corresponding point of the horocycle is
 \[H^C(t)=\left(\frac{2t}{2-C^2+t^2},\frac{-C^2+t^2}{2-C^2+t^2 } \right).\]
(Cf.~Figure \ref{fig:figHEP01}.(f).) Then we can check
\[\dist_{\mathrm{BCK}}(H^C(t),P^C(t))=\dist_{\mathrm{BCK}}(P^C(t),F^C).\]
Indeed, both sides clear to
\[=\arcosh\frac{C^2(2-C^2 )+(1-C^2)t^2}{2C\sqrt{(1-C^2)(C^2+(1-C^2)t^2)}}
=\artanh\frac{C^2C^2+(1-C^2)t^2}{C^2(2-C^2 )+(1-C^2)t^2}.\]
\begin{commentx}
\[=\arsinh\frac{C^2C^2+(1-C^2)t^2 }{2C\sqrt{(1-C^2)(C^2+(1-C^2)t^2)}}
=\ln\frac{\sqrt{C^2+(1-C^2)t^2}}{C\sqrt{1-C^2}}.\]
\end{commentx}
On the other hand, $l_t$ and the perpendicular bisector of $H^C(t)$ and $F^C$ can meet at most at one point,
 therefore  $X=P^C(t)$ is the only $h$-point on $l_t$ such that  $\dist_{\mathrm{BCK}}(H^C(t),X)=\dist_{\mathrm{BCK}}(X,F^C)$ holds.
This proves that Killing's characterization is valid.

(The directrix horocycle is not to be confused with the directrix line $(3 C^2-2)y\stackrel{\leftrightarrow}=C^2$, the polar of the focus $F^C$
 with respect to the $h$-elliptic hyperbola.
That has some geometric significance, see Story \cite{Sto} or Izmestiev \cite{Izm} for more on that,
 but it does not even necessarily yield a proper hyperbolic line.)

\snewpage
\textbf{Geometric elements related synthetically.}
The horocycle of equation
\begin{equation}
(x^2+y^2-1)+(y-1)^2\equiv x^2 +2y^2-2y =0.
\plabel{eq:shor}
\end{equation}
 can be called the supporting horocycle of the canonical $h$-elliptic parabola. Cf.~Figure \ref{fig:figHEP01}.(a).
(Other name for `horocycle' is `paracyle', a very expressive name, unfortunately not used much in English.)
Synthetically, it yields the smallest horodisk which contains the $h$-elliptic parabola.
The supporting horocycle meets the $h$-elliptic parabola at a single distinguished point, at the vertex of $h$-elliptic parabola.
(Here, in the canonical case, it is $(0,0)$.)
The $h$-line connecting the vertex and the asymptotic point of the $h$-elliptic parabola is it axis;
  which is obviously the singe symmetry axis of the $h$-elliptic parabola.
(Here, in the canonical case, it is the line $x=0$.)
The supporting distance band $\frac{x^2}{C^2}+ y^2-1\leq0$
 is synthetically the closure of the union of the distance lines who intersect the $h$-elliptic parabola
 and have the same axis as the $h$-elliptic parabola.
(Other name for `distance line' is `hypercycle', a name which is sometimes used.)
(n the canonical picture this can be overviewed very simply, as there  $h$-elliptic parabola is the $C$-squeezed
 version of the supporting horocycle and the supporting distance band is the $C$-squeezed
 version of the unit disk.
The supporting half-plane $y\geq0$ is half-plane which contains the $h$-elliptic parabola and its
 boundary line is tangent at its vertex.
The supporting half distance band is the intersection of the supporting distance band and the supporting half-plane.
(Cf.~Figure \ref{fig:figHEP01}.(b).)
The focus can be characterized as the intersection points of the supporting horocycle and the boundary of the
 supporting distance band $h$-perpendicularly projected the the axis.
(In canonical picture the intersection points are $\left(\pm\frac{2C\sqrt{1-C^2}}{2-C^2},\frac{C^2}{2-C^2}\right)$;
 by centrality, the Euclidean perpendicular projection is faithful regarding the model.
Cf.~Figure \ref{fig:figHEP01}.(c).)
The directrix horocycle \eqref{eq:dhor} can be characterized as horocycle through the asymptotic point of the $h$-elliptic parabola
 and through the vertices of supporting half distance band.
(In the canonical picture the latter vertices are $(\pm C,0)$. Cf.~Figure \ref{fig:figHEP01}.(d).)
Other characterizations of the focus and the directrix horocycle relative to each other are as follows:
The distance between the focus and the vertex is the focal distance.
(Here, in the canonical case, it is $\artanh\frac{C^2}{2-C^2}$.)
If we translate the supporting  horocycle by the focal distance along any of its axes, away from the asymptotic point,
 then we obtain the directrix horocycle.
Conversely, the directrix horocycle reflected to tangent at the vertex intersects the axis in the focus.
(Cf.~Figure \ref{fig:figHEP01}.(e).)
From the earlier discussions it must be clear that the asymptotic point, the supporting horocycle, the vertex, the axis, the
 (non-asymptotic) focus, the directrix horocycle can all be simply characterized synthetically in relation to the $h$-elliptic parabola.

Also, if a (directrix) horocycle and a (focus) point in its interior are given, then they can be moved into canonical position in order to
 yield a canonical $h$-ellipse, etc.
Ultimately, $h$-elliptical parabolas can be imagined to be quite similar to ordinary Euclidean parabolas
 analytically and synthetically.

\snewpage
\textbf{Notable distances and measuring the $h$-elliptic parabola up.}
Before starting, let us remark here that $\dist_{\mathrm{BCK}}((0,0),(t,0))=\dist_{\mathrm{BCK}}((0,0),(0,t))=\artanh t$.

\begin{figure*}[b]
   \caption*{Fig.~\ref{fig:figHEP02}\quad Notable distances and measuring the $h$-elliptic parabola  up. $\rightarrow$  }
%\phantomcaption
\end{figure*}

  \begin{figure}[htp]
   \begin{subfigure}[b]{2.5in}
    \includegraphics[width=2.5in]{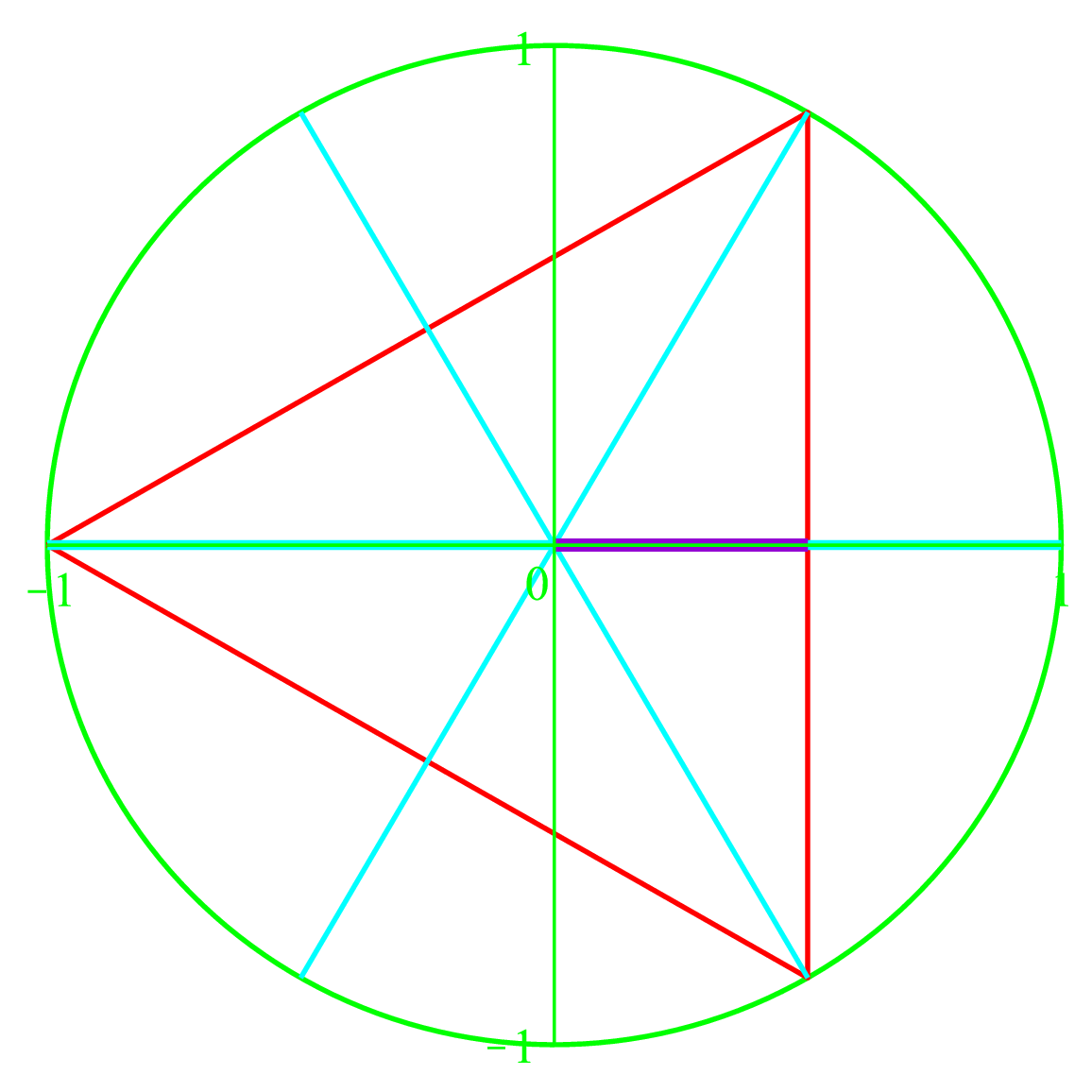}
    \caption*{\ref{fig:figHEP02}(a)  $\frac12\ln3=\arctan \frac12$ [violet].\\\phantom{  $ \sqrt2   \frac{\sqrt2}2 $}}
  \end{subfigure}
  \begin{subfigure}[b]{2.5in}
    \includegraphics[width=2.5in]{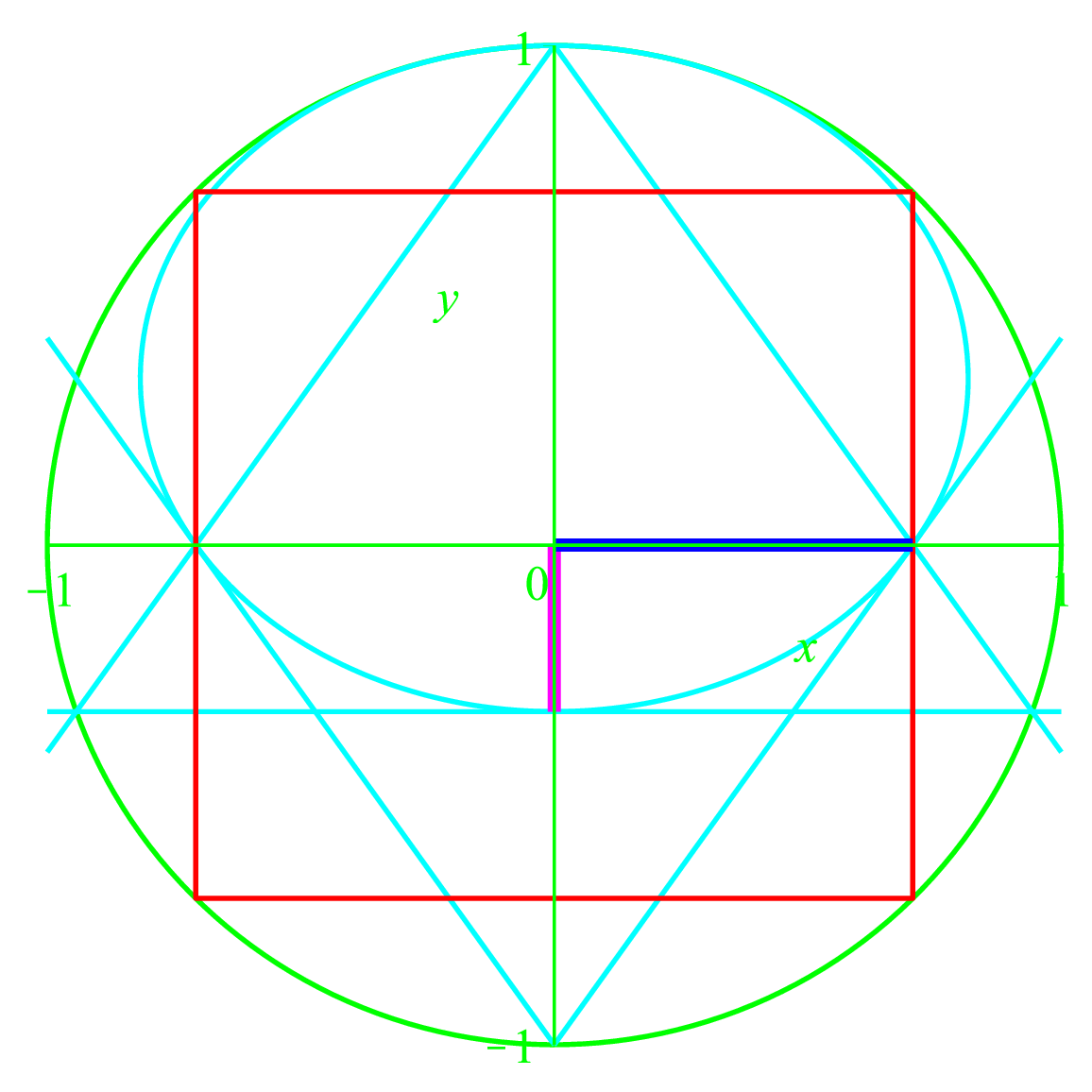}
    \caption*{\ref{fig:figHEP02}(b)  $\ln(1+\sqrt2)=\arctan\frac{\sqrt2}2=\sinh 1$\\ {}[blue],
    (and $\frac12\ln2=\arctan\frac13$ [magenta]).}
  \end{subfigure}
  \\
   \begin{subfigure}[b]{2.5in}
    \includegraphics[width=2.5in]{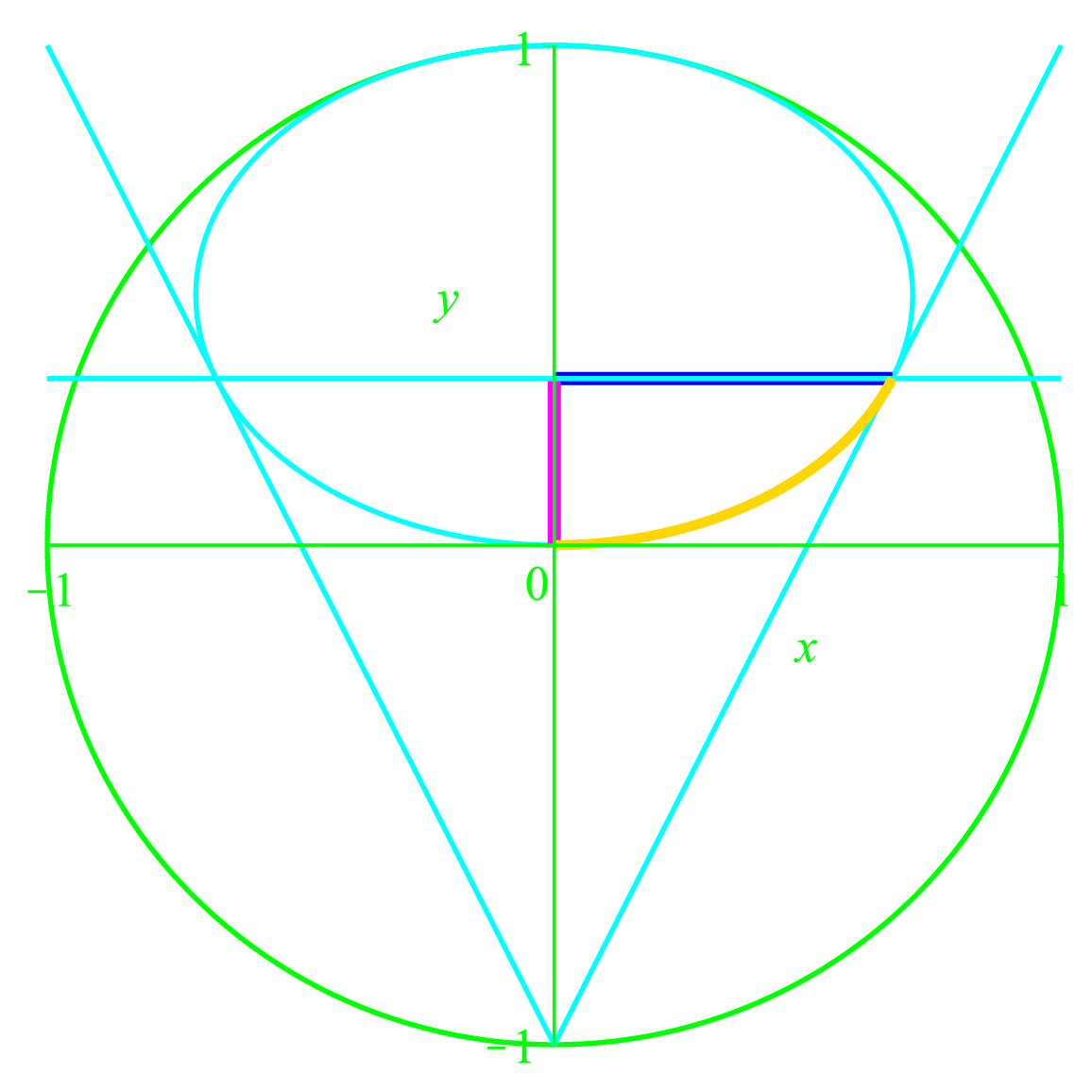}
    \caption*{ \ref{fig:figHEP02}(c) $\frac12\ln2=\arctan\frac13$ [magenta],\\
     $\ln(1+\sqrt2)=\arctan\frac{\sqrt2}2=\sinh 1$ [blue],\\
     and $1$ (in arc length) [gold]. \phantom{$\sqrt{\frac{C^2}{C^2}}$}}
  \end{subfigure}
  \begin{subfigure}[b]{2.5in}
    \includegraphics[width=2.5in]{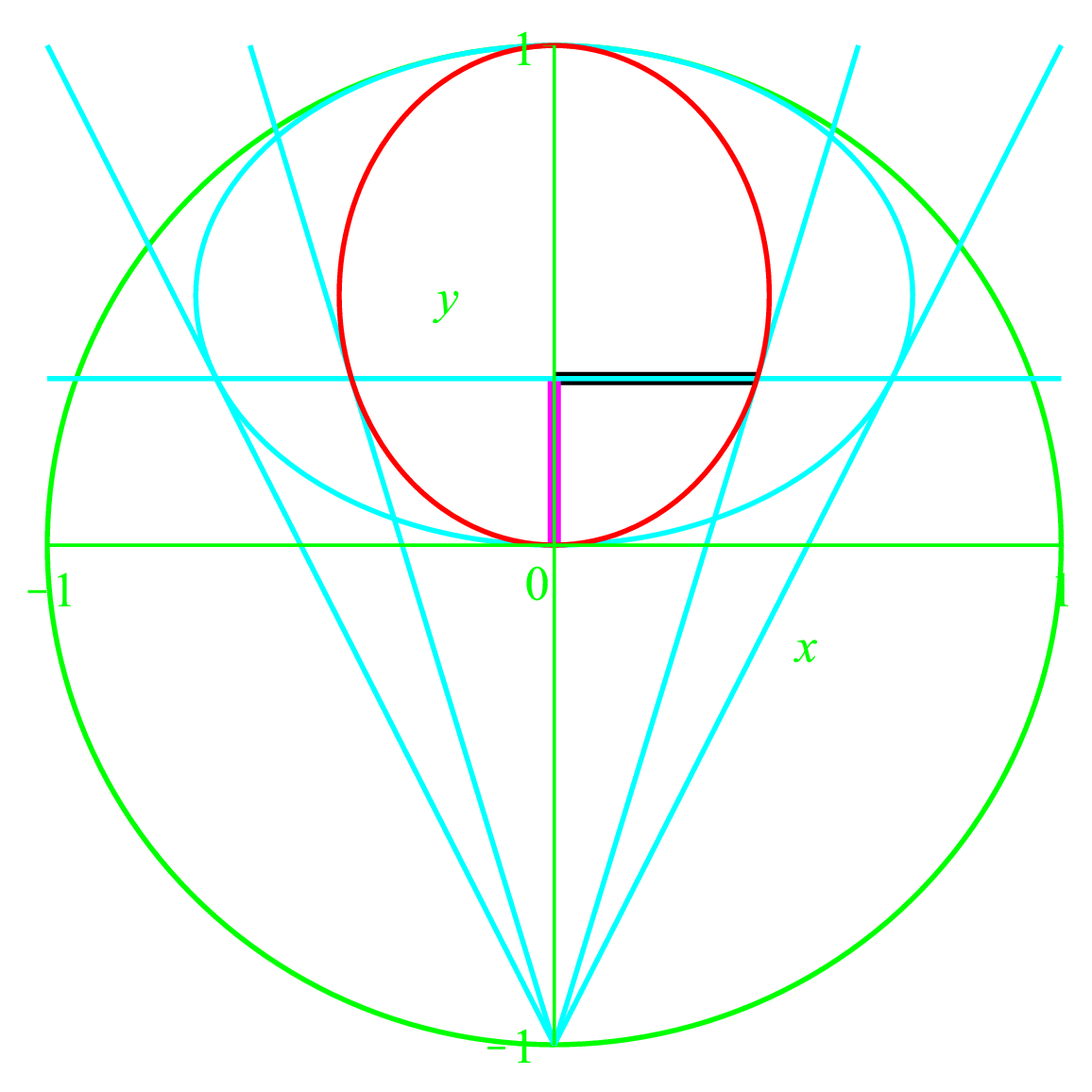}
    \caption*{\ref{fig:figHEP02}(d)
    $\frac12\ln2=\arctan\frac13$ [magenta]\\
    $\artanh(\frac{\sqrt2}2C)=\arsinh\sqrt{\frac{C^2}{2-C^2}}$ [black].\\
    \phantom{$\arctan\frac{\sqrt2}2$}
    }
  \end{subfigure}
     \begin{subfigure}[b]{2.5in}
    \includegraphics[width=2.5in]{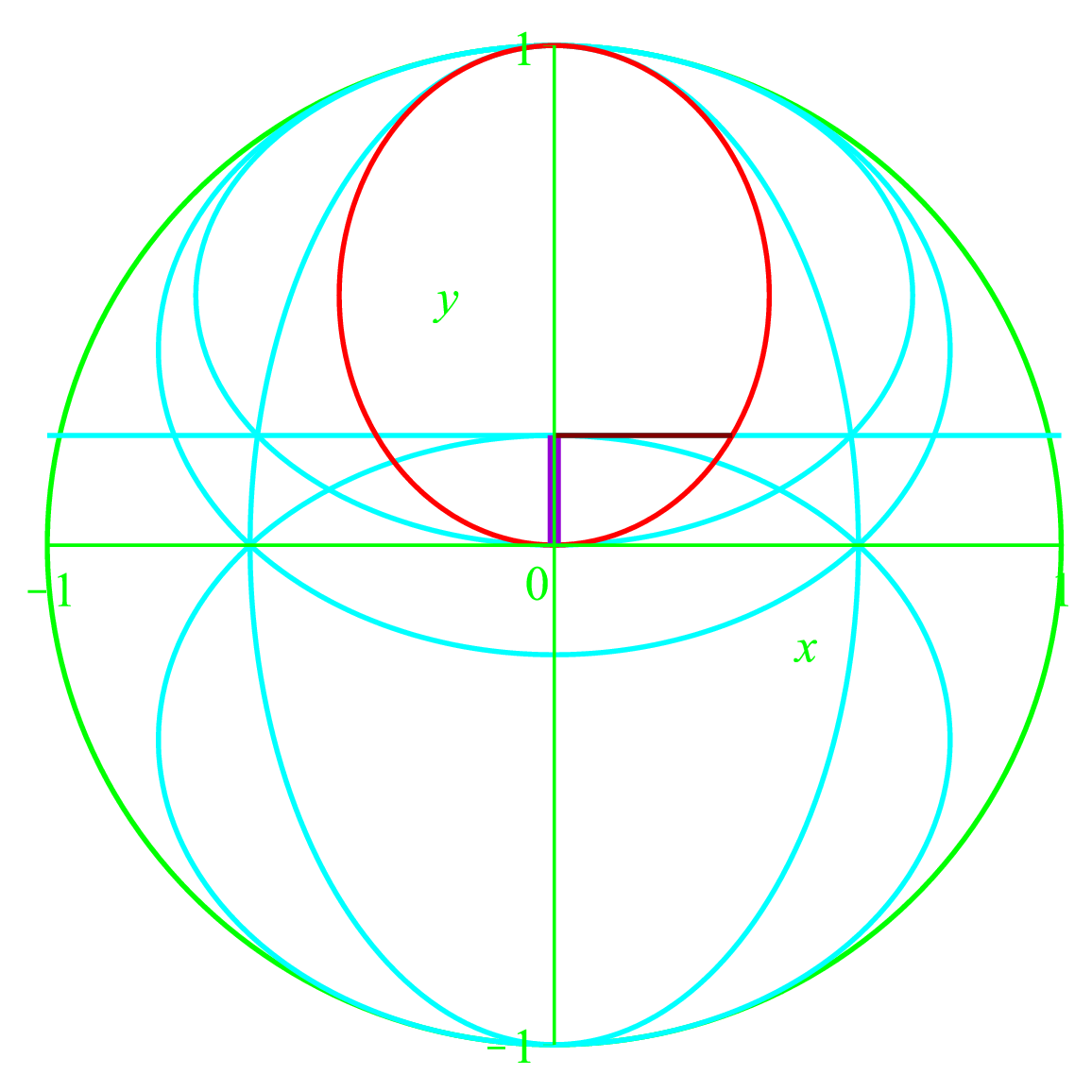}
    \caption*{\ref{fig:figHEP02}(e)  $\artanh \frac{C^2}{2-C^2}=\ln\frac1{\sqrt{1-C^2}}$ [violet],\phantom{$\sqrt{\frac11}$}\\
      $ \artanh C^2 $ [brown].  }
  \end{subfigure}
  \begin{subfigure}[b]{2.5in}
    \includegraphics[width=2.5in]{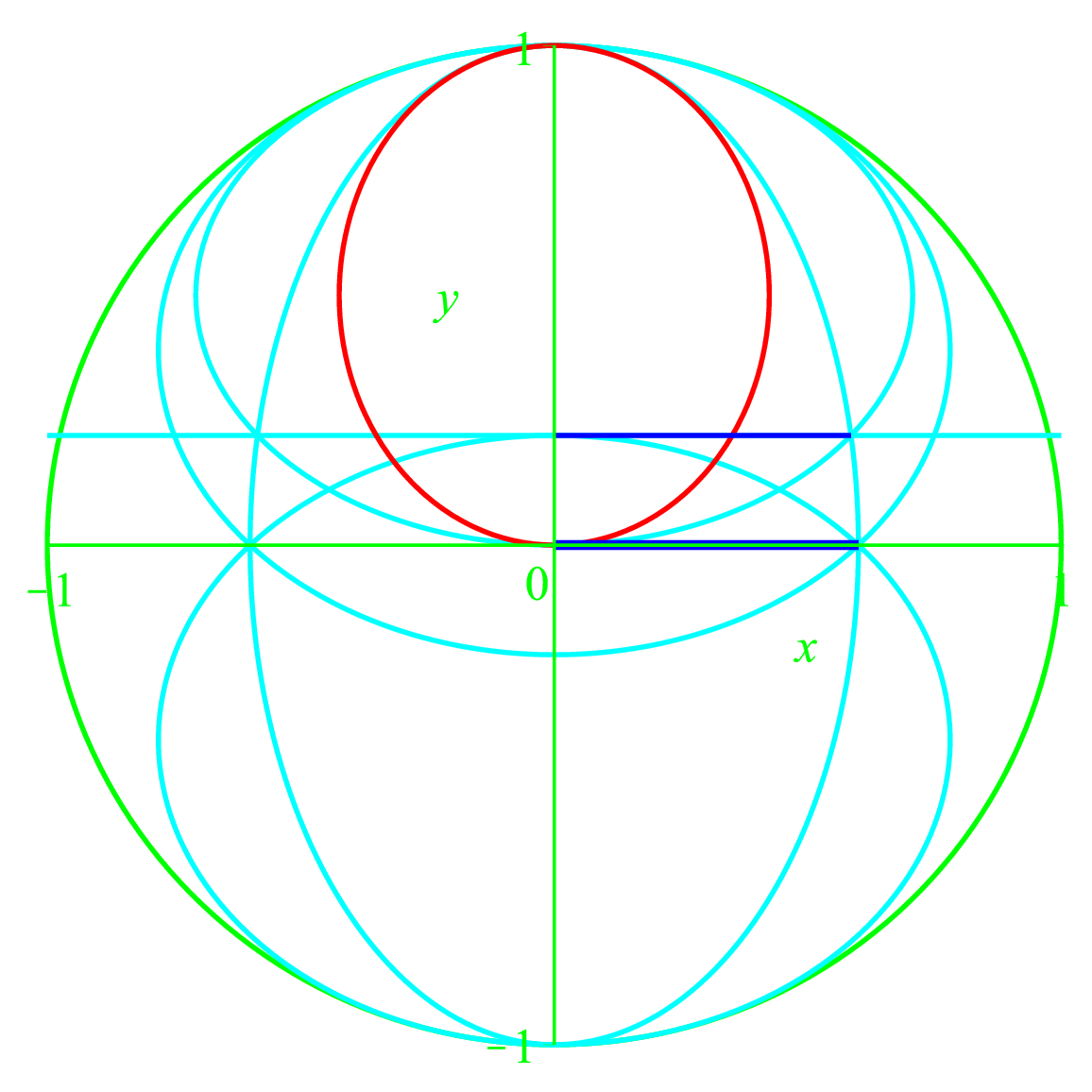}
    \caption*{\ref{fig:figHEP02}(f)  $\artanh C=\ln\sqrt{\frac{1+C}{1-C}}$ [blue].\\\phantom{$ \artanh C^2 $ [brown]}}
  \end{subfigure}
%   \caption*{Fig.~\ref{fig:figHEP02} Notable distances and measuring the $h$-elliptic parabola [red] up.  }
\phantomcaption
\plabel{fig:figHEP02}
\end{figure}

Earlier, we have noted that the distance $\frac12\ln3=\arctan \frac12$ can be characterized geometrically (with respect to $k=1$),
 using the asymptotic triangle, cf.~Figure \ref{fig:figHEP02}(a).
Another such distance is $\ln(1+\sqrt2)=\arctan\frac{\sqrt2}2=\sinh 1$, which can be   characterized geometrically
 using the asymptotic square cf.~Figure \ref{fig:figHEP02}(b).
However this distance, and also the distance $\frac12\ln2=\arctan\frac13$ can be characterized simply using a horocycle.
Indeed, if we take a horocycle and any external asymptotic point, then we can draw the two tangent lines to the horocycle,
 then in the resulting chordal bisectorial configuration the above mentioned
  distances can be characterized as Figure \ref{fig:figHEP02}(c) shows.
(If the horocycle is the canonical one, and the external asymptotical point is $(-1,0)$,
 then the tangent lines are $\pm2x-y-1=0$, the tangent points are $(\pm\frac23,\frac13)$, and the polar of the
 external asymptotic point with respect to the horocycle is $y=1/3$.)
This is also particularly nice in terms of characterizing the distance (scale) because the
 horocyclic arc facing the (any) external asymptotical point   has length exactly $2\cdot1$.
(See $\Len_\hyp(\partial^\cup E_\eta^1)=2\sqrt{\frac{2\eta}{1-\eta}}$ later (\eqref{eq:emon}) with $\eta= 1/3$.)

The characteristic distance $\frac12\ln2=\arctan\frac13$ is also present with respect to any $h$-elliptic parabola
 with a similar setup but the external asymptotic point be chosen as the asymptotic point of the axis
 not in the (closure) of the $h$-elliptic parabola (the ``anti-axial asymptotic point''), see Figure \ref{fig:figHEP02}(d).
This configuration of Figure \ref{fig:figHEP02}(d) can be used to measure $C$ without referring to the supporting distance band.
In Figure \ref{fig:figHEP02}(e) we see the focal distance $\artanh \frac{C^2}{2-C^2}=\ln\frac1{\sqrt{1-C^2}}$
 and the ``classical parameter'' $\artanh C^2$.
In Figure \ref{fig:figHEP02}(f) the radius $\artanh C$ of the supporting distance band is shown.
(This radius corresponds to the minor semiaxis length if the $h$-elliptic parabola is
 considered as a limiting object of $h$-ellipses.)
In synthetic terms, the hyperbolic cosine of the radius is equal to the exponential of the focal distance
($\cosh \artanh C=\frac1{\sqrt{1-C^2}}=\exp\artanh\frac{C^2}{2-C^2}$).
\snewpage
\section{The area}\plabel{sec:area}
\textbf{The computation of the area difference.}
Let
\begin{equation}
 E^C =\left\{(x,y)\,:\,\frac{x^2}{C^2}+2y^2-2y\leq0\right\},\plabel{eq:Edef}
\end{equation}
 and
\[B^C =\left\{(x,y)\,:\,\frac{x^2}{C^2}+ y^2-1\leq0\quad\text{ and }\quad 0\leq y\right\}.\]
Our primary objective here is computation of the hyperbolic area of $E^C\setminus B^C$.

It is well-known that the area density in the BCK model is
\[\mathrm{da}_{\hyp}(x,y)=\frac1{(1-x^2-y^2)^{3/2}}\,|\mathrm dx\wedge \mathrm dy|.\]
Using this, we can deal with the area problem effectively.

Set
\[B^C_\eta=B^C\cap\left\{(x,y)\,: y\leq \eta\right\}\]
 (cf.~Figure \ref{fig:figHEP05}(a)).
Then $B^C_\eta$ is a distance band belonging to the segment $\{0\}\times[0,\eta]$.
For a warm-up, we can compute its hyperbolic area:
\begin{align*}
\Area_{\hyp}(B^C_\eta)=&\int_{y=0}^\eta\int_{x=-C\sqrt{1-y^2}}^{C\sqrt{1-y^2}} \frac1{(1-x^2-y^2)^{3/2}}\,\mathrm dx\,\mathrm dy
\\=&\int_{y=0}^\eta \left[\frac{x}{(1-y^2)\sqrt{1-x^2-y^2}}\right]_{x=-C\sqrt{1-y^2}}^{C\sqrt{1-y^2}}\,\mathrm dy
\\=&\int_{y=0}^\eta\frac{2C}{ (1-y^2)\sqrt{1-C^2}}\,\mathrm dy
\\=&\frac{2C}{\sqrt{1-C^2}}\artanh \eta.
\end{align*}
\begin{commentx}
\[= \frac{2C}{\sqrt{1-C^2}}\arsinh \frac{\eta}{\sqrt{1-\eta^2}}
= \frac{2C}{\sqrt{1-C^2}}\arcosh \frac{1}{\sqrt{1-\eta^2}}
=\frac{2C}{\sqrt{1-C^2}}\ln\sqrt{\frac{1+\eta}{1-\eta}}\]
\end{commentx}
(Here and in what follows, when the result of an integration can be used for a primitive function
 we do note invoke the primitive function separately.)

This result, of course, does not require the full knowledge of the density function:
The area must be proportional to the hyperbolic length of the segment $\{0\}\times[0,\eta]$,
 which is $\arctan\eta$, and finding out the coefficient $\frac{2C}{\sqrt{1-C^2}}$ is not that hard,
 even if  a little analysis is involved in its derivation from the corresponding arc length formula.
In fact, both Lobachevsky and Bolyai were able to compute such this area  using their rectangular coordinate system:
They knew that the are signed area between an arc element and the axis is $(\sinh q)\,\mathrm dp$.
(For comparison to our setting: In the present situation, the rectangular coordinate system is attached to the line $x=0$, thus
 $p=(\sgn y)\dist_{\mathrm{BCK}}((0,0),(0,y))=\artanh y$ and $q=-(\sgn x)\dist_{\mathrm{BCK}}((0,y),(x,y))=-\artanh\frac{x}{\sqrt{1-y^2}}$.)
Regarding the distance lines here, we have $q=\pm \artanh C$ (constants).
This implies that the area to compute is $2\cdot \sinh(\artanh C)\cdot\arctan \eta=2\cdot \frac{C}{\sqrt{1-C^2}}\cdot\arctan \eta$.
(We  remark that Lobachevsky was particularly fascinated by analytical formulae coming from hyperbolical area computations.)
(In what follows, `known to Lobachevsky and Bolyai' can be understood, for the sake of simplicity,
 as reference to Lobachevsky \cite{Lob2} and Bolyai \cite{Bolyai4}.)

Next we compute the area of
\begin{equation}
E^C_\eta=  E^C\cap\left\{(x,y)\,: y\leq \eta\right\}
\plabel{eq:Edef2}
\end{equation}
$\subset B^C_\eta$  (cf.~Figure \ref{fig:figHEP05}(a)). Here
\begin{align*}
\Area_{\hyp}(E^C_\eta)=&\int_{y=0}^\eta\int_{x=-C\sqrt{2y(1-y)}}^{C\sqrt{2y(1-y)}} \frac1{(1-x^2-y^2)^{3/2}}\,\mathrm dx\,\mathrm dy
\\=&\int_{y=0}^\eta \left[\frac{x}{(1-y^2)\sqrt{1-x^2-y^2}}\right]_{x=-C\sqrt{2y(1-y)}}^{C\sqrt{2y(1-y)}}\,\mathrm dy
\\=&\int_{y=0}^\eta\frac{2C\sqrt{2y }}{ (1-y^2)\sqrt{1+y -2C^2y }}\,\mathrm dy
\\=&\frac{2C}{\sqrt{1-C^2}}\artanh\left(
{\frac { \sqrt {2\eta}\sqrt { 1-{C}^{2}}}{\sqrt { 1+\eta-2\,{C}^{2}\eta}
}}
\right)
-2\arctan\left( {\frac {C \sqrt {2\eta}}{\sqrt { 1+\eta-2\,{C}^{2}\eta}}}\right).
\end{align*}
\begin{commentx}
\[=\frac{2C}{\sqrt{1-C^2}}\arsinh\left( \sqrt{1-C^2}\sqrt{\frac{2\eta}{1-\eta}} \right)-2\arcsin\left( C\sqrt{\frac{2\eta}{1+\eta}} \right) \]
\[= \frac{2C}{\sqrt{1-C^2}}\arcosh\left(
 \frac {\sqrt { 1+\eta-2\,{C}^{2}\eta}} { \sqrt {1-\eta} }
\right)
-2\arccos\left(  \frac {\sqrt { 1+\eta-2\,{C}^{2}\eta}}{  \sqrt {1+\eta}}  \right)\]
\[= \frac{2C}{\sqrt{1-C^2}}\ln\left(
 \frac {\sqrt { 1+\eta-2\,{C}^{2}\eta} +  \sqrt {2\eta}\sqrt { 1-{C}^{2}} } { \sqrt {1-\eta} }
\right)
 -2\arctan\left( {\frac {C \sqrt {2\eta}}{\sqrt { 1+\eta-2\,{C}^{2}\eta}}}\right)\]
\end{commentx}
Rewriting $\frac{2C}{\sqrt{1-C^2}}\artanh\ldots$ into $\frac{2C}{\sqrt{1-C^2}}\ln\ldots$,
 using either $\artanh \eta=\ln\frac{\sqrt{1+\eta}}{\sqrt{1-\eta}}$, or,
\[\artanh \frac {\sqrt{a}}{\sqrt{A}} =\ln\frac{ \sqrt{A}+ \sqrt{a}}{\sqrt{A -a }}\qquad\qquad\qquad  (0<a<A), \]
\begin{commentx}
From this, using the addition formula
\[\artanh A -\artanh B= \artanh\left(\frac{A-B}{1-AB}\right),\]
 or
\[\arsinh A -\arsinh B= \arsinh\left(A\sqrt{1+B^2}-B\sqrt{1+A^2}\right)\]
\end{commentx}
 respectively, one can conclude that
\begin{multline*}
\Area_{\hyp}(B^C_\eta\setminus E^C_\eta)=\Area_{\hyp}(B^C_\eta)-\Area_{\hyp}(E^C_\eta)
\\= \frac{2C}{\sqrt{1-C^2}}\ln
 \frac { \sqrt {1+\eta} } {\sqrt { 1+\eta-2\,{C}^{2}\eta} +  \sqrt {2\eta}\sqrt { 1-{C}^{2}} }
 +2\arctan  {\frac {C \sqrt {2\eta}}{\sqrt { 1+\eta-2\,{C}^{2}\eta}}}.
\end{multline*}
\begin{commentx}
\[
= \frac{2C}{\sqrt{1-C^2}}\artanh\left(
{\frac {\sqrt{\eta} \left( 2\,{C}^{2}\eta+2\,{C}^{2}-\eta-2 \right) }{\sqrt{\eta}+ \left( \eta+1
 \right) \sqrt {2}\sqrt {1-{C}^{2}}\sqrt {1+\eta-2\,{C}^{2}\eta  }}
}
\right)
+2\arctan\left( {\frac {C\sqrt {2\eta}}{\sqrt {1+\eta-2\,{C}^{2}\eta}}}\right).
\]
\[ =\frac{2C}{\sqrt{1-C^2}}\arsinh\left({\frac {\eta \left( 2\,{C}^{2}\eta-\eta+2\,{C}^{2}-2 \right) }{
\left(\sqrt {\eta} \sqrt {2}\sqrt {1-{C}^{2}}+\eta\sqrt {1+\eta-2\,{C}^{2}\eta} \right) \sqrt {1+\eta}}}\right)
+2\arcsin\left( C\sqrt{\frac{2\eta}{1+\eta}}\right)\]
\[= \frac{2C}{\sqrt{1-C^2}}\ln\left(
 \frac { \sqrt {1+\eta} } {\sqrt { 1+\eta-2\,{C}^{2}\eta} +  \sqrt {2\eta}\sqrt { 1-{C}^{2}} }
\right)
 +2\arctan\left( {\frac {C \sqrt {2\eta}}{\sqrt { 1+\eta-2\,{C}^{2}\eta}}}\right)\]
\end{commentx}

Taking the limit $\eta\nearrow1$, we find
\[\Area_{\hyp}(B^C\setminus E^C)=\frac{2C}{\sqrt{1-C^2}}\ln\frac1{2\sqrt{1-C^2}}+2\arctan \frac{C}{\sqrt{1-C^2}}.\]
\begin{commentx}
\[=\frac{2C}{\sqrt{1-C^2}}\artanh\left( \frac{4C^2-3}{5-4C^2} \right) +2\arctan \frac{C}{\sqrt{1-C^2}}\]
\[=\frac{2C}{\sqrt{1-C^2}}\arsinh\left( \frac{C^2-\frac34}{\sqrt{1-C^2}} \right) +2\arcsin  C.\]
\end{commentx}
Thus, we can conclude that the half-distance band $B^C$ majorizes $h$-elliptical parabolical disk $E^C$ only by finite area, indeed.
\begin{commentx}
Note that
\[\artanh \frac{4C^2-3}{3-4C^2}=\artanh \frac{C^2}{2-C^2}-\artanh \frac35=\artanh \frac{C^2}{2-C^2}-\ln 2, \]
 and
\[\arctan \frac{C}{\sqrt{1-C^2}}=\arcsin C.\]
\end{commentx}
Somewhat rewritten, we have
\begin{equation}
\Area_{\hyp}(B^C\setminus E^C) =\frac{2C}{\sqrt{1-C^2}}\left(\artanh \frac{C^2}{2-C^2}-\ln 2\right)+2\arcsin C.
\plabel{eq:finar}
\end{equation}
\snewpage

\textbf{The interpretation of the area difference.}
Let $B^C_{\mathrm{up}\,\omega}$ be the set which is $B^C$ but translated by hyperbolic distance $\omega$ upward along the $y$-axis.
Note that the translated distance band is
\[B^C_{\mathrm{up}\,\omega}=\left\{(x,y)\,:\,\frac{x^2}{C^2}+ y^2-1\leq0\quad\text{ and }\quad \tanh \omega\leq y\right\}.\]
(Unfortunately, this invariant notation follows a philosophy different from the notation of $B^C_\eta$.
In fact the latter one could have been denoted more invariantly as $B^C_{\text{lin-cut}\, \omega}$ with $\omega=\artanh\eta$.)
Next, we ask for a half distance band $B^C_{\mathrm{up}\,\alpha(C)}$
 such that the area of  $B^C_{\mathrm{up}\,\alpha(C)}$ and $E^C$ should ``be equal''.
This, of course, is understood so that
\[\int_{x^2+y^2<1} \left(\chi_{ B^C_{\mathrm{up}\,\alpha(C)}}(x,y)-\chi_{E^C}(x,y) \right)\,\mathrm{da}_{\hyp}(x,y)=0,\]
 where $\chi$ denotes `characteristic function'.
Using \eqref{eq:finar} and the area formula for distance band segments, which tells
\[\int_{x^2+y^2<1} \left(\chi_{B^C}(x,y)-\chi_{B^C_{\mathrm{up}\,\omega}}(x,y) \right)\,\mathrm{da}_{\hyp}(x,y)
=\frac{C}{\sqrt{1-C^2}}\omega,\]
 we can answer this question immediately:
\[\alpha(C)= \artanh \frac{C^2}{2-C^2} -\ln 2+ \frac{\arcsin C}C\sqrt{1-C^2} .\]
Here the geometrically interesting quantity  $ \alpha(C)$, which is the distance of
 the ``supporting point'' of $B^C_{\mathrm{up}\,\alpha(C)} $ from the vertex of $h$-elliptic parabola.
This is not exactly the previously expected focal distance $\artanh \frac{C^2}{2-C^2}$ (even if it is nor far from that):

Although $\artanh \alpha(C)$ is monoton increasing in $C$; we find that
\[\lim_{C\searrow0 }\artanh \alpha(C)=1-\ln 2,\]
 and
\[\lim_{C\nearrow1 }\artanh \alpha(C)-\artanh \frac{C^2}{2-C^2}=-\ln 2,\]
 but, in general, the deviation
\[1-\ln2>\artanh \alpha(C)-\artanh \frac{C^2}{2-C^2}>-\ln 2.\]
 is bounded and monotone decreasing.
However,
\[ \lim_{C\nearrow1 } \frac{\Area_{\hyp}\left(B^C_{\mathrm{up}\,\artanh \frac{C^2}{2-C^2}}\setminus E^C\right)}
{\Area_{\hyp}\left(E^C\setminus B^C_{\mathrm{up}\,\artanh \frac{C^2}{2-C^2}}\right)}=1.\]
Thus, the half distance band $B^C_{\mathrm{up}\,\artanh \frac{C^2}{2-C^2}}$ supported
 at the focus point is not so bad after all for an approximation of $E^C$.
Nevertheless,
$\alpha(C)=\artanh \frac{C^2}{2-C^2}$ occurs only for $C=0.801986\ldots$

Let us mention that $\eqref{eq:finar}$ can be interpreted a bit more geometrically.
Let $A^C$ be the (partly) asymptotic triangle with vertices $(\pm C,0)$ and $(0,1)$.
This is given as
\[A^C=\left\{(x,y)\,:\, |x|\leq C(1-y)\quad\text{ and }\quad 0\leq y\right\}.\]
Its area can be computed as
\begin{align*}
\Area_{\hyp}(A^C )&=\int_{y=0}^1\int_{x=-C(1-y)}^{C (1-y)} \frac1{(1-x^2-y^2)^{3/2}}\,\mathrm dx\,\mathrm dy
\\&=\int_{y=0}^1\left[\frac{x}{(1-y^2)\sqrt{1-x^2-y^2}}\right]_{x=-C(1-y)}^{C(1-y)}\,\mathrm dy
\\&=\int_{y=0}^1  2\,{\frac {C}{ \left(  1+y \right) \sqrt { \left( 1-y \right)  \left( 1+y-{C}^{2}+{C}^{2}y \right) }}}
      \,\mathrm dy
\\&=\left[-2\arctan\frac{C\sqrt{1-y}}{\sqrt{1+y-C^2+C^2y}}\right]_{y=0}^1
\\&=2\arctan\frac{C}{\sqrt{1-C^2}}=2\arcsin C.
\end{align*}
\begin{commentx}
\[=\left[-2\arcsin\left( C\sqrt{\frac{1-y}{1+y}}\right)\right]_{y=0}^1
  =\left[-2\arccos\frac{\sqrt{1+y-C^2+C^2y}}{\sqrt{1+y}}\right]_{y=0}^1\]
\end{commentx}
This was, of course, a very uninformed computation again; computing the area of triangles
 is just much simpler from the angular defect:
At the non-asymptotic vertices, the angle is
$\Pi(\artanh C)=\arccos\tanh\artanh C=\arccos C=\frac\pi2-\arcsin C$.
That makes the area, that is the angular defect $\pi-0-2(\frac\pi2-\arcsin C)=2\arcsin C$.

Now we can state that
\begin{equation}
E^C=B^C_{\mathrm{up}\,\artanh \frac{C^2}{2-C^2}}+ B^C_{3/5}-A^C\qquad\text{in hyperbolic area}.
\plabel{eq:equ1}
\end{equation}
(Cf.~Figure \ref{fig:figHEP03}.)
This  can be understood   as
\[\int_{x^2+y^2<1} \left(\chi_{B^C_{\mathrm{up}\,\artanh \frac{C^2}{2-C^2}} }(x,y)
+\chi_{B^C_{3/5}}(x,y)
-\chi_{A^C}(x,y)
-\chi_{E^C}(x,y) \right)\,\mathrm{da}_{\hyp}(x,y)=0.\]

 \begin{figure}[htp]
   \begin{subfigure}[b]{2.5in}
    \includegraphics[width=2.5in]{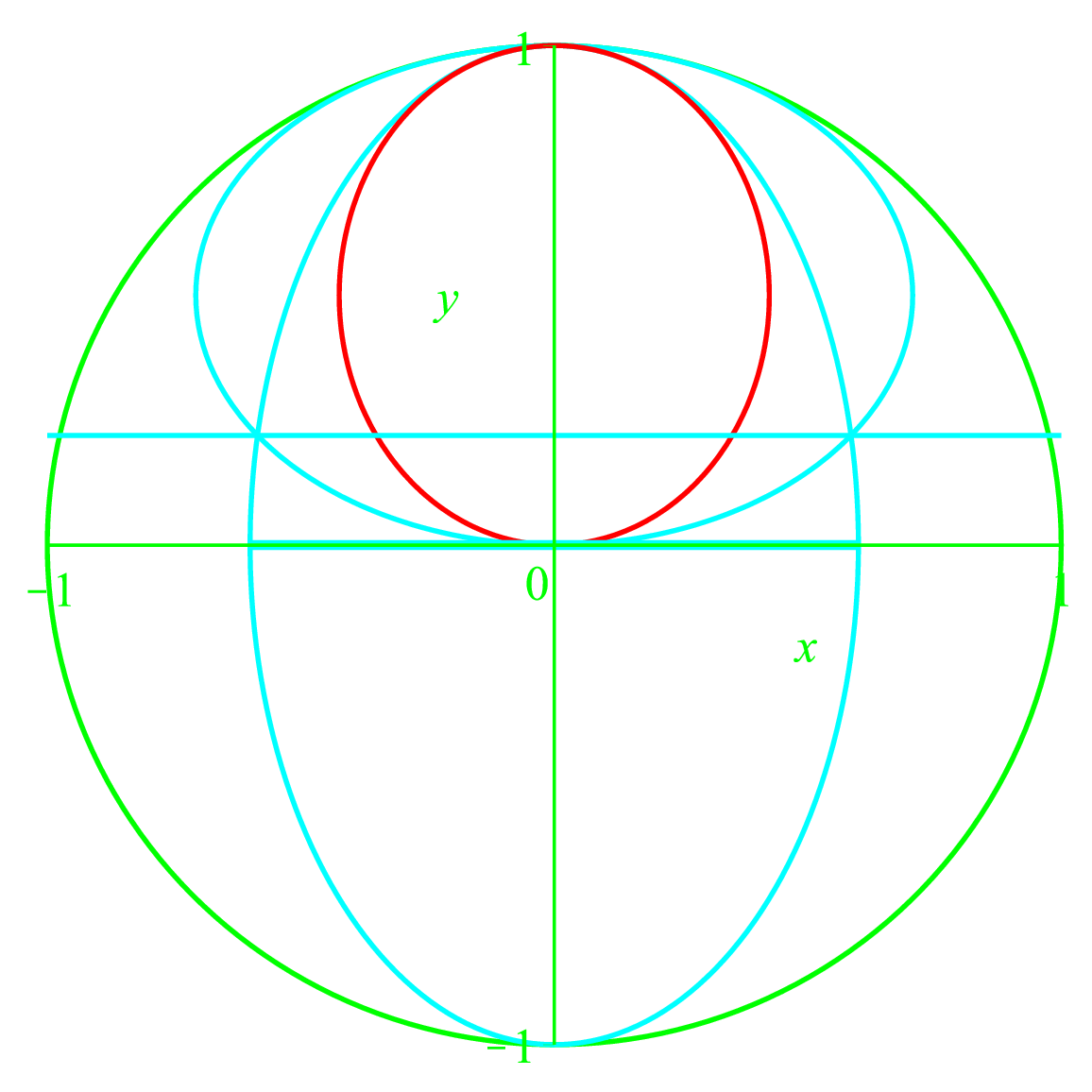}
    \caption*{\ref{fig:figHEP03}(a)  $\partial E^C_{\phantom{/}}$ [red].}
  \end{subfigure}
  \begin{subfigure}[b]{2.5in}
    \includegraphics[width=2.5in]{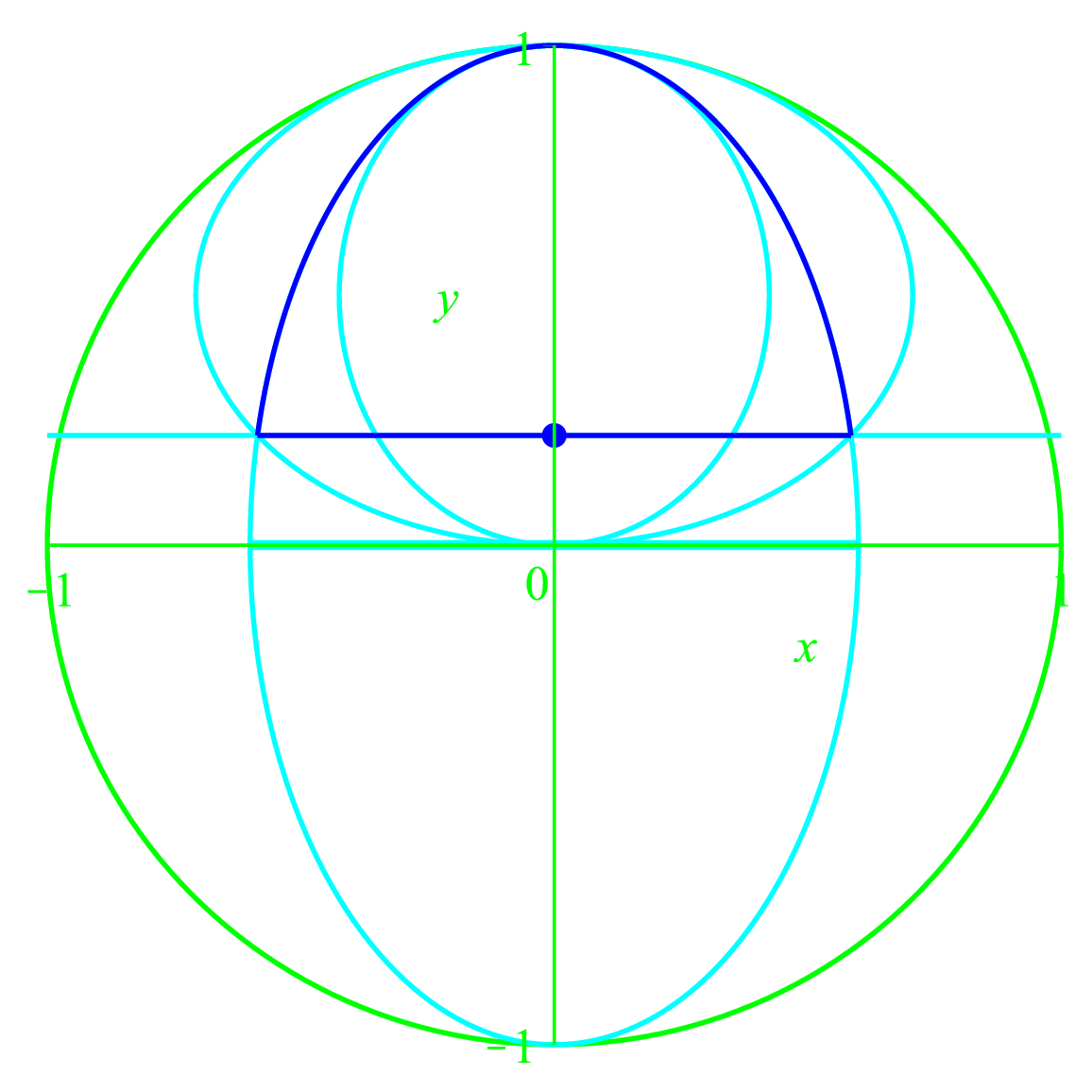}
    \caption*{\ref{fig:figHEP03}(b)  $\partial B^C_{\mathrm{up}\,\artanh C^2/(2-C^2)}$ [blue].}
  \end{subfigure}
  \\

   \begin{subfigure}[b]{2.5in}
    \includegraphics[width=2.5in]{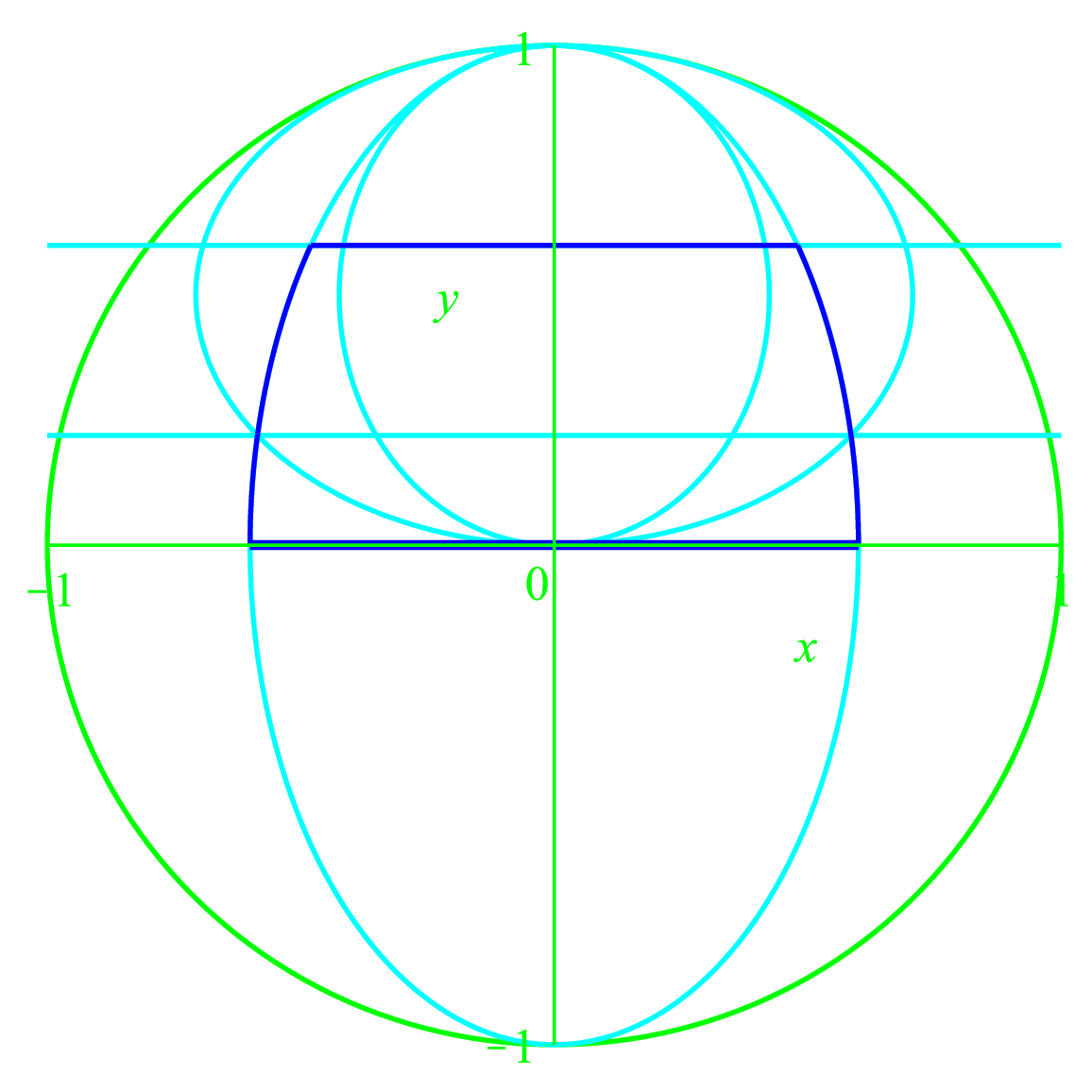}
    \caption*{ \ref{fig:figHEP03}(c)  $\partial B^C_{3/5}$ [blue].}
  \end{subfigure}
  \begin{subfigure}[b]{2.5in}
    \includegraphics[width=2.5in]{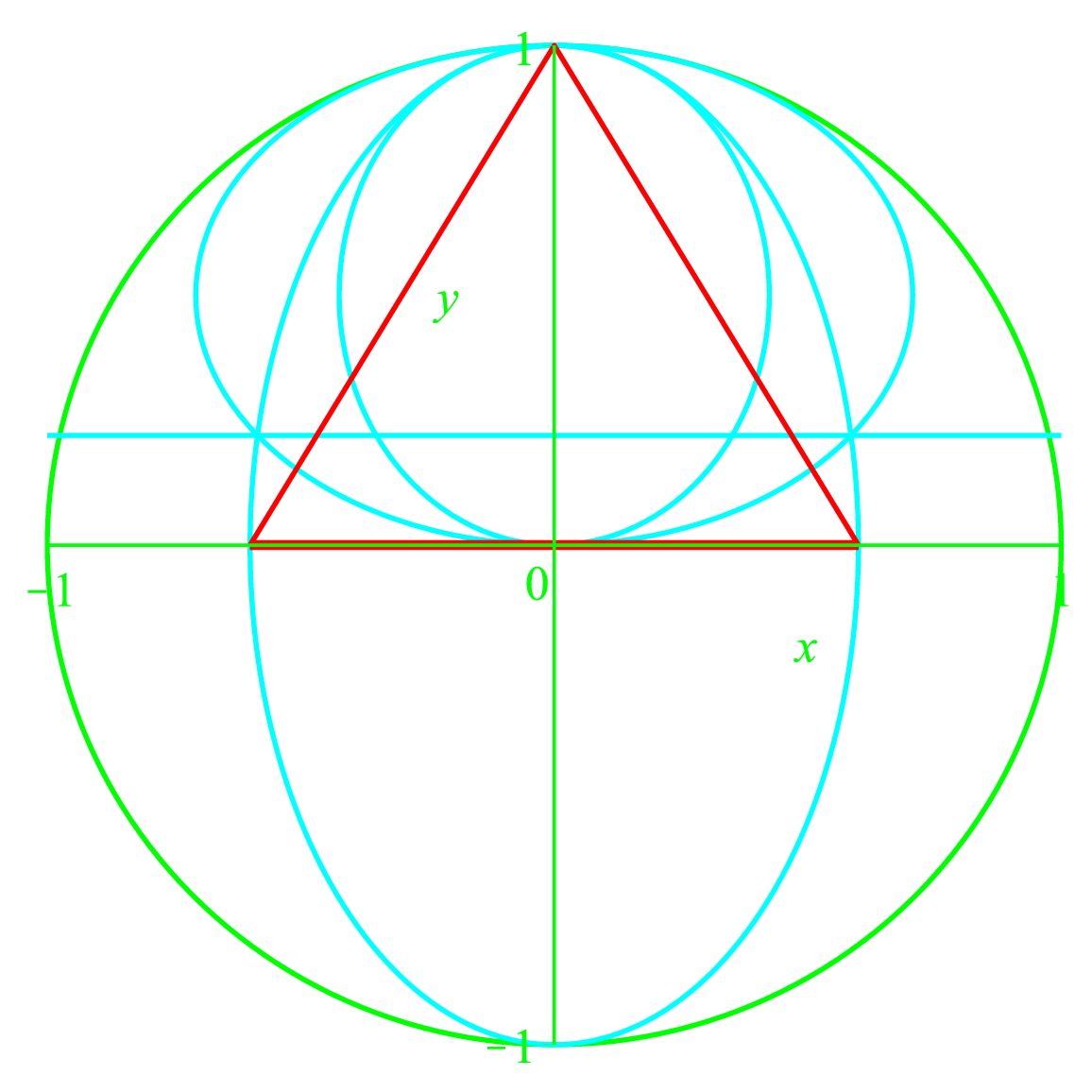}
    \caption*{\ref{fig:figHEP03}(d)  $\partial A^C_{\phantom{/}}$  [red].}
  \end{subfigure}

   \caption*{Fig.~\ref{fig:figHEP03} To
   $E^C=B^C_{\mathrm{up}\,\artanh  C^2/(2-C^2)}+ B^C_{3/5}-A^C$ in hyperbolic area.
   }
\phantomcaption
\plabel{fig:figHEP03}

   \begin{subfigure}[b]{2.5in}
    \includegraphics[width=2.5in]{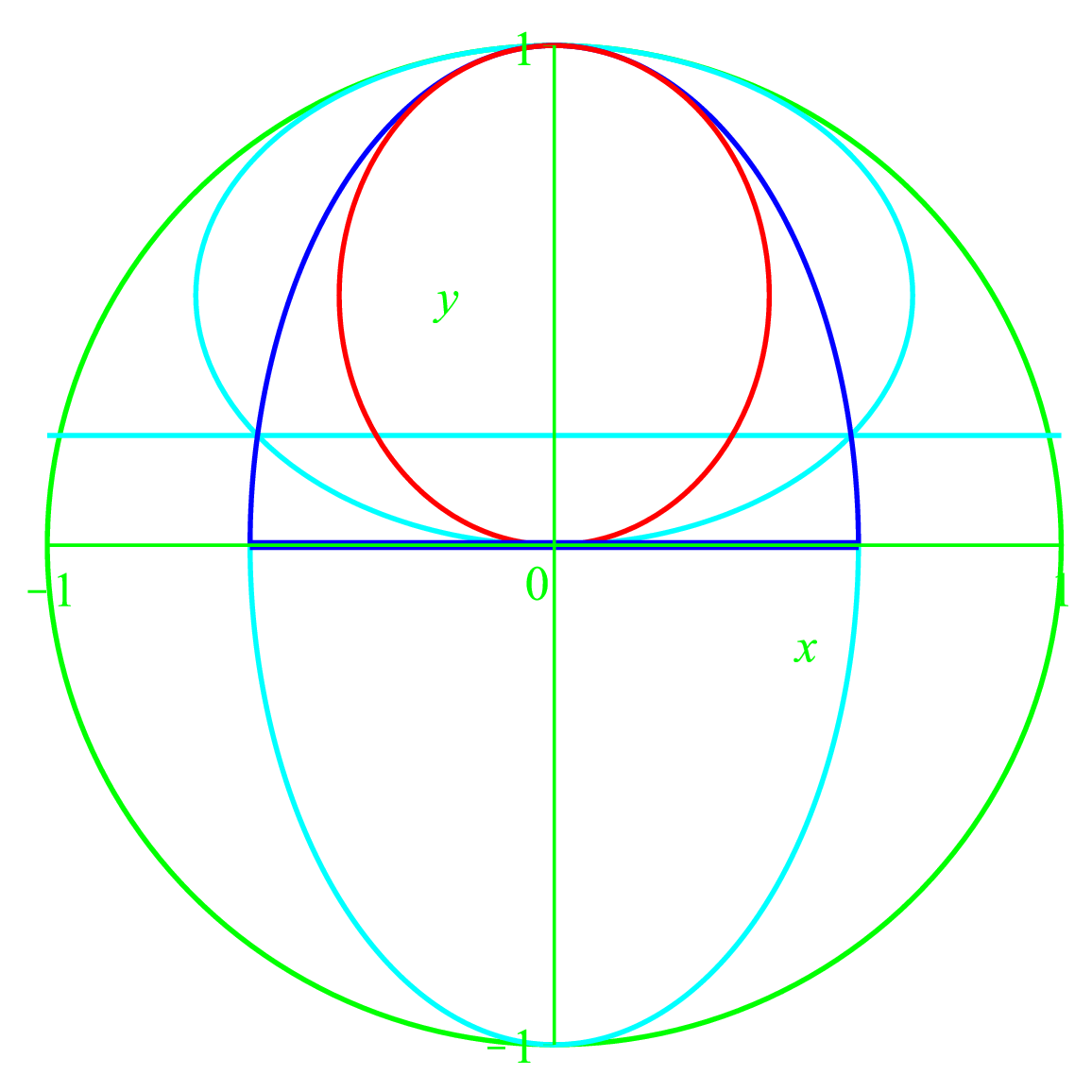}
    \caption*{\ref{fig:figHEP04}(a)  $\partial E^C$ [red], $\partial B^C$ [blue] }
  \end{subfigure}
  \begin{subfigure}[b]{2.5in}
    \includegraphics[width=2.5in]{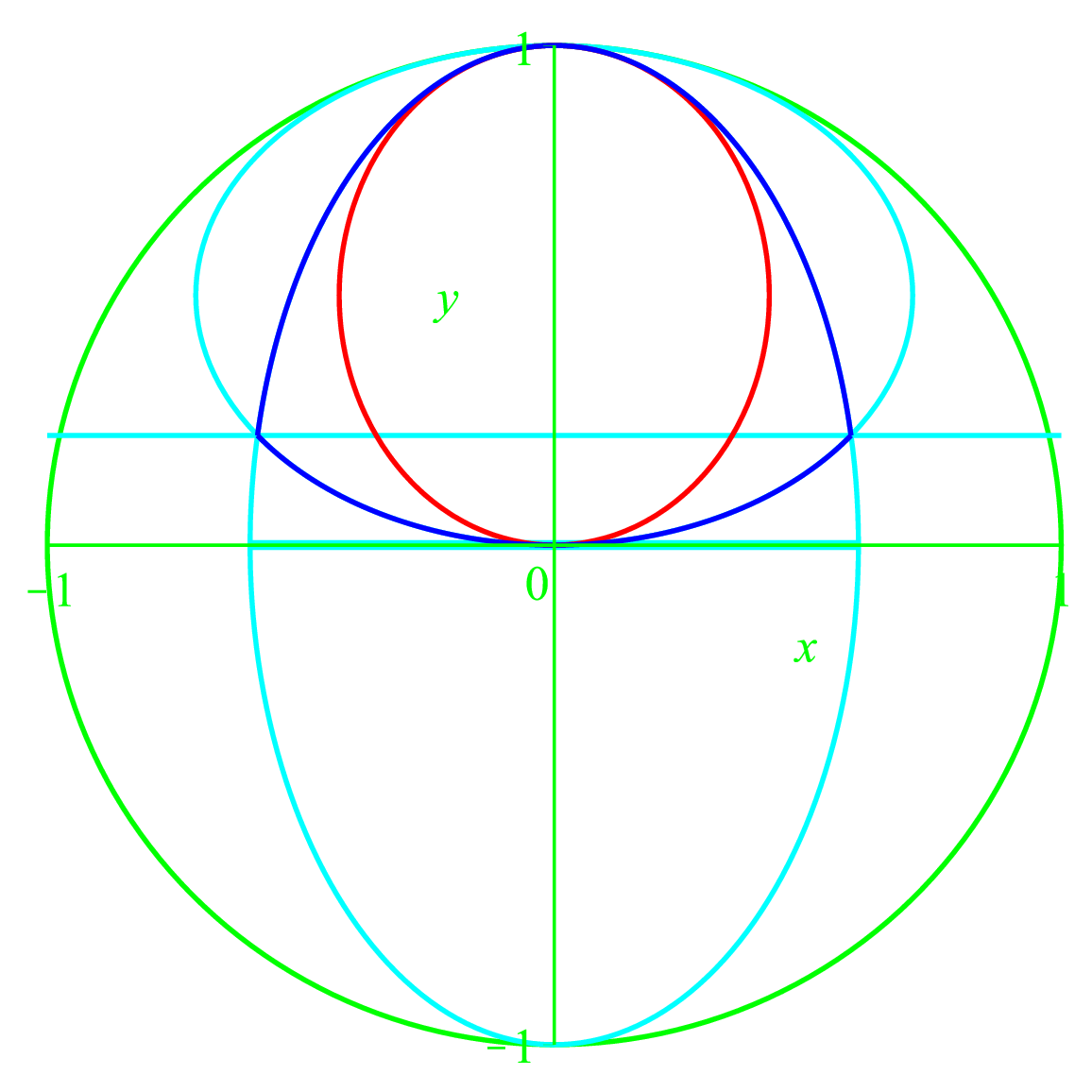}
    \caption*{\qquad\ref{fig:figHEP04}(b)  $\partial E^C$ [red], $\partial  D^C $ [blue].}
  \end{subfigure}
  \\
   \caption*{Fig.~\ref{fig:figHEP04} To $E^C=D'^C_{\mathrm{up}\,1-\ln 2 }\qquad\text{in area}$ .
   }
\phantomcaption
\plabel{fig:figHEP04}
\end{figure}

We  must note that the elements with finite area can be placed quite freely.
For example, $B^C_{\mathrm{up}\,\artanh \frac{C^2}{2-C^2}}+ B^C_{3/5}$ can be replaced by
 $ B^C_{  \mathrm{up}\,\artanh \frac{C^2}{2-C^2} -\ln 2}$; also, $A^C$
 can placed with its basis to be any of the perpendicular diameters of the supporting
 distance band, and with asymptotic direction upward or downward.
\snewpage

There is, however, an interpretation which is nicer than \eqref{eq:equ1}.
Firstly, we can extend the definition \eqref{eq:Edef} (and \eqref{eq:Edef2}) to the case of $C=1$.
It yields the supporting horodisk (with some cutoff at $\eta$).
Let
\[D^C=B^C\cap E^1,\]
 cf.~Figure \ref{fig:figHEP04}.
It must be obvious that $D^C$ approximates $E^C$ better than $B^C$.
The non-trivial question is whether this new approximation is nicer or not.

One can see that
\begin{align}
\Area_{\hyp}(D^C\setminus E^C)&=\Area_{\hyp}(B^C\setminus E^C)-\Area_{\hyp}(B^C\setminus D^C)\notag\\
&=\Area_{\hyp}(B^C\setminus E^C)-\Area_{\hyp}(B^C_{C^2/(2-C^2)}\setminus E^1_{C^2/(2-C^2)})\notag\\
&=\Area_{\hyp}(B^C\setminus E^C)-\Area_{\hyp}(B^C_{C^2/(2-C^2)})+ \Area_{\hyp}( E^1_{C^2/(2-C^2)}).\plabel{eq:triox}
\end{align}
\leaveout{
 \begin{figure}[htp]
   \begin{subfigure}[b]{2.5in}
    \includegraphics[width=2.5in]{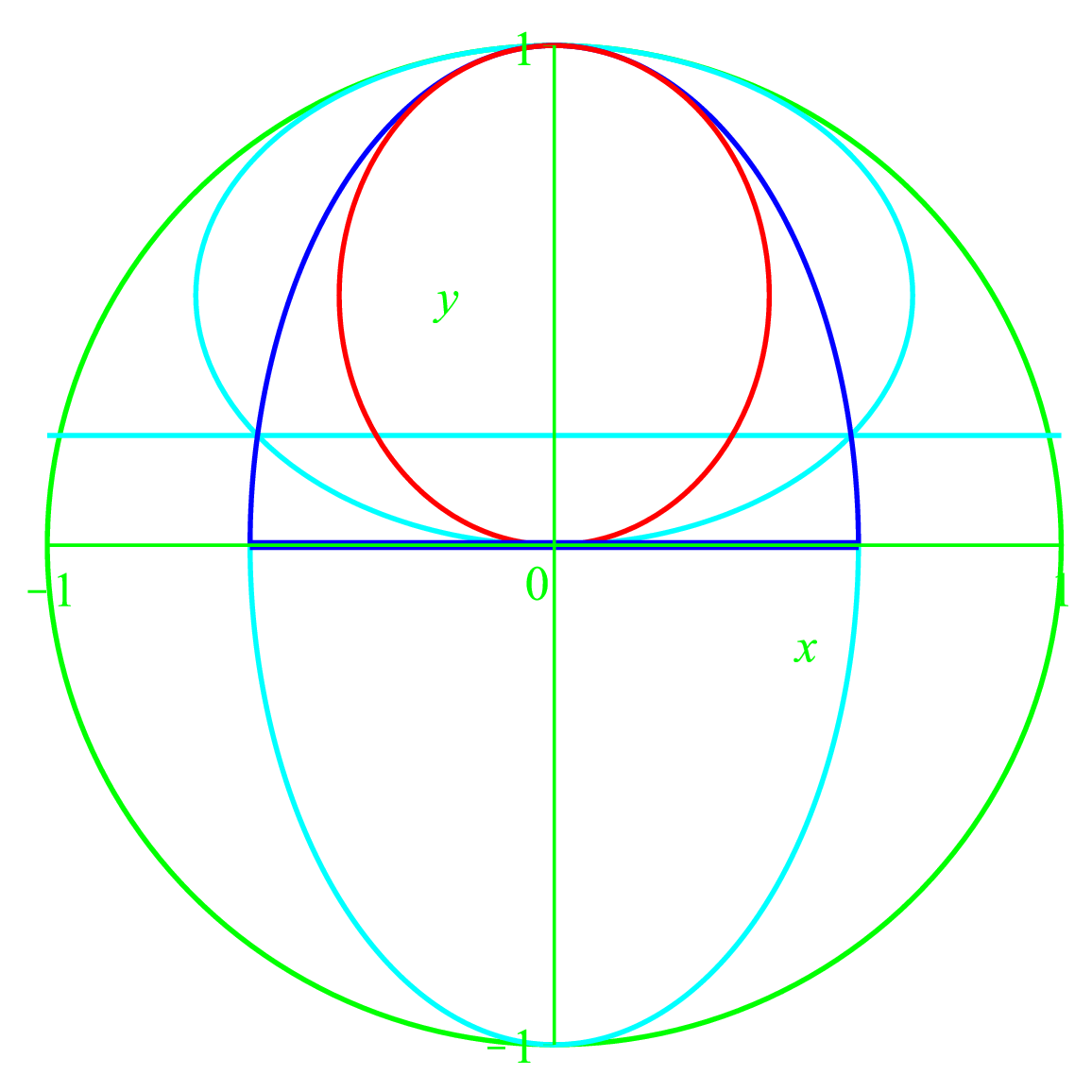}
    \caption*{\ref{fig:figHEP04}(a)  $\partial E^C$ [red], $\partial B^C$ [blue] }
  \end{subfigure}
  \begin{subfigure}[b]{2.5in}
    \includegraphics[width=2.5in]{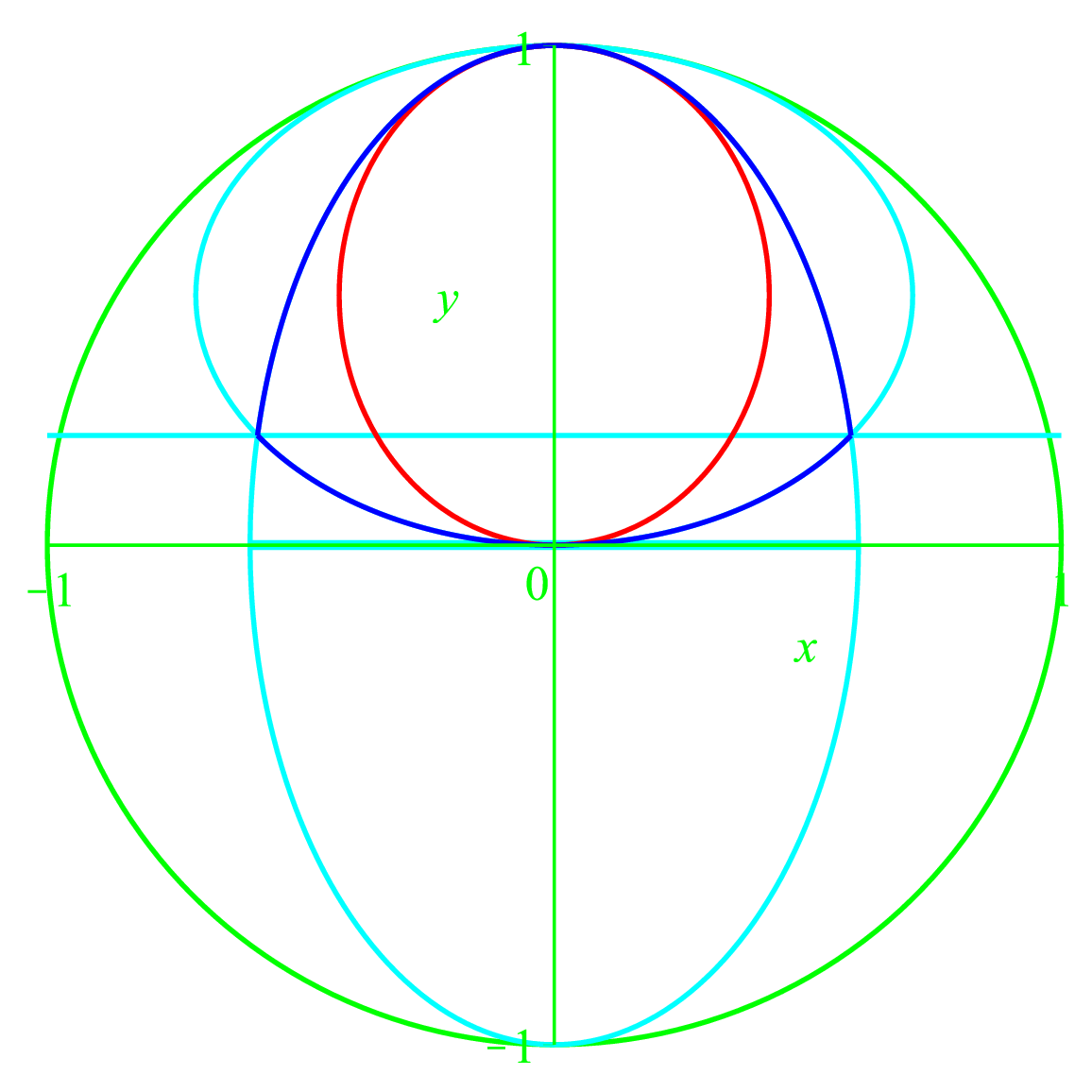}
    \caption*{\qquad\ref{fig:figHEP04}(b)  $\partial E^C$ [red], $\partial  D^C $ [blue].}
  \end{subfigure}
  \\
   \caption*{Fig.~\ref{fig:figHEP04} To $E^C=D'^C_{\mathrm{up}\,1-\ln 2 }\qquad\text{in area}$ .
   }
\phantomcaption
\plabel{fig:figHEP04}
\end{figure}
}

We have seen that $\Area_{\hyp}(B^C_{C^2/(2-C^2)})=\frac{2C}{\sqrt{1-C^2}}\artanh \frac{C^2}{2-C^2}$.
We can also compute
\begin{align*}
\Area_{\hyp}(E^1_\eta)&=\int_{y=0}^\eta\int_{x=- \sqrt{2y(1-y)}}^{ \sqrt{2y(1-y)}} \frac1{(1-x^2-y^2)^{3/2}}\,\mathrm dx\,\mathrm dy
\\&=\int_{y=0}^\eta \left[\frac{x}{(1-y^2)\sqrt{1-x^2-y^2}}\right]_{x=- \sqrt{2y(1-y)}}^{ \sqrt{2y(1-y)}}\,\mathrm dy
\\&=\int_{y=0}^\eta\frac{2 \sqrt{2y }}{ (1-y^2)\sqrt{1  - y }}\,\mathrm dy
\\&=2
{\frac { \sqrt {2\eta} }{\sqrt { 1-\eta }
}}
-2\arctan {\frac {  \sqrt {2\eta}}{\sqrt { 1 -\eta }}} .
\end{align*}
\begin{commentx}
\[=2{\frac { \sqrt {2\eta} }{\sqrt { 1-\eta }}}-2\arcsin {\frac {  \sqrt {2\eta}}{\sqrt { 1 +\eta }}}
=2{\frac { \sqrt {2\eta} }{\sqrt { 1-\eta }}}-2\arcsin {\frac {  \sqrt {1-\eta}}{\sqrt { 1 +\eta }}}
 .\]
\end{commentx}
The computation above was sort of unnecessary again.
Firstly, we could have used   $\Area_{\hyp}(E^1_\eta)=\lim_{C\nearrow1}\Area_{\hyp}(E^C_\eta)$ utilizing our
 earlier (more complicated) computation.
Secondly, the area of the horocircular segment was already known to Lobachevsky and Bolyai.
(It is the difference of a horocircular sector and a symmetric triange asymptotic in a vertex).
In any case, that yields,
\[\Area_{\hyp}(E^1_{C^2/(2-C^2)})=2\frac{C}{\sqrt{1-C^2}}-2\arctan\frac{C}{\sqrt{1-C^2}}.\]

Comparing \eqref{eq:finar}, and \eqref{eq:triox}, and the areas we have just considered, we find,
\begin{equation}
\Area_{\hyp}(D^C \setminus E^C) =\frac{2C}{\sqrt{1-C^2}}(1-\ln 2).
\plabel{eq:finarp}
\end{equation}
Let $D^C_{\mathrm{up}\,\omega}$ be defined as $D ^C $ but translated by hyperbolic distance $\omega$ upward along the $y$-axis.
Then the area difference between  $D ^C $ and  $ D^C_{\mathrm{up}\,\omega}$ is $\frac{2C}{\sqrt{1-C^2}}\omega$
 (as we can think about their relation by taking a cut at $y=\eta\gg0$ and translating the lower part, and then eliminating the
 doubled distance band segment).
This implies that
\begin{equation}
E^C=D^C_{\mathrm{up}\, 1-\ln 2 }\qquad\text{in area}.
\plabel{eq:equ2}
\end{equation}
It is a surprisingly simple relationship that $E^C$ has the same area as $D^C$ translated by hyperbolic distance $ 1-\ln 2 $ upward.
It is notable, because the translation distance does not depend on $C$.
\snewpage

 \section{Circumference}\plabel{sec:cf}

\textbf{The computation of the circumference difference.}
Comparing the circumference of $E^C$ and $B^C$ is theoretically a little more complicated than the area,
 where the area $B^C\setminus E^C$ could have been computed in a direct manner.
For the circumference some kind of particular interpretation is needed.
We will simply use regularization by with cutoffs $y\leq\eta$ as $\eta\nearrow1$.
Such cut-offs were already used in the computation of area but there they were for the convenience of computation only.
Here the role of the cut-offs will be conceptual, signalling one particular approach but one which is quite reasonable.

Let us explain our regularization process in detail.
In terms of the boundary, let
\[\partial^\cup E^C_\eta=\partial E^C_\eta\cap\{(x,y)\,:\, y\leq\eta\},\]
and
\[\partial^\cup B^C_\eta=\partial B^C_\eta\cap\{(x,y)\,:\, y\leq\eta\},\]
 moreover, let
\[M^C_\eta=[-C\sqrt{1-\eta^2},-C\sqrt{2\eta(1-\eta)}]\times\{\eta\}
\cup[C\sqrt{2\eta(1-\eta)},C\sqrt{1-\eta^2}]\times\{\eta\},\]
 cf.~Figure \ref{fig:figHEP05}.
Then, regarding the arc length $\Len_{\hyp}$, we find
\[\Len_{\hyp}\left(\partial^\cup B^C_\eta\right)-\Len_{\hyp}\left(\partial^\cup E^C_\eta\right)+\Len_{\hyp}(M^C_\eta)=
\Len_{\hyp}\left(\partial  B^C_\eta\right)-\Len_{\hyp}\left(\partial  E^C_\eta\right)\geq 0.\]
(The last inequality is valid because we deal with bounded convex domains, one contained in the other.)
Now, as long as $\lim_{\eta\nearrow1}\Len_{\hyp}(M^C_\eta)=0$ and
\begin{equation}
\lim_{\eta\nearrow1}\Len_{\hyp}\left(\partial^\cup B^C_\eta\right)-\Len_{\hyp}\left(\partial^\cup E^C_\eta\right) =G^C
\plabel{eq:arclim}
\end{equation}
 exists, we can confidently say that the difference of the circumference of $B^C$ and $E^C$ is  $G^C$.
In what follows, we carry out these computations.
 \begin{figure}[ht]
   \begin{subfigure}[b]{2.5in}
    \includegraphics[width=2.5in]{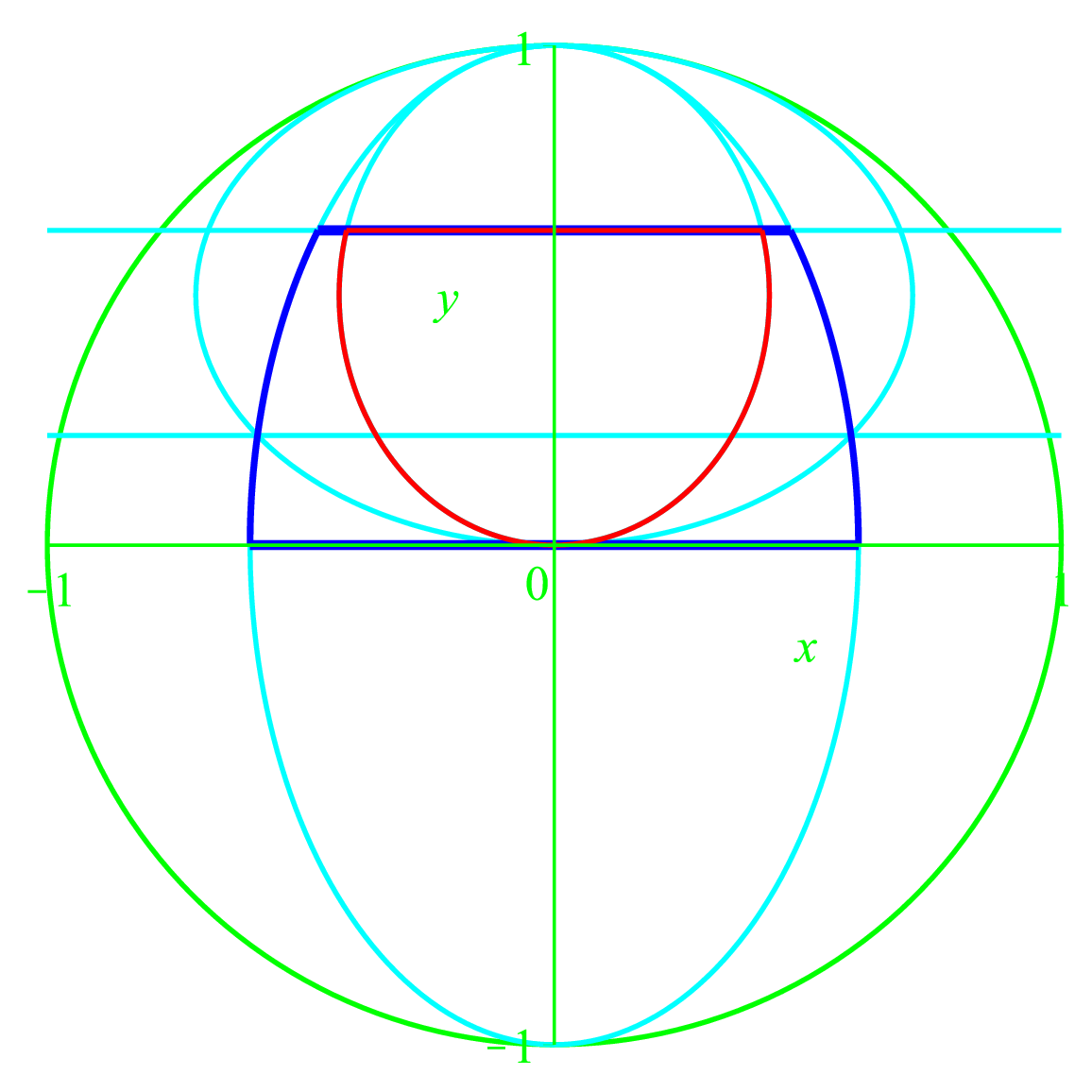}
    \caption*{\ref{fig:figHEP05}(a)  $\partial E^C_\eta$ [red], $\partial B^C_\eta$ [blue].\\ \phantom{$M^C_\eta$}}
  \end{subfigure}
  \begin{subfigure}[b]{2.5in}
    \includegraphics[width=2.5in]{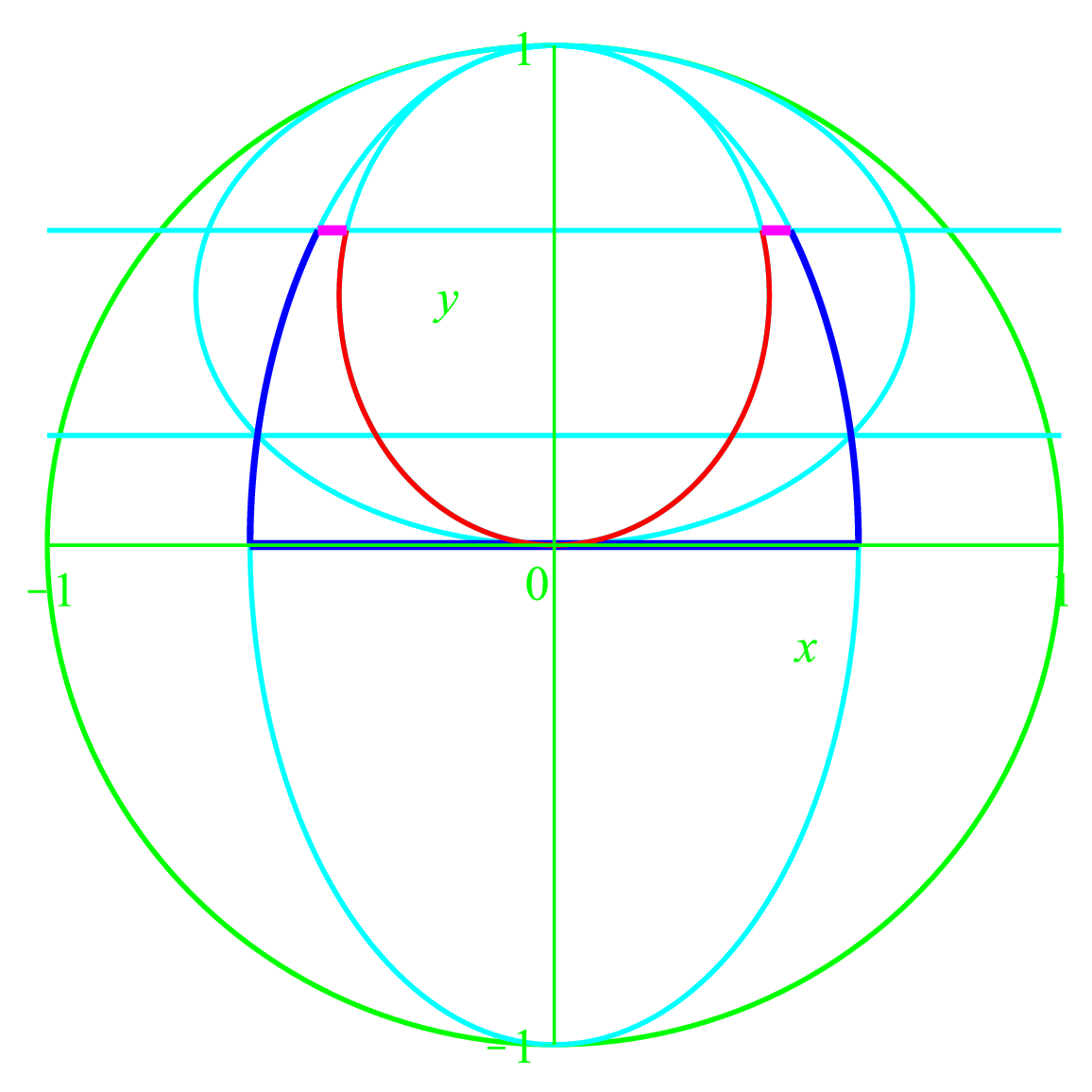}
    \caption*{\qquad\ref{fig:figHEP05}(b)  $\partial^\cup E^C_\eta$ [red], $\partial^\cup B^C_\eta$ [blue],\\
    \phantom{art}\qquad\qquad  $M^C_\eta$ [magenta]. }
  \end{subfigure}
  \\
   \caption*{Fig.~\ref{fig:figHEP05} Lineal cut-offs.}
\phantomcaption
\plabel{fig:figHEP05}

\end{figure}
\snewpage

Firstly, let us note that
\begin{align*}
\Len_{\hyp}\left(M^C_\eta\right)&=2\dist\left((C\sqrt{1-\eta^2},\eta),(C\sqrt{2\eta(1-\eta)},\eta) \right)
\\&=2\arcosh{\frac {\sqrt {\eta+1}-{C}^{2} \sqrt {2\eta}}{\sqrt {1-C^2 }\sqrt {1+\eta-2\,{C}^{2}\eta}}}
\\&=2\artanh {\frac { \left( \sqrt {\eta+1}- \sqrt {2\eta} \right) C}{\sqrt {\eta+1}-{C}^{2}\sqrt {2\eta}}}.
\end{align*}
\begin{commentx}
\[=2\arsinh{\frac { \left( \sqrt {\eta+1}- \sqrt {2\eta} \right) C}{\sqrt{1-C^2 }\sqrt{1+\eta-2\,{C}^{2}\eta}}}
  =2\ln{\frac { (1+C)\left( \sqrt {\eta+1}- C\sqrt {2\eta} \right)}{\sqrt{1-C^2 }\sqrt{1+\eta-2\,{C}^{2}\eta}}}\]
\end{commentx}
This implies that
\[\lim_{\eta\nearrow1}\Len_{\hyp}\left(M^C_\eta\right)=0,\]
 indeed.
Therefore we can proceed by examining \eqref{eq:arclim}.

It is well-known that the arc length element (squared) in the BCK model is given by
\[\text{``}(\mathrm ds_{\hyp})^2\text{''}=\frac{(1-y^2)(\mathrm dx)^2+2xy\,\mathrm dx\,\mathrm dy+(1-x^2)(\mathrm dy)^2}{(1-x^2-y^2)^2}.\]
I. e., the hyperbolic length of a $C^1$  arc $\gamma:(a,b)\rightarrow \{(x,y)\,:\,x^2+y^2<1\}$ is
\begin{align*}
\Len_{\hyp}(\gamma)&=
\int_{t=a}^b \sqrt{\frac{(1-( \gamma(t)_2)^2)(\dot \gamma(t)_1)^2+2  \gamma(t)_1   \gamma(t)_2 \dot \gamma(t)_1 \dot \gamma(t)_2
+(1- ( \gamma(t)_1)^2)(\dot \gamma(t)_2)^2}{(1-(  \gamma(t)_1)^2-(  \gamma(t)_2)^2)^2}}\,\mathrm dt
\\&=\int_{t=a}^b\frac{\sqrt{(1-|\gamma(t)|^2)|\dot\gamma(t)|^2+(\gamma(t)\dot\gamma(t))^2 }}{1-|\gamma(t)|^2}\,\mathrm dt.
\end{align*}

Applying this to map $\gamma_{B^C,\eta}: t\in(0,\eta)\mapsto (C\sqrt{1-t^2},t)$,  we find
\[\Len_{\hyp}(\gamma_{B^C,\eta})=\int_{t=0}^\eta\,\frac1{\sqrt{1-C^2}}\cdot\frac1{1-t^2}\,\mathrm dt=\frac1{\sqrt{1-C^2}}\artanh \eta. \]
Again, using the arc length formula was not really necessary.
This arc length is naturally proportional to the hyperbolic length the segment $\{0\}\times[0,\eta]$, which is $\artanh\eta$,
 and finding out the coefficient to the axial length is not hard.
In, fact it was already known by Lobachevsky and Bolyai that the arc length of a (one-sided) distance line is
  $\cosh d$ times the axial  length, where $d$ is characteristic distance (radius) of the distance line.
(And this is the actual cornerstone of their analytic investigations.)
In the present situation $d=\artanh C$, so $\cosh d=\cosh \artanh C=\frac1{\sqrt{1-C^2}}$, indeed.

Now, it follows that
\begin{equation}
\Len_{\hyp}\left(\partial^\cup B^C_\eta\right)=2\frac1{\sqrt{1-C^2}}\artanh\eta+2\artanh C.\plabel{eq:bmon}
\end{equation}
\begin{commentx}
\[=2\frac1{\sqrt{1-C^2}}\arsinh\frac{\eta}{\sqrt{1-\eta^2}}+2\arsinh \frac{C}{\sqrt{1-C^2}}\]
\[=2\frac1{\sqrt{1-C^2}}\arcosh\frac{1}{\sqrt{1-\eta^2}}+2\arcosh \frac{1}{\sqrt{1-C^2}}\]
\[=2\frac1{\sqrt{1-C^2}}\ln\frac{\sqrt{1+\eta}}{\sqrt{1-\eta}}+2\ln \frac{\sqrt{1+C}}{\sqrt{1-C}}\]
\end{commentx}
%\snewpage

Applying the arc length formula to the map $\gamma_{E^C,\eta}: t\in(0,\eta)\mapsto (C\sqrt{2t(1-t)},t)$,  we find
\begin{align}\notag
\Len_{\hyp}\left(\partial^\cup E^C_\eta\right)&=2\Len_{\hyp}(\gamma_{E^C,\eta} )
\notag\\&=2\int_{t=0}^\eta \,{\frac {\sqrt{{C}^{2}+2\,t-3\,{C}^{2}t}}{ \left( 1-t \right)  \left(1+t-2\,{C}^{2}t \,\right) \sqrt{2t}} }
 \,\mathrm dt
\notag\\&=2\left(\frac1{\sqrt{1-C^2}}\artanh   \sqrt {{\frac {2\eta   \left( 1-C^2 \right) }{{C}^{2}+2\,\eta-3\,{C}^{2}\eta}}}
 -\artanh   \sqrt {{\frac {2\eta   \left( 1-C^2 \right)^2 }{{C}^{2}+2\,\eta-3\,{C}^{2}\eta}}}\right).
\plabel{eq:Ecal}
\end{align}
\begin{commentx}
\[=2\left(\frac1{\sqrt{1-C^2}}\arsinh\left(\frac{\sqrt{1-C^2}}{C}\sqrt{\frac{2\eta}{1-\eta}}\right)
 -\arsinh\left(\frac{\sqrt{1-C^2}}{C}\sqrt{\frac{2\eta(1-C^2)}{1+\eta-2C^2\eta}}\right)\right)\]
\[=2\left(\frac1{\sqrt{1-C^2}}\arcosh  \sqrt  {\frac {{C}^{2}+2\,\eta-3\,{C}^{2}\eta}{ C^2   \left( 1- \eta\right) }}
 -\arcosh   \sqrt {{ \frac {{C}^{2}+2\,\eta-3\,{C}^{2}\eta} {C(1+\eta-2C^2\eta) } }}\right)\]
\[=2\left(\frac1{\sqrt{1-C^2}}\ln\frac {\sqrt{{C}^{2}+2\,\eta-3\,{C}^{2}\eta}
+\sqrt{1-C^2}\sqrt{2\eta}
}{ \sqrt  {C^2   \left( 1- \eta\right) }}
-\artanh   \sqrt {{\frac {2\eta   \left( 1-C^2 \right)^2 }{{C}^{2}+2\,\eta-3\,{C}^{2}\eta}}}\right)\]
\end{commentx}
Rewriting $\frac1{\sqrt{1-C^2}}\artanh\ldots$ into $\frac1{\sqrt{1-C^2}}\log\ldots$, similarly as before, we find
\begin{multline*}
\Len_{\hyp}\left(\partial^\cup B^C_\eta\right)-\Len_{\hyp}\left(\partial^\cup E^C_\eta\right)=
%\\=
2 \frac1{\sqrt{1-C^2}}\ln\frac{ \sqrt  {C^2   \left( 1+ \eta\right) }} {\sqrt{{C}^{2}+2\,\eta-3\,{C}^{2}\eta}+\sqrt{1-C^2}\sqrt{2\eta}}
\\
+2\artanh   \sqrt {{\frac {2\eta   \left( 1-C^2 \right)^2 }{{C}^{2}+2\,\eta-3\,{C}^{2}\eta}}}+2\artanh C  .
\end{multline*}
\begin{commentx}
\begin{multline*}
=2\frac{1}{\sqrt{1-C^2}}\artanh
\frac{{C}^{2}\eta-\left(  1+\eta \right)\sqrt {2\eta \left( 1-C^2 \right)   }\sqrt { {C}^{2}+2\,\eta-3\,{C}^{2}\eta}}
{{C}^{2}+2\,\eta+2\,{\eta}^{2}-2\,{C}^{2}\eta -2\,{C}^{2}{\eta}^{2}}
\\
+2\artanh  \sqrt {{\frac {2\eta   \left( 1-C^2 \right)^2 }{{C}^{2}+2\,\eta-3\,{C}^{2}\eta}}}+2\artanh  C.
\end{multline*}
\begin{multline*}
=2\frac{1}{\sqrt{1-C^2}}\arsinh {\frac { \left( 3\,{C}^{2}\eta+2\,{C}^{2}-2\,\eta-2 \right) \eta}{C\sqrt {1+\eta}
\, \left( \sqrt {2}\sqrt { \eta \left( 1-C^2\right)  }+
\eta\sqrt {2\,\eta+{C}^{2}-3\,{C}^{2}\eta} \right) }}
\\
+2\arsinh\left(\frac{\sqrt{1-C^2}}{C}\sqrt{\frac{2\eta(1-C^2)}{1+\eta-2C^2\eta}}\right)+2\arsinh \frac{C}{\sqrt{1-C^2}}.
\end{multline*}
\end{commentx}

Consequently,
\begin{align}
G^C&\equiv\lim_{\eta\nearrow1}\Len_{\hyp}\left(\partial^\cup B^C_\eta\right)-\Len_{\hyp}\left(\partial^\cup E^C_\eta\right)
 \plabel{eq:lofi}
 \\&=2\frac1{\sqrt{1-C^2}}\ln\frac{C}{2\sqrt{1-C^2}}
 +2\artanh\sqrt{1-C^2} +2\artanh C.
 \notag
\end{align}
\begin{commentx}
\[=2\frac1{\sqrt{1-C^2}}\artanh\frac{5C^2-4}{4-3C^2}
 +2\artanh\sqrt{1-C^2} +2\artanh C\]
\[ =2\frac1{\sqrt{1-C^2}}\arsinh\frac{5C^2-4}{2C\sqrt{1-C^2}}
 +2\arsinh \frac{\sqrt{1-C^2}}{C} +2\arsinh \frac{C}{\sqrt{1-C^2}}\]
\end{commentx}
This is finally the formula for ``$\Len_{\hyp}\left(\partial^\cup B^C\right)-\Len_{\hyp}\left(\partial^\cup E^C\right)$''.
We remark that
\[\ln\frac{C}{2\sqrt{1-C^2}}=\artanh\frac{C^2}{2-C^2}+\ln C-\ln 2.\]
\begin{commentx}
We remark that
\begin{align*}
\artanh\frac{5C^2-4}{4-3C^2}&=\artanh (2C^2-1)-\artanh\frac35\\
&\equiv\left(\artanh\frac{C^2}{2-C^2}+\ln C\right)-\ln 2.
\end{align*}
\end{commentx}

Using standard methods elementary analysis, one can show that
$G^C$ is strictly increasing in $C\in(0,1)$ and
$\lim_{C\searrow0}G^C=0$ and $\lim_{C\nearrow1}G^C=+\infty$.
%\snewpage

\textbf{The interpretation of the circumference difference.}
Again, we may look after $B^C_{\mathrm{up}\,\beta(C)}$ such that
``$\Len_{\hyp}\left(\partial^\cup B^C_{\mathrm{up}\,\beta(C)}\right)-\Len_{\hyp}\left(\partial^\cup E^C\right)$''$=0$.
As, in general,
``$\Len_{\hyp}\left(\partial^\cup B^C\right)-\Len_{\hyp}\left(\partial^\cup B^C_{\mathrm{up}\,\omega}\right)$''$
=2\frac1{\sqrt{1-C^2}}\omega$, one has
\[\beta(C)= \ln\frac{C}{2\sqrt{1-C^2}} + \sqrt{1-C^2}\left(\artanh\sqrt{1-C^2} +\artanh C\right).\]
One can show that
$ \beta(C)$ is strictly increasing in $C\in(0,1)$ and
$\lim_{C\searrow0} \beta(C)=0$ and $\lim_{C\nearrow1} \beta(C)=+\infty$.
But $\lim_{C\nearrow1} \beta(C)-\artanh\frac{C^2}{2-C^2}=-\ln 2$.

One can, predictably, approximate better by $D^C$.
One can see that
``$\Len_{\hyp}\left(\partial  B^C\right)-\Len_{\hyp}\left(\partial D^C\right)$''$
= \Len_{\hyp}\left(\partial^\cup B^C_{C^2/(2-C^2)}\right)-\Len_{\hyp}\left(\partial^\cup E^1_{C^2/(2-C^2)}\right)$.
From \eqref{eq:bmon},
\[\Len_{\hyp}\left(\partial^\cup B^C_{C^2/(2-C^2)} \right)=2\frac1{\sqrt{1-C^2}}\artanh \frac{C^2}{2-C^2} +2\artanh C.\]
It is easy to compute that
\begin{align}
\Len_{\hyp}\left(\partial^\cup E^1_\eta\right)&=2\Len_{\hyp}(\gamma_{E^1,\eta} )\plabel{eq:emon}
\\&=2\int_{t=0}^\eta \,{\frac {1}{ \left( 1-t \right)^2}\sqrt { {\frac {1-t}{2t}}}}\,\mathrm dt\notag
\\&=2\left[\sqrt{\frac{2t}{1-t}}\right]_{t=0}^\eta=2\sqrt{\frac{2\eta}{1-\eta}}.\notag
\end{align}
But this also follows from taking \eqref{eq:Ecal} to the limit $C\nearrow1$.
And also, it was already known to  Lobachevsky and Bolyai that
 to a horocyclic arc of chord length $2r$ $(=2\artanh\sqrt{\frac{2\eta}{1+\eta}})$ the corresponding arc length
 is $2\sinh r$ $(=2 \sqrt{\frac{2\eta}{1-\eta}})$.
In any case, this yields
\[\Len_{\hyp}\left(\partial^\cup E^1_{C^2/(2-C^2)}\right)=2 \frac{C}{\sqrt{1-C^2}}.\]

Then we have
\begin{align}
G'^{C}&\equiv\text{``$\Len_{\hyp}\left(\partial  D^C\right)-\Len_{\hyp}\left(\partial E^C\right)$''}\plabel{eq:Gprime}
\\&=\text{``$\Len_{\hyp}\left(\partial  B^C\right)-\Len_{\hyp}\left(\partial  E^C\right)$''}\notag
\\&\quad -\Len_{\hyp}\left(\partial^\cup B^C_{C^2/(2-C^2)}\right)+\Len_{\hyp}\left(\partial^\cup E^1_{C^2/(2-C^2)}\right)\notag
\\&=2\frac{\ln C-\ln2}{\sqrt{1-C^2}}+2\artanh\sqrt{1-C^2}+2 \frac{C}{\sqrt{1-C^2}}.\notag
\end{align}
Again, $G'^C$ is strictly increasing in $C\in(0,1)$ and
$\lim_{C\searrow0}G'^C=0$ and $\lim_{C\nearrow1}G'^C=+\infty$.

Now, we look after $D^C_{\mathrm{up}\,\widehat{\beta}(C)}$ such that
``$\Len_{\hyp}\left(\partial  D^C_{\mathrm{up}\,\widehat{\beta}(C)}\right)-\Len_{\hyp}\left(\partial  E^C\right)$''$=0$.
As, again,
``$\Len_{\hyp}\left(\partial  D^C\right)-\Len_{\hyp}\left(\partial  D^C_{\mathrm{up}\,\omega}\right)$''$
=2\frac1{\sqrt{1-C^2}}\omega$, one has
\[\widehat{\beta}(C)=\ln C-\ln 2 +\sqrt{1-C^2} \artanh\sqrt{1-C^2} + C .\]
One can show that
$ \widehat{\beta}(C)$ is strictly increasing in $C\in(0,1)$ and
$\lim_{C\searrow0} \widehat{\beta}(C)=0$ and $\lim_{C\nearrow1} \widehat{\beta}(C)=1-\ln2$.
The situation is not so nice as in the case of the area, but there is at least some similarity.
Unfortunately, we have not managed to find a practically simple replacement object for  $\partial E^C$.

Yet, as a main point, we have managed to find answers regarding the circumference difference between $E^C$ and $B^C$.

There might remain, however, at this point, one nagging question in this matter.
We have managed to find the circumference difference only using a particular type of regularization.
We may wonder, whether there is an other method which has the same stability as the area difference computation.
Before providing a particular answer, we need some preliminaries:
\snewpage
\section{Background: The Study--de Sitter plane geometry}\plabel{sec:sds}
As hyperbolic geometry is the geometry of the interior of the unit circle in the ambient projective space of the BCK model,
 the Study--de Sitter geometry is the geometry of the exterior of the unit circle in the ambient projective space of the BCK model.
Considering the plane, in terms of differential geometry, we let
\[\text{``}(\mathrm ds_{\SdS})^2\text{''}=-\frac{(1-y^2)(\mathrm dx)^2+2xy\,\mathrm dx\,\mathrm dy+(1-x^2)(\mathrm dy)^2}{(1-x^2-y^2)^2},\]
 but one should note that this metric becomes indefinite (Lorentzian) for $x^2+y^2>1$.
The associated area element is
\[\mathrm{da}_{\SdS}(x,y)=\frac1{(x^2+y^2-1)^{3/2}}\,|\mathrm dx\wedge \mathrm dy|.\]
(In the formulas above, we just changed the signs of the natural hyperbolic expressions
 in order not to compute with negative arc lengths and volumes later.)
One should note that the formulae above do not describe the geometry quite completely.
The reason is that the Study--de Sitter geometry is the exterior of unit circle in the ambient projective plane,
 therefore alternative coordinatization is needed around the points which are ideal in Euclidean view.
We leave this to the reader, or see the literature quoted below.

The  Study--de Sitter geometry was considered first by Study \cite{Stu} (1907) (cf.~A'Campo-Neuen,  Papadopoulos \cite{ACP}).
One should note  that the exterior points were also geometrized for metric computational purposes
 by Vörös \cite{V1} (1909) (cf.~G. Horváth \cite{GH12}).
The  customary naming `de Sitter spaces' is after the work of de Sitter \cite{dS} (1917),
 who considered these spaces in terms of general relativity (cf.~Coxeter \cite{Cxt}).
For a modern review, the reader is directed to Fillastre, Seppi \cite{FS}.

Let us assume that $K$ is a non-empty compact convex set in the BCK model of hyperbolic geometry (so it is in the interior of the unit circle).
The co-polar (complementary polar) set $\widetilde{K}$ is given as follows.
Let us consider all those lines in the ambient projective space which intersect $K$.
Then let us take their pole with respect to the absolute (that is the unit circle), the set of these poles is $\widetilde{K}$.
Then its boundary  $\partial\widetilde{K}$ corresponds to the polars of the tangent lines of $K$.
%(The word ``dual'' is also used; but in my viewpoint the dual lives in the dual projective space, while the
% polar lives in the original one).
Then $\widetilde{K}$ always contains a projective line containing only exterior points of the absolute
 (e.~g., the polar of any point of $K$), and this line is a homotopy retract of $\widetilde{K}$.

The relevant points for us are
\begin{equation}
\Area_{\hyp}(K)=-2\pi+\Len_{\SdS}(\partial\widetilde{K})
\plabel{eq:pc1}
\end{equation}
and
\begin{equation}
\Len_{\hyp}(\partial K)=\Area_{\SdS}(\widetilde{K}).
\plabel{eq:pc2}
\end{equation}
I.~e., the computation of area and circumference are dual to each other with respect to the polar.
The equalities \eqref{eq:pc1} and \eqref{eq:pc2} are in analogy to spherical geometry.
Relative to spherical geometry, apart from the sign conventions,
 what leads to the asymmetry in \eqref{eq:pc2} is that we consider not the polar here but he
 complementary polar (in order to deal with less mysteries in the computation of the area).
\snewpage

We illustrate  \eqref{eq:pc1} and \eqref{eq:pc2} in the case when $K=K^C$ is the origin centered disk of hyperbolic radius $R\equiv\artanh C$.
The boundary of the disk is given by the equation $x^2+y^2-C^2=0$ and the boundary of the co-polar set is $x^2+y^2-1/C^2=0$.
Then
\begin{align*}
	\Area_{\hyp}(K^C)&=\int_{x^2+y^2\leq C^2}\frac1{(1-x^2-y^2)^{3/2}}\,\mathrm dx\,\mathrm dy
	%\\&
    =\int_{\theta=0}^{2\pi}\int_{r=0}^C\frac1{(1-r^2)^{3/2}}\,r\,\mathrm dr\,\mathrm d\theta
	\\&=\int_{\theta=0}^{2\pi}\left[\frac1{\sqrt{ 1-r^2} }\right]_{r=0}^C \,\mathrm d\theta
	   =\int_{\theta=0}^{2\pi}\left(\frac1{\sqrt{ 1-C^2} }-1\right)\,\mathrm d\theta
	\\&=2\pi\left(\frac1{\sqrt{ 1-C^2} }-1 \right)\equiv2\pi\left(\cosh R -1 \right).
\end{align*}
(This area was already known to Lobachevsky and Bolyai, of course.)
The boundary of the co-polar is parametrized by $\gamma_{K,1/C}\colon t\in[0,\pi]\rightarrow(\frac 1C\cos t,\frac 1C\sin t))$.
Then
\[\Len_{\SdS}(\partial\widetilde K^{C})=
\Len_{\SdS}( \gamma_{K,1/C})=\int_{t=0}^{2\pi} \sqrt{\frac1{1-C^2 }}\, \mathrm dt
=2\pi\frac1{\sqrt{1-C^2}}\equiv 2\pi \cosh R.\]
This is in accordance to \eqref{eq:pc1}.
Similarly,
\[\Len_{\hyp}(\partial  K^{C})=
\Len_{\hyp}( \gamma_{K,C})=\int_{t=0}^{2\pi} \sqrt{\frac{C^2}{1-C^2 }}\, \mathrm dt
=2\pi\frac{C}{\sqrt{1-C^2}}\equiv 2\pi \sinh R.\]
(Again, this circumference was already known to Lobachevsky and Bolyai.)
Here,
\begin{align*}
	\Area_{\SdS}(\widetilde K^C)&=\int_{x^2+y^2\geq 1/C^2}\frac1{(x^2+y^2-1)^{3/2}}\,\mathrm dx\,\mathrm dy
	%\\&
    =\int_{\theta=0}^{2\pi}\int_{r=1/C}^\infty\frac1{(r^2-1)^{3/2}}\,r\,\mathrm dr\,\mathrm d\theta
	\\&=\int_{\theta=0}^{2\pi}\left[-\frac1{\sqrt{ r^2-1} }\right]_{r=1/C}^\infty  \,\mathrm d\theta
	  =\int_{\theta=0}^{2\pi} \frac C{\sqrt{ 1-C^2} } \,\mathrm d\theta
	\\&=2\pi\frac C{\sqrt{ 1-C^2} }   \equiv2\pi \sinh R  .
\end{align*}
This is in accordance to \eqref{eq:pc2}.

Explaining the workings of  \eqref{eq:pc1}, \eqref{eq:pc2}    in complete generality,  is beyond the scope of this paper.
However, at this point, our objectives, in a certain way, are less mathematical than psychological:
 we want to reaffirm the arc length difference computation in an alternative way.
For this reason, we can leave more place to mathematical imagination than before.

\snewpage
\section{Circumference II}\plabel{sec:cf2}

At this point, it is clear that a good candidate for ``$\Len_{\hyp}\left(\partial^\cup B^C\right)-\Len_{\hyp}\left(\partial^\cup E^C\right)$'' is
 $\Area_{\SdS}(\widetilde B^C\setminus \widetilde E^C)$.
From the dual conic \eqref{eq:dnhepm}, the polar conic (relative to the absolute) is given by
\begin{multline*}
	\begin{bmatrix} x\\y\\1\end{bmatrix}^\top
	\begin{bmatrix} 1\\&1\\&&-1\end{bmatrix}
	\begin{bmatrix} \frac1{C^2}&0&0\\0&2&-1\\0&-1&0\end{bmatrix}^{-1}
	\begin{bmatrix} 1\\&1\\&&-1\end{bmatrix}	
	\begin{bmatrix} x\\y\\1\end{bmatrix}
	\equiv \\
	\equiv\begin{bmatrix} x\\y\\1\end{bmatrix}^\top
	\begin{bmatrix} C^2&0&0\\0&0&1\\0&1&-2\end{bmatrix}
	\begin{bmatrix} x\\y\\1\end{bmatrix}
	\equiv {C}^{2}{x}^{2}+2\,y-2
	\stackrel{\leftrightarrow}=0.
	\plabel{eq:pnhepm}
\end{multline*}
If the equation is taken literally, then it yields the implicit equation for  $\partial\widetilde E^C$, cf.~Figure \ref{fig:figHEP06}(a).
Computing  $\partial\widetilde B^C$ is more complicated:
Any tangent line to $B^C$ is either tangent to the hypercyclic  arc or passes through the vertices.
The dual of the hypercycle is given by
\begin{multline*}
	\begin{bmatrix} x\\y\\1\end{bmatrix}^\top
	\begin{bmatrix} 1\\&1\\&&-1\end{bmatrix}
	\begin{bmatrix}  \frac1{C^2}\\&1\\&&-1 \end{bmatrix}\begin{bmatrix} x\\y\\1\end{bmatrix}^{-1}
	\begin{bmatrix} 1\\&1\\&&-1\end{bmatrix}	
	\begin{bmatrix} x\\y\\1\end{bmatrix}
	\equiv \\\equiv
	\begin{bmatrix} x\\y\\1\end{bmatrix}^\top
	\begin{bmatrix} C^2&&\\ &1&\\&&-1\end{bmatrix}
	\begin{bmatrix} x\\y\\1\end{bmatrix}
	\equiv {C}^{2}{x}^{2}+y^2-1
	\stackrel{\leftrightarrow}=0.
%	\plabel{eq:pnhbpm}
\end{multline*}
The polars of the points $(\pm C,0)=[\pm C,0,1]$ are the lines $\pm Cx+y-1\stackrel{\leftrightarrow}=0$.
Then  $\partial\widetilde B^C$ is combined out of these polars as  in Figure \ref{fig:figHEP06}(a).
Then  $ \widetilde E^C$ and $ \widetilde B^C$ are the non-convex ``exteriors'' of
 $\partial\widetilde E^C$ and $\partial\widetilde B^C$, respectively.
In order to see the picture from a different viewpoint, we can consider a $h$-collineation $\Xi$
(like a downward translation) which sends an interior point of
 $E^C$ to the origin; then $\Xi$ extends to the ambient projective space, and then the images of
 $\partial\widetilde E^C$ and $\partial\widetilde B^C$ will avoid the ideal line, see Figure \ref{fig:figHEP06}(b).

%\snewpage
\begin{figure}[htp]
   \begin{subfigure}[b]{5.0in}
    \includegraphics[width=5.0in]{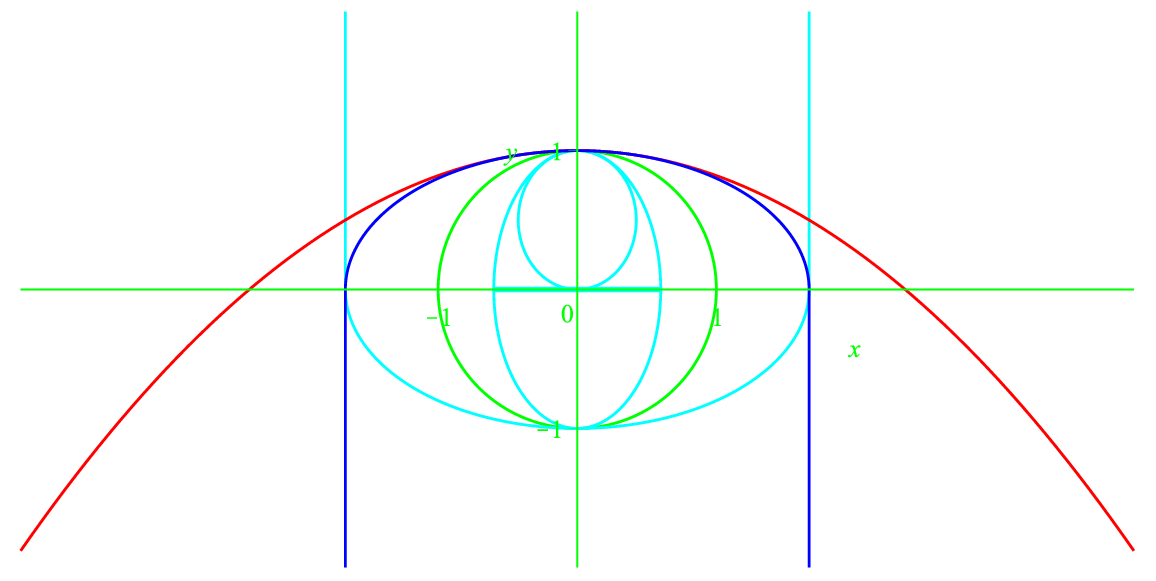}
    \caption*{\ref{fig:figHEP06}(a)  $\partial \widetilde{E}^C $ [red], $\partial \widetilde{B}^C$ [blue]}
  \end{subfigure}
   \begin{subfigure}[b]{2.5in}
    \includegraphics[width=2.5in]{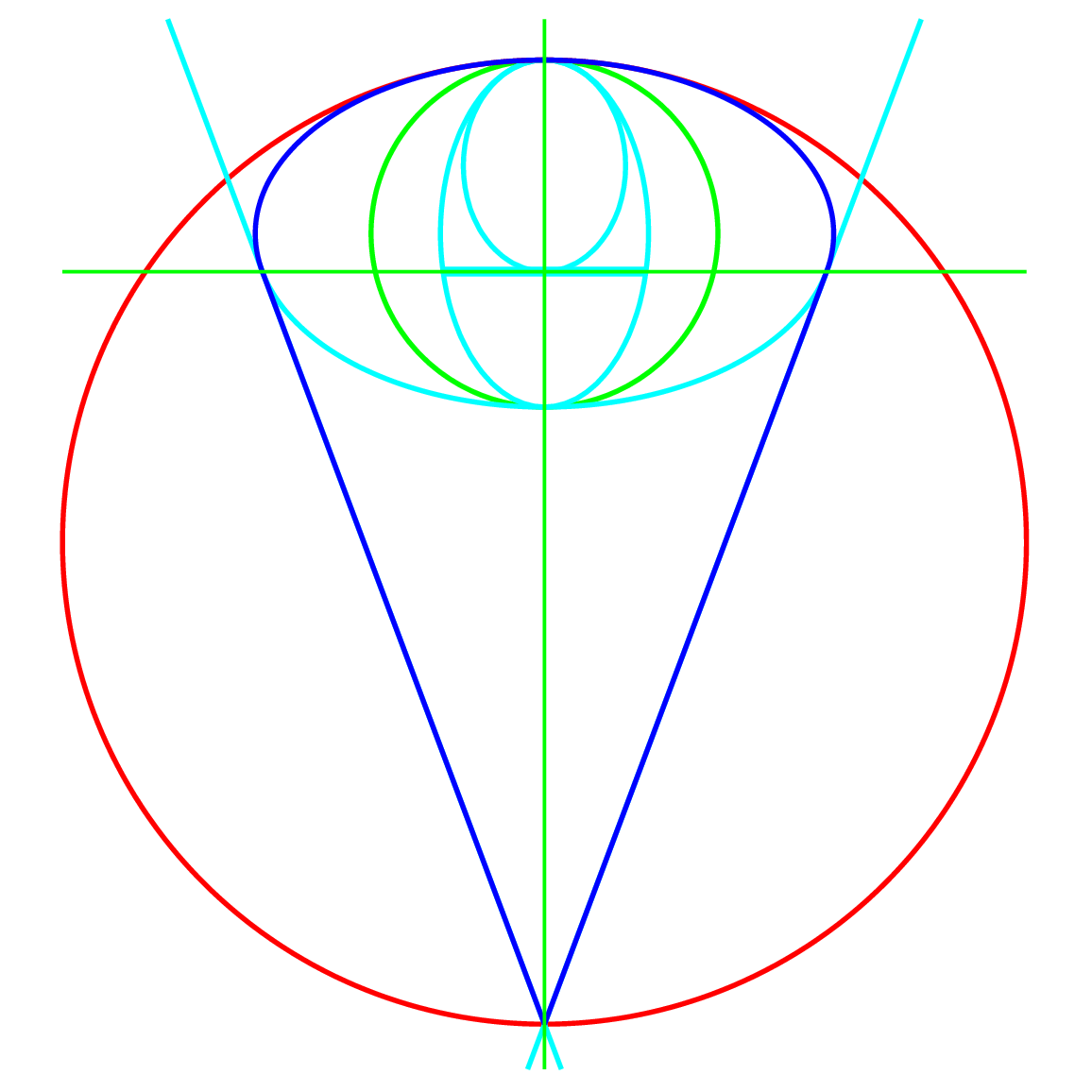}
    \caption*{\ref{fig:figHEP06}(b)  After a projective transformation.}
  \end{subfigure}
  \begin{subfigure}[b]{2.5in}
    \includegraphics[width=2.5in]{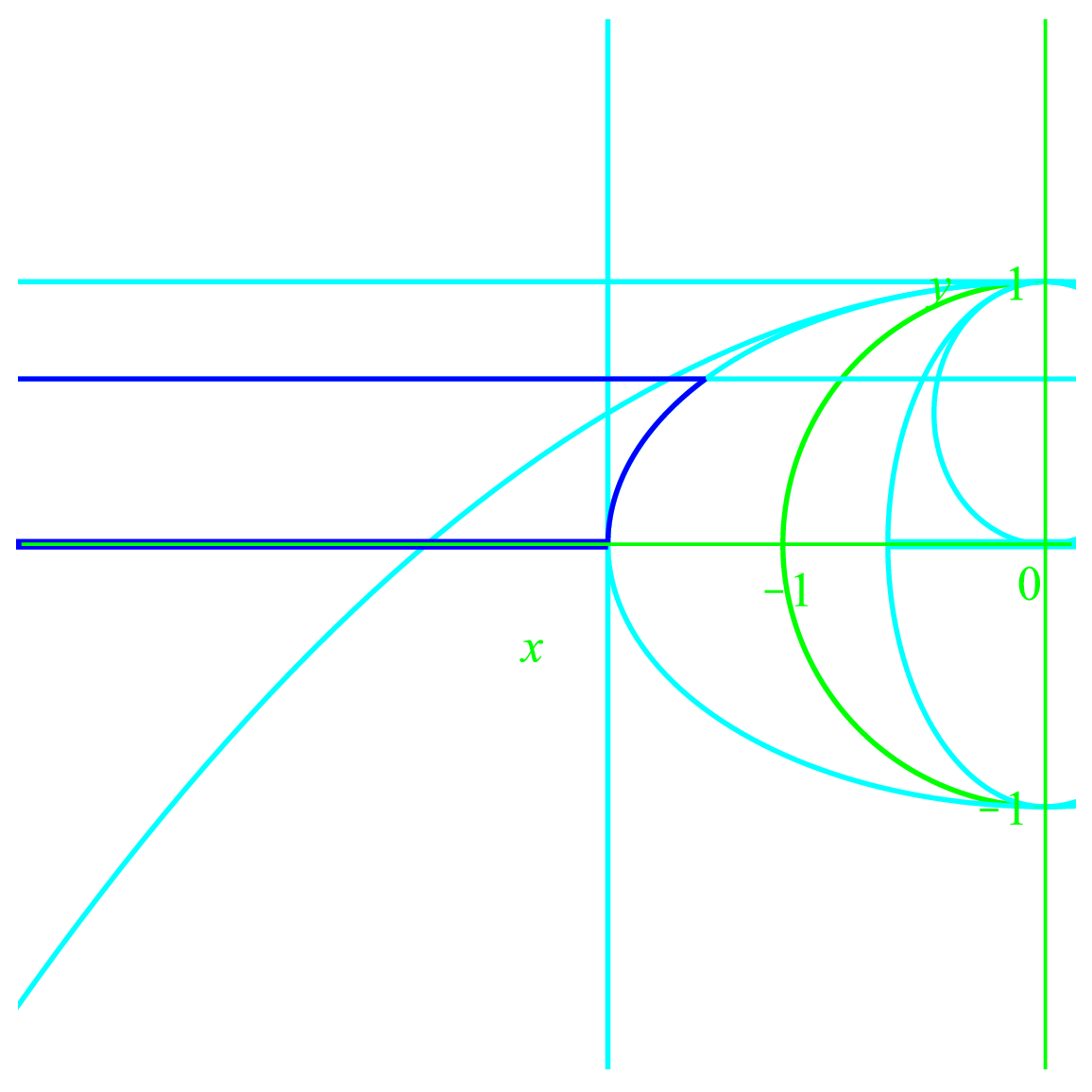}
    \caption*{\qquad\ref{fig:figHEP06}(c)    $\partial\,\nnn\widetilde{W}^C_\eta$ [blue].}
  \end{subfigure}
  \\
  \begin{subfigure}[b]{2.5in}
    \includegraphics[width=2.5in]{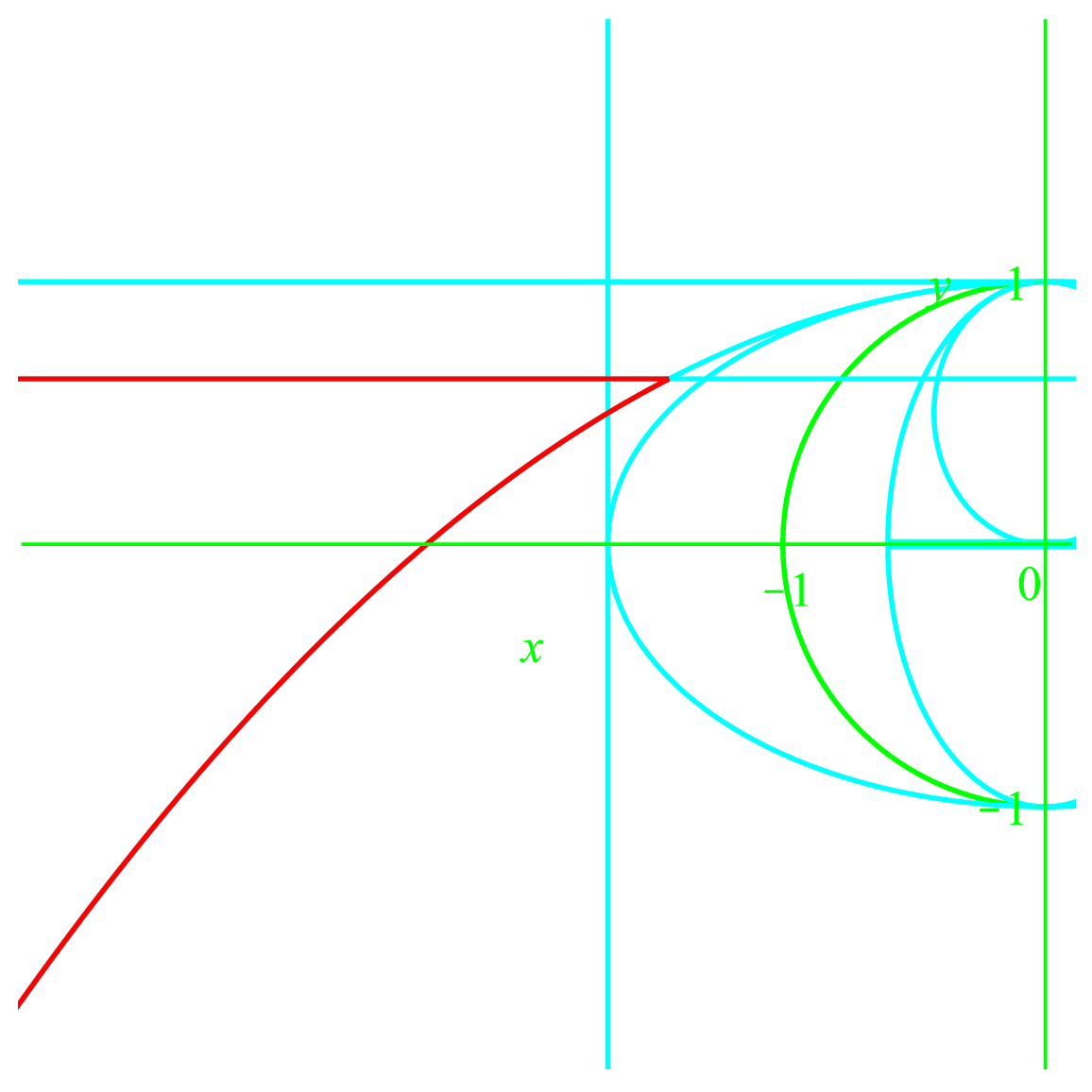}
    \caption*{\qquad\ref{fig:figHEP06}(d)    $\partial\,\nnn\widetilde{E}^C_\eta$ [red].}
  \end{subfigure}
  \begin{subfigure}[b]{2.5in}
    \includegraphics[width=2.5in]{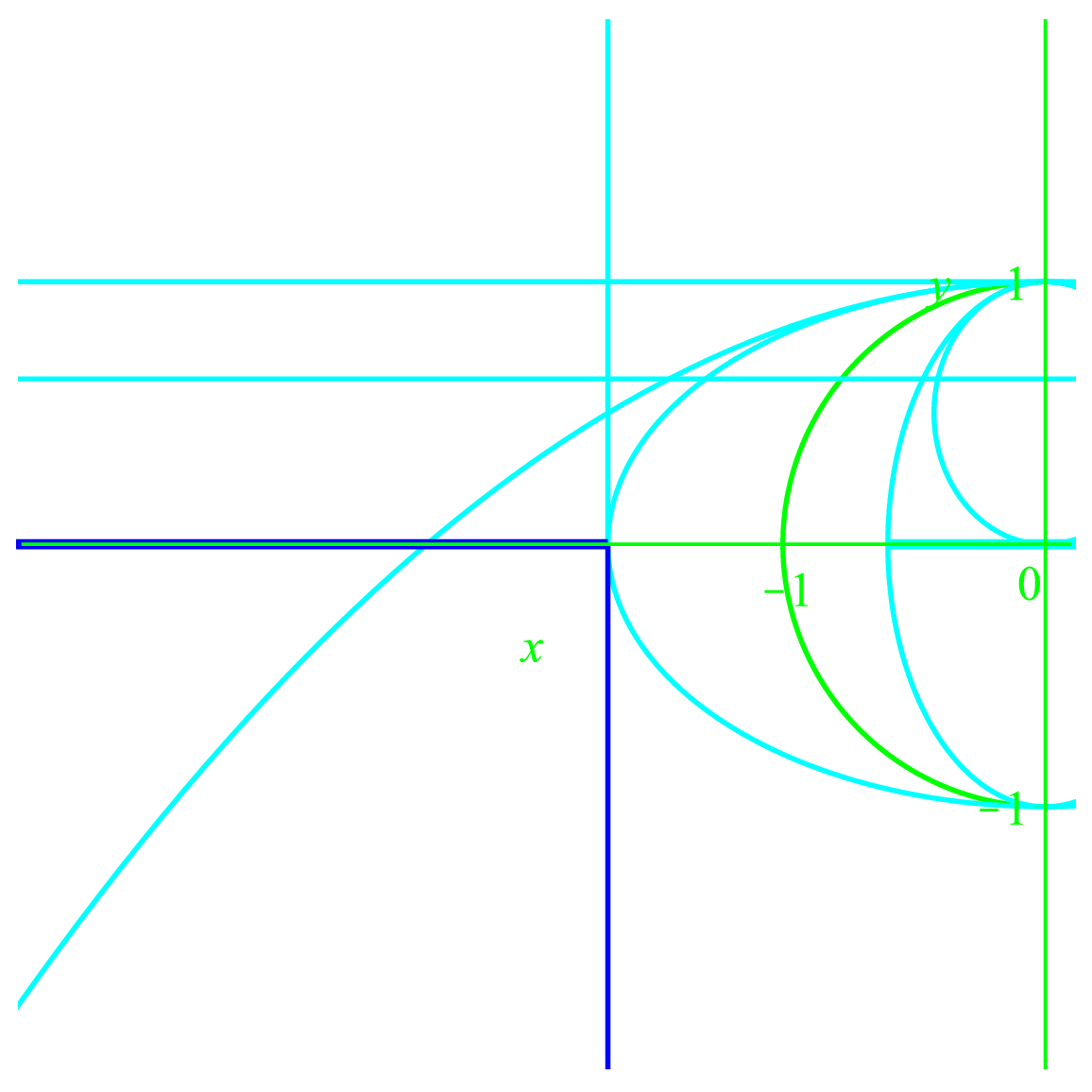}
    \caption*{\qquad\ref{fig:figHEP06}(e)    $\partial\,\nnn\widetilde{Z}^C$ [blue].}
  \end{subfigure}
  \\
   \caption*{Fig.~\ref{fig:figHEP06} The polar computation of the circumference difference
   }
\phantomcaption
\plabel{fig:figHEP06}
\end{figure}

In any case, our computational objective is to compute the SdS area between the red and blue curves in Figure \ref{fig:figHEP06}(a).
In order to make this computation, let us define, for $\eta\in(0,1)$, the auxiliary sets
\[\nnn\widetilde{W}^C_\eta=\{(x,y)\,:\, {C}^{2}{x}^{2}+y^2-1\geq0,\,0\leq y\leq\eta,\,x\leq 0\}, \]
 cf.~Figure \ref{fig:figHEP06}(c);
\[\nnn\widetilde{E}^C_\eta=\{(x,y)\,:\, {C}^{2}{x}^{2}+2\,y-2\leq0,\, y\leq\eta,\,x\leq 0  \}, \]
 cf.~Figure \ref{fig:figHEP06}(d);
\[\nnn\widetilde{Z}^C =\{(x,y)\,:\, x\leq -\tfrac1C,\,y\leq 0\}, \]
 cf.~Figure \ref{fig:figHEP06}(e).
Then $\ppp\widetilde{W}^C_\eta$, $\ppp\widetilde{E}^C_\eta$, $\ppp\widetilde{Z}^C $ can be defined similarly
 but reflected to the $y$-axis.
These sets are well away from the absolute, therefore their SdS area will be finite;
 and so is the SdS area of $(\widetilde B^C\setminus \widetilde E^C)\cap \{(x,y)\,:\, y\leq \eta\}$.
In fact,
\begin{align*}
\Area_\SdS(    (\widetilde B^C\setminus \widetilde E^C)\cap \{(x,y)\,:\, y\leq \eta\})
&=\Area_\SdS(\nnn\widetilde{W}^C_\eta )+\Area_\SdS(\nnn\widetilde{Z}^C )-\Area_\SdS(\nnn\widetilde{E}^C_\eta )\\
&\quad+\Area_\SdS(\ppp\widetilde{W}^C_\eta )+\Area_\SdS(\ppp\widetilde{Z}^C )-\Area_\SdS(\ppp\widetilde{E}^C_\eta )\\
&=2\Area_\SdS(\nnn\widetilde{W}^C_\eta )+2\Area_\SdS(\nnn\widetilde{Z}^C )-2\Area_\SdS(\nnn\widetilde{E}^C_\eta ).
\end{align*}
Let us now compute these areas.
\begin{align*}
\Area_\SdS(\nnn\widetilde{W}^C_\eta )
&=\int_{y=0}^\eta\int_{x=-\infty}^{-\frac1C\sqrt{1-y^2}}\frac1{(x^2+y^2-1)^{3/2}}\,\mathrm dx\,\mathrm dy
\\&=\int_{y=0}^\eta\left[\frac1{\sqrt{x^2+y^2-1}\left(\sqrt{x^2+y^2-1}-x\right)}\right]_{x=-\infty}^{-\frac1C\sqrt{1-y^2}}\,\mathrm dy
\\&=\int_{y=0}^\eta\frac{1-\sqrt{1-C^2}}{(1-y^2)\sqrt{1-C^2}}\,\mathrm dy
\\&=\frac{1-\sqrt{1-C^2}}{\sqrt{1-C^2}}\,\artanh \eta.
\end{align*}
Also,  and this actually extends to $\eta\in(-\infty,1)$,
\begin{align*}
\Area_\SdS(\nnn\widetilde{E}^C_\eta )
&=\int_{y=-\infty}^\eta\int_{x=-\infty}^{-\frac1C\sqrt{2(1-y)}}\frac1{(x^2+y^2-1)^{3/2}}\,\mathrm dx\,\mathrm dy
\\&=\int_{y=-\infty}^\eta\left[\frac1{\sqrt{x^2+y^2-1}\left(\sqrt{x^2+y^2-1}-x\right)}\right]_{x=-\infty}^{-\frac1C\sqrt{2(1-y)}}\,\mathrm dy
\\&=\int_{y=-\infty}^\eta\frac{C^2}{(1-y)2\sqrt{1-C^2(1+y)/2}(\sqrt{1-C^2(1+y)/2}+1)}\,\mathrm dy
\\&={\frac {1}{\sqrt {1-{C}^{2}}} \artanh   {\frac {\sqrt { 1-{C}^{2}}}{\sqrt {    1-C^2(1+\eta)/2}}} }
\\&\qquad-
 \artanh {\frac {2-C^2 + 2\sqrt {1- \,{C}^{2}(1+\eta)/2} }{ 2-C^2\eta+ 2\sqrt {1- \,{C}^{2}(1+\eta)/2} }}.
\end{align*}
\begin{commentx}
\[={\frac {1}{\sqrt{1-C^2}}\arsinh  {\frac { \sqrt{1-C^2}}{C\sqrt {(1-\eta)/2}}} }
-\arsinh   {\frac {2-{C}^{2}+2\sqrt{1-C^2(1-\eta)/2}}{2C\sqrt{(1-\eta)/2} \left(   1 +\sqrt{1-C^2(1-\eta)/2}\right) }} \]
 \[={\frac {1}{\sqrt {1-C^2}}\arcosh   {\frac {  \sqrt{1-C^2(1-\eta)/2}}{C\sqrt {(1-\eta)/2}}} }-
\arcosh   {\frac {2-{C}^{2}\eta+2\sqrt{1-C^2(1-\eta)/2}}{2C\sqrt {(1-\eta)/2} \left(   1+\sqrt{1-C^2(1-\eta)/2}
 \right) }}  \]
 \[={\frac {1}{\sqrt{1-C^2}}\ln    {\frac {  \sqrt{1-C^2(1-\eta)/2}+\sqrt {1-C^2} }{C \sqrt{(1-\eta)/2}}}   }
-\ln  {\frac {  1+\sqrt{1-C^2(1-\eta)/2}}{C\sqrt{(1-\eta)/2}}}  \]
\end{commentx}
However, for $\eta\in(-1,1)$,  and, in particular for $\eta\in(0,1)$, the case we are interested in,
\begin{equation}
\artanh {\frac {2-C^2+ 2\sqrt {1- \,{C}^{2}(1+\eta)/2} }{ 2-C^2\eta+ 2\sqrt {1- \,{C}^{2}(1+\eta)/2} }}=
\artanh  \sqrt {   1-C^2(1+\eta)/2} +\artanh\eta,
\plabel{eq:remo}
\end{equation}
 which is a standard application of the addition formula for $\artanh A+\artanh B=\artanh\frac{A+B}{1-AB}$;
 or comes via the transcription $\ln \frac{1+\sqrt {1- \,{C}^{2}(1+\eta)/2}}{C\sqrt{(1-\eta)/2}}$.
(But the RHS is not suited for use  in a primitive function in the integration above, due to domain problems.)
\snewpage

Furthermore,
\begin{align*}
\Area_\SdS(\nnn\widetilde{Z}^C )
&=\int_{x=-\infty}^{-\frac1C }\int_{y=-\infty}^0\frac1{(x^2+y^2-1)^{3/2}}\,\mathrm dy\,\mathrm dx
\\&=\int_{x=-\infty}^{-\frac1C }\left[   {\frac {y+\sqrt {{x}^{2}+{y}^{2}-1}}{\sqrt {{x}^{2}+{y}^{2}-1} \left(
{x}^{2}-1 \right) }}\right]_{y=-\infty}^0 \,\mathrm dx
 =\int_{x=-\infty}^{-\frac1C } \frac{1}{x^2-1}\,\mathrm dx
\\&=\left[-\arcoth x\right]_{x=-\infty}^{-\frac1C }=\artanh C.
\end{align*}

Consequently,
\begin{align*}
\Area_\SdS&(    (\widetilde B^C\setminus \widetilde E^C)\cap \{(x,y)\,:\, y\leq \eta\})=
\\&
=2\left(\Area_\SdS(\nnn\widetilde{W}^C_\eta )+\Area_\SdS(\nnn\widetilde{Z}^C )-\Area_\SdS(\nnn\widetilde{E}^C_\eta )\right)
\\&=2\Biggl(\frac{1 }{\sqrt{1-C^2}}\,\artanh \eta-
{\frac {1}{\sqrt {1-{C}^{2}}}{\artanh}   \,{\frac {\sqrt { 1-{C}^{2}}}{\sqrt {    1-C^2(1+\eta)/2}}} }
\\&\qquad+
{\artanh}{\frac {2-C^2 + 2\sqrt {1- \,{C}^{2}(1+\eta)/2} }{ 2-C^2\eta+ 2\sqrt {1- \,{C}^{2}(1+\eta)/2} }}
-\artanh \eta +\artanh C\Biggr)
\\&=2\Biggl(\frac{1 }{\sqrt{1-C^2}}\ln\frac{ C\sqrt{(1+\eta)/2}}{ \sqrt {   1-C^2(1+\eta)/2} + \sqrt{1-C^2}}
\\&\qquad
+\artanh    \sqrt {   1-C^2(1+\eta)/2}  +\artanh C\Biggr),
\end{align*}
\begin{commentx}
\begin{align*}
&=2\Biggl(\frac{1 }{\sqrt{1-C^2}}\artanh {\frac {{C}^{2}\eta-\sqrt {1-{C}^{2}}(1+\eta)2\sqrt {    1-C^2(1+\eta)/2}
 }{ 2+2\,\eta-{C}^{2}-2\,{C}^{2}\eta}}
\\&\qquad+\artanh \left( \,\sqrt {   1-C^2(1+\eta)/2}\right)
 +\artanh C\Biggr)
\end{align*}
\end{commentx}
 where we have used \eqref{eq:remo}, and we have transcribed $\frac1{\sqrt{1-C^2}} \artanh\ldots$
 into $\frac1{\sqrt{1-C^2}} \ln\ldots$ in the usual manner.
Then, taking the limit $\eta\nearrow1$, we obtain that
\begin{align*}
\Area_\SdS (     \widetilde B^C\setminus \widetilde E^C )
&=\lim_{\eta\nearrow1}\Area_\SdS(    (\widetilde B^C\setminus \widetilde E^C)\cap \{(x,y)\,:\, y\leq \eta\})
\\&=2\Biggl(\frac{1 }{\sqrt{1-C^2}}\ln\frac{C}{2\sqrt{1-C^2}}+\artanh \sqrt{1-C^2} +\artanh C\Biggr).
\end{align*}
Thus, ultimately, we reobtain \eqref{eq:lofi}, which is reassuring.
\snewpage

\snewpage
\section{Computation in the Beltrami--Poincaré half-plane model}\plabel{sec:bph}
Experience tells that some computations are easier in the Beltrami--Poincaré half-plane (BPh) model of the hyperbolic plane.
(In fact, they are sometimes laughably easy, which is somehow never quite the case with the BCK model.)
Therefore, it is reasonable to examine how things look like in the BPh model.

\textbf{The Beltrami--Poincaré half-plane model} will not be explained here in detail, there are plenty of expositions for that.
We just recall some technical matters for the sake of computations.
Considering proper hyperbolical, i.~e.~``interior'', points of the models,
 a point with coordinates $(x,y)$ in the BCK model transcribes to the point $(\check x,\check y)$
 in the BPh model, where
\begin{equation}
\check x=\frac{x}{1-y},\qquad\qquad\qquad
\check y=\frac{\sqrt{1-x^2-y^2}}{1-y}.
\plabel{eq:sat1}
\end{equation}
Conversely, the point $(\check x,\check y)$
 in the BPh model transcribes to the point $(x,y)$ in the BCK model, where
\begin{equation}
 x=\frac{2\check x}{ \check x^2+\check y^2+1},\qquad\qquad\qquad
 y=\frac{ \check x^2+\check y^2-1}{ \check x^2+\check y^2+1}.
 \plabel{eq:sat2}
\end{equation}
Then, in terms of the new coordinates, the area density is given as
\[\mathrm{d\check a}_{\hyp}  =\frac{|\mathrm d\check x\wedge \mathrm d\check y|}{\check y^2} ;\]
 and the arc length element (squared) is given by
\[\text{``}(\mathrm d\check s_{\hyp})^2\text{''}=\frac{ (\mathrm d\check x)^2+ (\mathrm d\check y)^2}{\check y^2}.\]

\textbf{The transcription of the canonical $h$-elliptic parabola and related objects.}
The description \eqref{eq:amat1} of the $h$-elliptical parabolical disk can be transcribed by substituting
\eqref{eq:sat2} into it.
It yields
$\dfrac{4((1-C^2) \check x^2-C^2\check y^2+C^2)}{(\check x^2+\check y^2+1)^2}\leq 0,$
 or, equivalently,
\[\frac{1-C^2}{C^2}\check x^2-\check y^2+1\leq0.\]
Considering that we deal only with $\check y>0$, this is the closed half interior of a hyperbola in Euclidean view.
Similarly, the half distance band \eqref{eq:amat2} transcribes to
\[ \left(\frac{\sqrt{1-C^2}}{C }\check x-\check y \right) \left(\frac{\sqrt{1-C^2}}{C }\check x+\check y \right)
\equiv\frac{1-C^2}{C^2}\check x^2-\check y^2 \leq0\quad\text{and}\quad\check x^2+\check y^2-1\geq0.\]
We can also note that the supporting horocycle \eqref{eq:shor} transcribes to
\[\check y^2-1=0,\]
 which is just
\[\check y=1\]
 because $\check y>0$ is required.
In this manner, we can transcribe $\partial E^C$, $\partial B^C$, $\partial D^C$ to the BPh model as
 $\partial\check E^C$, $\partial\check B^C$, $\partial\check D^C$.
Without going into further details, Figure \ref{fig:figHEP07} shows how the general configuration looks like in the BPh model.
We, however, mention that the focus $F^C=(0,\frac{C^2}{2-C^2})$ transcribes to
 $\check F^C=(0,\frac1{\sqrt{1-C^2}})$, which happens to be a focus of the of the respective hyperbola considered above.
In itself that would have been a weak reason for the designation of the focus, but it is nice to see now.
(And correspondingly, the horocycle directrix \eqref{eq:dhor} transcribes to $\check y=\sqrt{1-C^2}$, which is the directrix in Euclidean sense corresponding the
to the focus $\check F^C$ regarding the hyperbola.)

 \begin{figure}[ht]
   \begin{subfigure}[b]{2.5in}
    \includegraphics[width=2.5in]{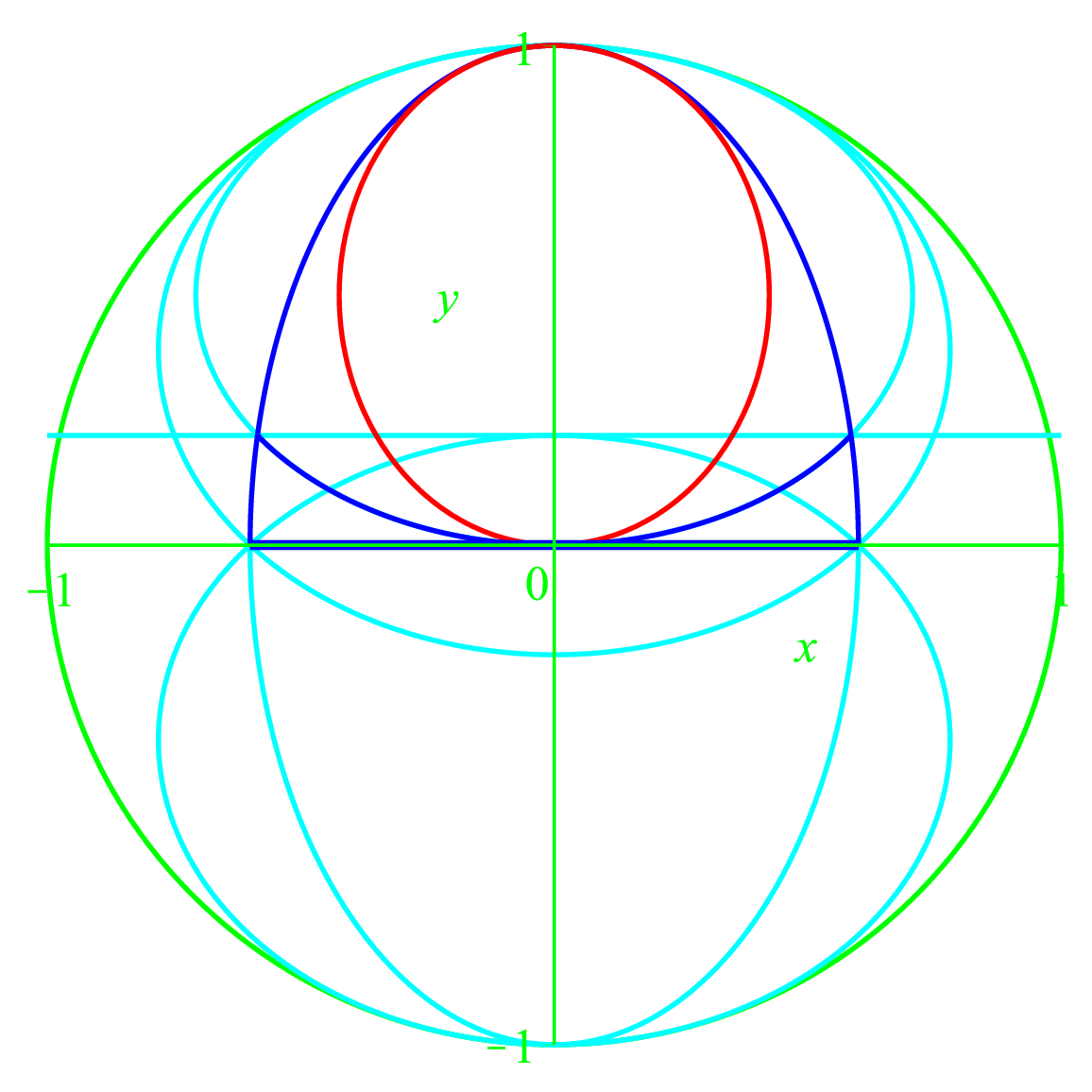}
     \caption*{\ref{fig:figHEP07}(a)  $\partial E^C $ [red], $\partial B^C $ [blue], $\partial D^C $ [blue].}
  \end{subfigure}
  \begin{subfigure}[b]{2.5in}
    \includegraphics[width=2.8in]{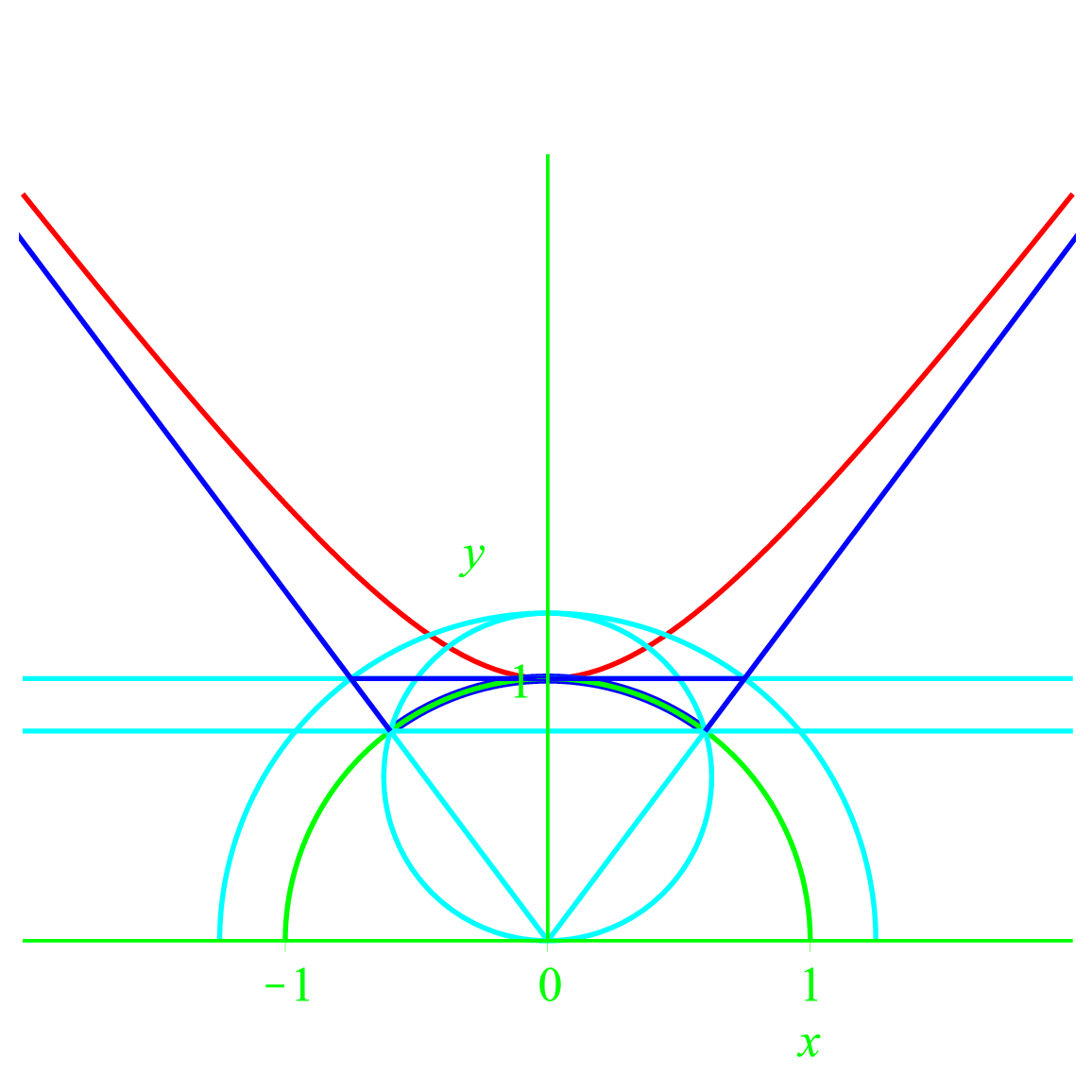}
   \caption*{\quad\ref{fig:figHEP07}(b)  $\partial\check E^C $ [red], $\partial\check B^C $ [blue], $\partial\check D^C $ [blue].}
  \end{subfigure}
  \\
   \caption*{Fig.~\ref{fig:figHEP07}  Transcription from the BCK model to the BPh model.
   }
\phantomcaption
\plabel{fig:figHEP07}

   \begin{subfigure}[b]{2.5in}
    \includegraphics[width=2.8in]{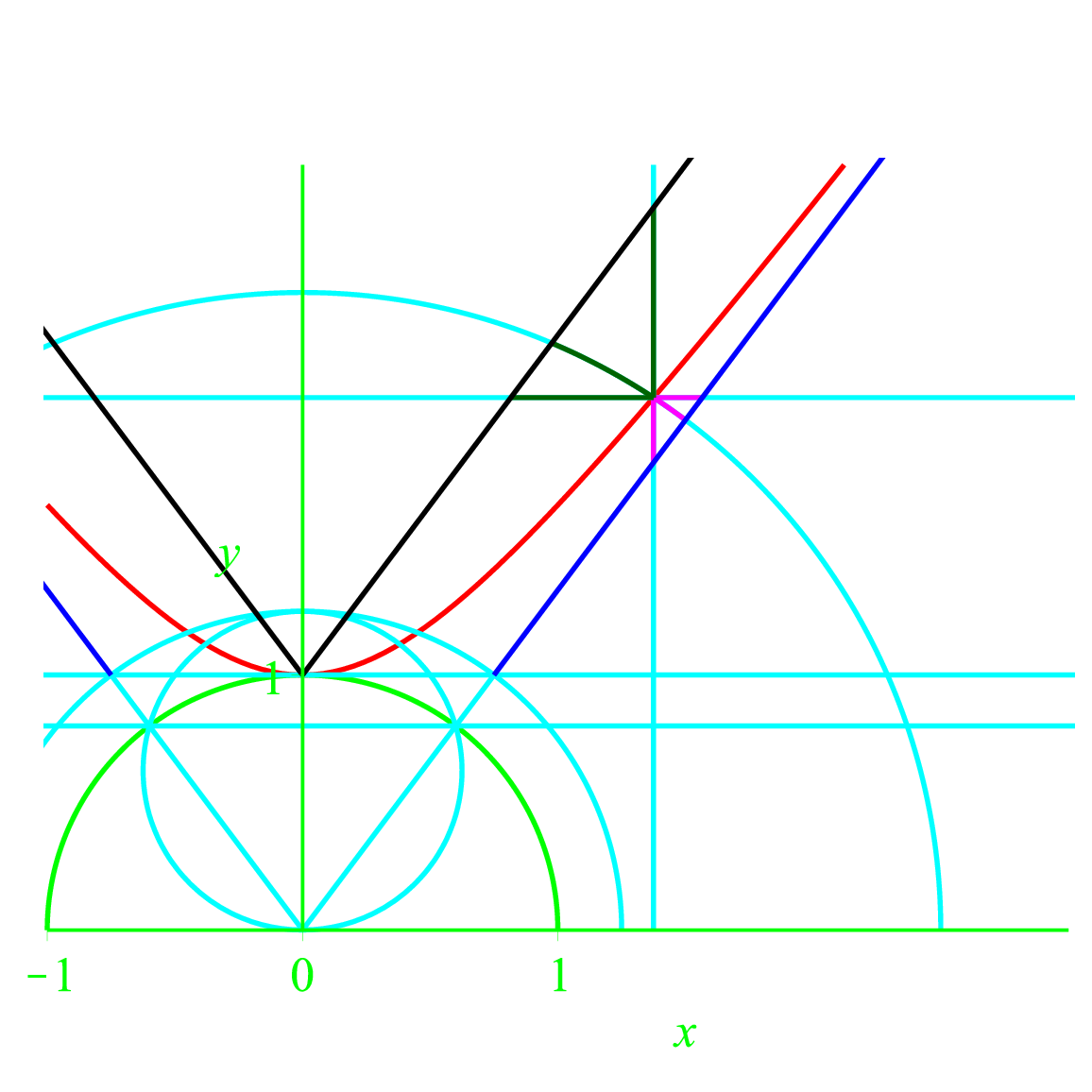}
     \caption*{\ref{fig:figHEP08}(a)  Variations in cutoff [magenta and\\ \phantom{aaaa} dark green],
     $\partial\check V^C $ [black]. }
  \end{subfigure}
  \begin{subfigure}[b]{2.5in}
    \includegraphics[width=2.5in]{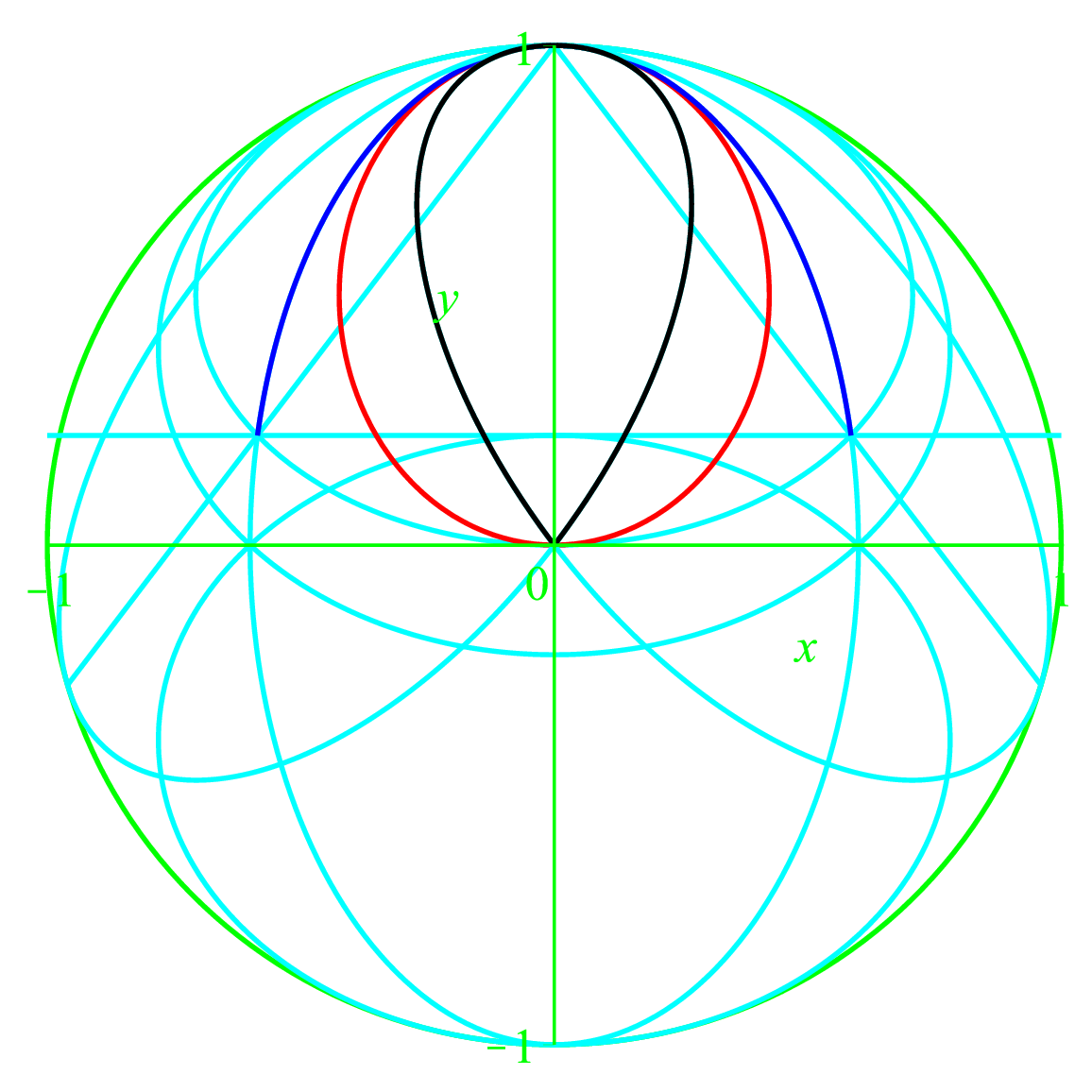}
   \caption*{\quad\ref{fig:figHEP08}(b)  $\partial V^C $ [black].\\ \phantom{$\partial\check V^C $}}
  \end{subfigure}
  \\
   \caption*{Fig.~\ref{fig:figHEP08}  Some variations.
   }
\phantomcaption
\plabel{fig:figHEP08}
\end{figure}

\snewpage

\textbf{The computation of the area difference in the BPh model.}
We concentrate only on the critical computation of  $\check\Area_\hyp(\check E^C\setminus\check D^C)$.
For $\vartheta\in(1+\infty)$, let
\[\check E^C_{  \vartheta}=\check E^C\cap\{(\check x,\check y)\,:\,\check y\leq \vartheta\}\]
 and
\[\check D^C_{  \vartheta}=\check D^C\cap\{(\check x,\check y)\,:\,\check y\leq \vartheta\}.\]
Note that we use horocyclical cutoffs here, thus, for example, $E^C_\eta$  does not  transcribe
 to $\check E^C_{\vartheta}$ with any $\vartheta$.
It does not, however, matter what cutoffs we use for an area computation.
Then
\begin{align*}
\check\Area_\hyp( \check D^C_\vartheta)&=
\int_{\check y=1}^{\vartheta}\int_{\check x=-\frac{ C}{\sqrt{1-C^2}} \check y }^{\frac{ C}{\sqrt{1-C^2}}\check y  }
\,\frac1{\check y^2}\, \mathrm d\check x\,\mathrm d\check y
\\&=\int_{\check y=1}^{\vartheta}\frac{ 2C}{\sqrt{1-C^2}}\,\frac{1}{\check y }\,\mathrm d\check y
\\&=\frac{ 2C}{\sqrt{1-C^2}} \ln \vartheta  .
\end{align*}
 and
\begin{align*}
\check\Area_\hyp( \check E^C_\vartheta)&=
\int_{\check y=1}^{\vartheta}\int_{\check x=-\frac{ C}{\sqrt{1-C^2}}\sqrt{\check y^2-1}}^{\frac{ C}{\sqrt{1-C^2}}\sqrt{\check y^2-1} }
\,\frac1{\check y^2}\, \mathrm d\check x\,\mathrm d\check y
\\&=\int_{\check y=1}^{\vartheta}
\frac{ 2C}{\sqrt{1-C^2}}
\,\frac{\sqrt{\check y^2-1}}{\check y^2}\,\mathrm d\check y
\\&=\frac{ 2C}{\sqrt{1-C^2}}\left(-\frac{\sqrt{\vartheta^2-1}}{\vartheta}+\artanh \frac{\sqrt{\vartheta^2-1}}{\vartheta}\right)
\\&\equiv\frac{ 2C}{\sqrt{1-C^2}}\left(-\frac{\sqrt{\vartheta^2-1}}{\vartheta}+\arcosh\vartheta\right)
\\&\equiv\frac{ 2C}{\sqrt{1-C^2}}\left(-\frac{\sqrt{\vartheta^2-1}}{\vartheta}+
\ln\left(\vartheta+\sqrt{\vartheta^2-1}\right)\right).
\end{align*}
Thus
\begin{align*}
\check\Area_\hyp( \check D^C_\vartheta\setminus \check E^C_\vartheta)&=
 \check\Area_\hyp( \check D^C_\vartheta)-\check\Area_\hyp(  \check E^C_\vartheta)
\\&=\frac{ 2C}{\sqrt{1-C^2}}\left( \frac{\sqrt{\vartheta^2-1}}{\vartheta} -
\ln\left(\vartheta+\sqrt{\vartheta^2-1}\right)+\ln\vartheta\right)
\\&=\frac{ 2C}{\sqrt{1-C^2}}\left(   \sqrt{1-\frac1{\vartheta^2}}-
\ln\left(1+ \sqrt{1-\frac1{\vartheta^2}} \right) \right).
\end{align*}
Therefore,
\[\check\Area_\hyp( \check D^C \setminus \check E^C )
=\lim_{\vartheta\nearrow+\infty}\check\Area_\hyp( \check D^C_\vartheta\setminus \check E^C_\vartheta)
=\frac{ 2C}{\sqrt{1-C^2}}\left(1 -\ln2 \right).\]
Note that the computation above is more accessible by standard calculus than the one in the BCK setting.
\snewpage

\textbf{The computation of the arc length difference in the BPh model.}
For $\vartheta\in(1+\infty)$, let
\[\partial^\cup\check E^C_{  \vartheta}=\partial\check E^C\cap\{(\check x,\check y)\,:\,\check y\leq \vartheta\}\]
 and
\[\partial^\cup\check D^C_{  \vartheta}=\partial\check D^C\cap\{(\check x,\check y)\,:\,\check y\leq \vartheta\}.\]
It is easy to see that the hypercyclic part of $\partial^\cup\check D^C_{  \vartheta}$ restricted to $\check x\geq0$ has arc length
\[\int_{\check y=1}^\vartheta\sqrt{1+
\left(\frac{\mathrm d}{\mathrm d\check y}\left(\frac{C}{\sqrt{1-C^2}}\check y\right)\right)^2}\,\frac{\mathrm d\check y}{\check y}
=\int_{\check y=1}^\vartheta \frac{1}{\sqrt{1-C^2}} \,\frac{\mathrm d\check y}{\check y}
=\frac{1}{\sqrt{1-C^2}}\log\vartheta,\]
 while the horocyclic part of $\partial^\cup\check D^C_{  \vartheta}$ restricted to $\check x\geq0$ has arc length
\[\int_{\check x=0}^{\frac{C}{\sqrt{1-C^2}}}1\,\frac{\mathrm d\check x}{ 1}=\frac{C}{\sqrt{1-C^2}}.\]
Therefore,
\[\check\Len_\hyp\left(\partial^\cup\check D^C_{  \vartheta}\right)=2 \frac{C}{\sqrt{1-C^2}}+2\frac{1}{\sqrt{1-C^2}}\log\vartheta.\]

On the other hand,
\begin{align*}
\check\Len_\hyp\left(\partial^\cup\check E^C_{  \vartheta}\right)&=
2\int_{\check y=1}^\vartheta\sqrt{1+
\left(\frac{\mathrm d}{\mathrm d\check y}\left(\frac{C}{\sqrt{1-C^2}}\sqrt{\check y^2-1}\right)\right)^2}\,\frac{\mathrm d\check y}{\check y}
\\&=2\int_{\check y=1}^\vartheta\sqrt{\frac{C^2+\check y^2-1}{(1-C^2)(\check y^2-1)}
}\,\frac{\mathrm d\check y}{\check y}
\\&=2\Biggl( \frac{1}{\sqrt{1-C^2}}\artanh\sqrt{\frac{ \vartheta^2-1 }{C^2+\vartheta^2-1}}
-\artanh\sqrt{\frac{(1-C^2)(\vartheta^2-1)}{C^2+\vartheta^2-1}}\Biggr).
\end{align*}
\begin{commentx}
\[=2\Biggl( \frac{1}{\sqrt{1-C^2}}\arsinh\sqrt{\frac{ \vartheta^2-1 }{C^2}}
-\arsinh\sqrt{\frac{(1-C^2)(\vartheta^2-1)}{C^2\vartheta^2}}\Biggr)\]
\[=2\Biggl( \frac{1}{\sqrt{1-C^2}}\arcosh\sqrt{\frac{C^2+\vartheta^2-1}{ C^2 }}
-\arcosh\sqrt{\frac{C^2+\vartheta^2-1}{C^2\vartheta^2}}\Biggr)\]
\[=2\Biggl(\frac1{\sqrt{1-C^2}}\ln\frac{\sqrt{C^2+\vartheta^2-1}+\sqrt{\vartheta^2-1} }{C}
-\ln\frac{\sqrt{1-C^2}\sqrt{\vartheta^2-1}+\sqrt{C^2+\vartheta^2-1} }{C\vartheta}\Biggr)\]
\end{commentx}
Rewriting $\frac1{\sqrt{1-C^2}}\artanh\ldots$ to $\frac1{\sqrt{1-C^2}}\ln\ldots$ as before, we find
\begin{multline*}
\check\Len_\hyp\left(\partial^\cup\check D^C_{  \vartheta}\right)-\check\Len_\hyp\left(\partial^\cup\check E^C_{  \vartheta}\right)=\\
=2\left(\frac1{\sqrt{1-C^2}}\left(C+\ln\frac{C\vartheta}{\sqrt{C^2+\vartheta^2-1}+\sqrt{\vartheta^2-1} }\right)  +\artanh\sqrt{\frac{(1-C^2)(\vartheta^2-1)}{C^2+\vartheta^2-1}} \right).
\end{multline*}

Therefore
\begin{align*}
\text{``$\check\Len_\hyp(\check D^C )-\check\Len_\hyp(\check E^C )$''}&
=\lim_{\vartheta\nearrow+\infty}\check\Len_\hyp\left(\partial^\cup\check D^C_{  \vartheta}\right)
-\check\Len_\hyp\left(\partial^\cup\check E^C_{  \vartheta}\right)
\\&=2\left(\frac1{\sqrt{1-C^2}}\left(C+\ln\frac{C }{ 2}\right)  +\artanh\sqrt{ 1-C^2  } \right).
\end{align*}
This yields the same as \eqref{eq:Gprime}.
(Here the computation was not as simple as in the case of the area,
 but it was a bit simpler than in the BCK setting.)

Note, however, that we used  horocyclyc cutoffs, and not lineal ones, therefore we have not reproduced the computation of $G'^C$.
On the other hand, the BPh model is particularly advantageous to address this discrepancy.
One can see that if we take some particular cutoffs though the endpoint of $\partial^\cup\check E^C_{  \vartheta}$
 as in Figure \ref{fig:figHEP08}(a) (the magenta arcs), the resulting cutoff points  on the corresponding
 hypercyclic arcs are not far from each other as $\vartheta\nearrow+\infty$; in fact, their distance from each other,
 and the lenghths of the cutoff-arcs go to $0$ (considering the distortion $\frac1{\check y}$ relative to the Euclidean case).
This applies at least to horocyclic cutoffs through the asymptotic point of the $h$-elliptic parabola,
 lineal cutoffs perpendicular to the  $h$-elliptic parabola,
 and lineal cutoffs through the asymptotic point of the $h$-elliptic parabola.
%Indeed, distances on the hypercyclic are very easy to estimate  on the hypercyclic arc.
We do not carry out the estimates explicitly, but the reader can now easily see that all of the above mentioned
 cutoff types lead to the same arc length difference limit.

In fact, the observations above are valid even if the hypercyclic arc is subjected to horocyclic translation
 through the asymptotic point of the $h$-elliptic parabola, cf.~Figure \ref{fig:figHEP08}(a) (the dark green arcs).
This leads to the idea of
\snewpage

\textbf{An inner approximation.}
We can take the hypercyclic arcs of $\partial D^C$ ($\partial\check D^C$) and
 apply horocylic translations through the   asymptotic point of the $h$-elliptic parabola,
 such that the endpoints of hypercyclic arcs get translated to the vertex of the $h$-elliptic parabola.
Then we can take convex closure yielding  $  V^C$ ($ \check V^C$).
It is unwieldy in the BCK model,
%\begin{multline*}
\begin{align*}
V^C=\bigl\{(x,y)\,:&\,{x}^{2}+2\,C(\sqrt {1-{C}^{2}}\,x+Cy)(y-1)\leq0,\\
&\,{x}^{2}+2\,C(-\sqrt {1-{C}^{2}}\,x+Cy)(y-1)\leq0
\bigr\},
\end{align*}
%\end{multline*}
 cf.~Figure \ref{fig:figHEP08}(b); but it is easy to describe in the BPh model,
\[\check V^C=\left\{(\check x,\check y)\,:\, \check y\geq 1+\frac{\sqrt{1-C^2}}{C}|\check x|\right\},\]
 cf.~Figure \ref{fig:figHEP08}(a).
In this latter setting, we can compute easily.

In any way, the general pattern is that
\[B^C\supset D^C\supset E^C\supset V^C,\]
 and similarly in the BPh model.

Regarding the area, it is easy to see that
\[
%\begin{align*}
\check\Area_\hyp( \check D^C \setminus \check V^C)
%&
=\int_{\check y=1}^{\infty}\frac{ 2C}{\sqrt{1-C^2}}\,\frac{1}{\check y^2 }\,\mathrm d\check y
%\\&
=  \left[-\frac{ 2C}{\sqrt{1-C^2}}\,\frac{1}{\check y }\right] _{\check y=1}^{\infty}=\frac{ 2C}{\sqrt{1-C^2}}.
%\end{align*}
\]
Therefore
\[\Area_\hyp( \check E^C \setminus \check V^C)=\frac{ 2C}{\sqrt{1-C^2}}\ln 2.\]
One can also see that
\[\check E^C=\check V^C_{\mathrm{up}\,-\ln2}\quad\text{in area}.\]

In terms of circumference, we can easily see, using the horocyclic cutoff, that
\[\text{``$\check\Len_\hyp(\check D^C )-\check\Len_\hyp(\check V^C )$''}=\frac{ 2C}{\sqrt{1-C^2}}  \]
 (we get difference in the horocyclic arc).
Therefore
\begin{align*}
\widehat G^C\equiv
\text{``$\check\Len_\hyp(\check E^C )-\check\Len_\hyp(\check V^C )$''}
=2\left(\frac1{\sqrt{1-C^2}}\left(\ln\frac{ 2}{C }\right)  -\artanh\sqrt{ 1-C^2  } \right).
\end{align*}
Here $\widehat G^C$ is monotone increasing in $C\in(0,1)$;  $\lim_{C\searrow0}\widehat G^C=0$ and $\lim_{C\nearrow1}\widehat G^C=+\infty$.
Taking
\[\widehat \beta(C)=\ln\frac{ 2}{C }  -\sqrt{1-C^2}\artanh\sqrt{ 1-C^2  }, \]
 one can see that
\[\text{``$\check\Len_\hyp(\check E^C )-\check\Len_\hyp(\check V^C_{\mathrm{up}\,-\widehat \beta(C)} )$''}=0\]
is valid using the horocyclic and standard lineal cutoff.
Here $\widehat \beta(C)$ is monotone increasing in $C\in(0,1)$;
 $\lim_{C\searrow0}\widehat \beta(C)=0$ and $\lim_{C\nearrow1}\widehat \beta(C)=\ln 2$.
\snewpage

\section{Reflections}\plabel{sec:ww}
\textbf{What we did ``wrong''} is easy to point out:
We have done a lot of unnecessary, complicated computations in the BCK model.
More precisely:
In approaching the area and circumference problem of $E^C$, we did have ``essential'' and ``inessential'' computations.
The ``essential'' ones are  related to $E^C_\eta$ directly, and the ``inessential'' ones are related to
 areas and circumferences of hypercyclic and horocyclic segments, etc.
In those latter ones we could have referred to common knowledge known already to Lobachevsky and Bolyai;
 or if one is dissatisfied with the mathematical rigor of those presentations, one can refer to various
 introductory works.
Therefore the results of inessential computations could have been obtained from the
 combination very classical knowledge and the basic distance formulae of the models.
The case of the essential computations is different, but even there
 the area computation is much each easier to be engineered back from the BPh model.
This leads to the matter of the intellectual honesty of the (essential) computations:
%While the analytic computations presented were relatively straightforward (as checking
% primitive functions is not hard), this simplicity is somewhat deceiving.
While primitive functions are relatively easy to check, getting them is not so easy.
They are accessible by classical methods and  computer algebra software, but obtaining them can be messy.
Ultimately, they are quite unrealistic to be claimed to be done by hand, or using naive intuition.
(To get a measure, the reader can try to compute
$\Area_\SdS(\nnn\widetilde{W}^C_\eta)$, $\Area_\SdS(\nnn\widetilde{E}^C_\eta)$, $\Area_\SdS(\nnn\widetilde{V}^C_\eta)$
 interchanging the order of integration in terms of $x$ and $y$.
Admittedly, the computation is more complicated, but not immensely so, thus one can try.)
A relatively minor issue but which affects the actual face of the computations is as follows:
Primarily, integrals were given mainly in terms of `$\artanh$'.
Those expressions could have been written in terms  `$\ln$'; which would be more advantageous in terms of addition formulae,
 and asymptotic estimations.
Another possibility is to write those expressions in terms of `$\arsinh$', where
 the resulted expressions are often simpler (this can be another fun assignment for the reader),
 and we do not have to check the domain of the function as we should always do in case of
 `$\artanh$' (although that is mostly easy to do).
In order to put it differently:
Instead of uniformity we could have been more opportunistic for the sake of nice expressions.
Ultimately, it can be said that our presentation involving the computations in the BCK model
 do not satisfy the efficiency requirements of mathematical presentations.

\snewpage

\textbf{What we did fine.}
Despite the criticisms offered in the previous paragraph,
 I think we did well by covering the situation in the BCK model extensively.
And the reason is the following:
Understanding follows different ways than the efficiency principles of mathematical communications.
The BCK model is natural for the study of hyperbolic conics,
 thus it is useful to get a grip on it in an accessible situation,
 even if certain computations are simpler in another model.
In order to understand the proportions of the argument, it is useful to see
 even the elements of the non-essential computations.
Constantly applying to references may confuse the depth of the argument,
 and it is annoying to pretend that one should remember the content the classics readily from the start.
The BCK model is theoretically superior in terms of considering the collineation and convexity properties,
 therefore it is hard avoid mentioning it entirely, anyway.
Moreover, the polar computation of the circumference requires the context of the ambient projective space.
Using `$\artanh$' a lot is a reasonable choice:
Hyperbolic trigonometric functions are more in the analogy of trigonometric functions and
 spherical geometry; yet `$\artanh$' is arithmetically closer to $\ln$ than `$\arsinh$',
 and it is more in the spirit of the BCK model.
On the other hand, we also did well by stepping out of the BCK model:
One should not be ignorant about alternative approaches, especially if they are more effective.
Hyperbolic geometry illustrates particulary well that  it is very useful to be opportunistic in mathematics;
 changes in the viewpoint and adopting synthetic ideas can expedite computations very much.
While Lobachevsky's hope that considering multiple viewpoints in hyperbolic geometry will
 lead to great advances in calculus is now considered too naive;
 at least for computations related to hyperbolic geometry, using multiple viewpoints, synthetic or not, may indeed enhance
 the effectiveness of the computations, or rather can help to avoid certain pitfalls in less manageable cases.
Moreover, ultimately, and that also applies to research mathematics,
 it is more realistic to have a somewhat wider knowledge to select from relative to what can be presented in brief terms.

\textbf{What we did learn mathematically.}
Although no explicit comparisons were made, what we have learned underscores
 that there are three tiers of non-degenerate conics from analytic viewpoint:
In the first tier we have the uniform conics (the circle in the Euclidean case, the cycles in the hyperbolic case),
 these are the nicest ones; in the second tier we have some strange but analytically treatable ones
 (like the parabola in the Euclidean case, and the $h$-elliptic parabola in the hyperbolic case);
 and, finally, in the third tier, we have the analytically more challenging ones
 (like the proper ellipses in the Euclidean and hyperbolic cases).
\snewpage

\textbf{What we did in general.}
We have hopefully learned certain things about some hyperbolic conics.
(Hopefully not too little, hopefully not too much.)
But the main point is different.
Mathematics has technically grown far away from popular textbooks.
This is certainly good, as it is, say, in the case of consumer electronics.
Nevertheless, this separation can be formidable.
However, there is a class of mathematical statements which, for the technically skilled,
 can be termed as ``observations'', that are not entirely trivial ones but can be explored well with a given technical apparatus.
Now, I think that exposing technical information can be of value, and, if possible, should be pursued.
In that, however, neither imitating elementary mathematics nor playing for unexplainable elegance
 is of value, but viewing things from multiple angles is.
Experience tells that mathematical understanding is more likely to come from a diluted exposition, where
 this latter expression does not mean lack of content but giving the opportunity of reflection
 for the reader about the advantages and disadvantages, possibilities and limits of an approach.
We live in an electronic age, but mathematical expositions often follow old instincts,
 like emphasising classical problems and/or the elegance of a particular approach.
In the first case one is not necessarily convinced about the relative importance of classical questions,
 and in the second case one cannot necessarily appreciate the beauty of an approach in lack of comparison.
(Therefore the most elegant exposition can be boring or weightless in itself.)
Nor very elegant approaches to pure theory prepare for applications.
Research articles can be marvels of treating a problem, but do  not necessarily expose
 discomfitures and alternative approaches.
Also, while effectiveness is to be praised in mathematics, space constraints often limit the scope of investigations.
Therefore, I believe that discussions of technical but less than challenging problems
 exposed from multiples angles, yet leaving something to aesthetical judgment of the reader, may be of interest.

\snewpage

\end{document}